\documentclass[19pt]{article}
\usepackage[utf8]{inputenc}
\usepackage[english]{babel}
\usepackage{indentfirst}
\usepackage{misccorr}
\usepackage{color}
 \usepackage{amsmath}
\usepackage{amssymb}

\def\R{{\mathbb R}}
\def\aff{\mbox{aff}}

\newtheorem{lemma}{Lemma}

\newtheorem{definition}{Definition}

\newtheorem{theorem}{Theorem}

\newtheorem{remark}{Remark}

\newtheorem{proposition}{Proposition}

\newtheorem{corollary}{Corollary}

\begin{document}

\title{The Conic Tverberg plus minus}

\author{E. Kolpakov}

\date{}

\maketitle

\tableofcontents

\begin{abstract} We prove a Tverberg type theorems: Given a set $A\subset \R^d$ and an affine subspaces $L_1,...,L_m\subset \R^d$. In some conditions there is a partition $A_1, ..., A_m, B_1, ..., B_m$ of $A$ such that the sets $L_i\cup A_i\cup B_i$ for $i=1,...,m$ have common barycentre with negative masses in $B_i$ and positive masses in $A_i$. The case $B_1=...=B_m=L_1=...=L_m=\varnothing$ is Tverberg's classical theorem.
The case $L_1=...=L_m=\varnothing$ gives a partial answer to the question of Barany-Soberon about the possible number of negative masses.
\end{abstract}

%\keywords{Tverberg's theorem, Roudneff's theorem, positive cone, sign conditions}
%\subjclass[2000]{Primary 52A37, secondary 52B05}

\section{Introduction and result}

Denote the convex hull of a set $ P $ by $ \left< P \right> $.

\begin{theorem}[Tverberg]\label{Tverbergs_theorem} Let $X$ be a  set of $(m-1)(d+1)+1$ points in $\R^d$.
Then there is a partition $X_1, X_2, ... ,X_m$ of $X$
such that $\bigcap\limits_{i=1}^m  \left< X_i\right> \neq\varnothing$.\end{theorem}

Theorem \ref{Tverbergs_theorem} was given in [Tv66].
Consider some generalizations of this Theorem.

For an affine subspace $L\subset\R^d$ and a finite subset $X\subset\R^d$ denote by $[L, X)$ the set of barycentres $\sum\limits_{x\in L\cup X}\alpha_xx$
with $\sum\limits_{x\in L\cup X}\alpha_x=1$ and   $\alpha_x\geq0$ for all $x\in X$.
Without conflict with this notation we denote by $[0;+\infty)$ the set of non-negative real number.

The following definition was given in [Ro01]:
\begin{definition} A collection of affine subspaces of $\R^d$ is in general
position at infinity if the hyperplane at infinity does not contain any vertex of the projective
arrangement defined by those subspaces. In other words, each subcollection of these affine
subspaces which crosses at infinity also crosses in $\R^d$.\end{definition}

Note that the convex hull of a set $ P $ is equal to $[\varnothing, P)$.
Then the following theorem is a generalization of Theorem \ref{Tverbergs_theorem} :

\begin{theorem}[Roudneff]\label{Roudneffs_theorem} Let $X$ be a finite set of points in $\R^d$ and $L_1, L_2, ... ,L_m$ be a collection of affine subspaces of $\R^d$ in general position at infinity.
If\newline $|X|+\sum\limits_{i=1}^{m}(1+\dim L_i)\geq(m-1)(d+1)+1$, then there is a partition $X_1, X_2, ... ,X_m$ of $X$
such that $\bigcap\limits_{i=1}^m[L_i, X_i)\neq\varnothing$.\end{theorem}

 Theorem \ref{Roudneffs_theorem} was given in [Ro01].

 \begin{definition}
  A partition $\{A_1,\ldots,A_r\}$ of subset $A\subset\R^d$ is called {\sl proper} if $1\le |A_j|\le d+1$ for every $j \in [r]$. Here and in what follows $[r]$ stands for the set $\{1,\ldots,r\}$.
 \end{definition}

 The Theorem \ref{Tverbergs_theorem} has a form which was considered in [BS16]:

 \begin{theorem}\label{th:tverb}
 Let $X$ be a  set of $(m-1)(d+1)+1$ points in $\R^d$.
 Suppose that the coordinates of the points of the set $X$ (altogether $d((m-1)(d+1)+1)$ real numbers) are algebraically independent.
 Then there is a proper partition of $X$ into sets $X_1,\ldots,X_r$ such that the solution $\{\alpha(x)\}_{x\in X}$ of the following system:

 \begin{equation}\label{aff}
z=\sum_{x\in X_j}\alpha(x)x \quad \mbox{ and } \quad 1=\sum_{x\in X_j}\alpha(x) \quad \mbox{ for all }j \in \{1,...,r\}.
\end{equation}
do not contains negative numbers.
\end{theorem}

The main result of [BS16] is a Theorem~\ref{th:pm}:

\begin{theorem}\label{th:pm}
Assume $k \in \{0,1,\ldots,r-1\}$. Under the conditions of Theorem~\ref{th:tverb} there is a (proper) partition of $X$ into $r$ parts $X_1, ..., X_r$ so that in the unique solution to (\ref{aff}) $\alpha(x) <0$ for exactly $k$ elements $x \in X$.
\end{theorem}

This Theorem was proved in [BS16] in a stronger form:

\begin{theorem}\label{th:tverbM}
Under the conditions of Theorem~\ref{th:tverb} let $M\subset A$ be a set of size at most $r-1$ such that $\left< M\right> \cap \left< A\setminus M\right>=\varnothing$.
Then there is a partition of $X$ into $r$ parts $X_1, ..., X_r$ such that in (\ref{aff}) $\alpha(x)<0$ if and only if $x \in M$.
\end{theorem}
%As generalization of

For an affine subspace $L\subset\R^d$ and a finite subsets $A, B\subset\R^d$ such that $L\cup A\neq\varnothing$
denote by $C(L, A, B)$ the set of barycentres $\sum\limits_{x\in L\cup A\cup B}\alpha_xx$
with

$\bullet$ $\sum\limits_{x\in L\cup A\cup B}\alpha_x=1$;

$\bullet$ $\alpha_x\geq0$ for all $x\in A$;

$\bullet$ $\alpha_x\leq0$ for all $x\in B$.

\begin{definition}
  A subset $A\subset \R^d$ is divisible for a collection $L_1, L_2, ... ,L_m$ of affine subspaces of $\R^d$ with $k$ negative masses if there is a partition $A_1, A_2, ... ,A_m, B_1, B_2, ... ,B_m$ of $A$
such that

 $\bullet$ $\sum\limits_{i=1}^m|B_i|=k$;

  $\bullet$ $\bigcap\limits_{i=1}^mC(L_i, A_i, B_i)\neq\varnothing$.
\end{definition}

The main results of the present article are Theorems~\ref{main_result} and~\ref{main_result_2} :

\begin{theorem}\label{main_result}
 Let $A$ be a finite set of points in $\R^d$ and $L_1, L_2, ... ,L_m$ be a collection of affine subspaces of $\R^d$ in general position at infinity. Suppose that

 $\bullet$ $|A|+\sum\limits_{i=1}^{m}(1+\dim L_i)\geq (m-1)(d+1)+1$;

%$\bullet$ $k\leq\max((m-1)(d+1)-\sum\limits_{i=1}^m(1+\dim L_i);0)$.

$\bullet$ $|A|\geq k+|\{i|L_i=\varnothing\}|$.

 Then $A$ is divisible for a collection $L_1, L_2, ... ,L_m$ of affine subspaces of $\R^d$ with $k$ negative masses if at least one of the conditions holds:

 1) $k\leq \min\limits_{\mathfrak{p}\subset\{1,...,m\}, |\mathfrak{p}|\leq \max(0;m-d-1)}\sum\limits_{i\in \mathfrak{p}}\lfloor \frac{d-\dim L_i}{2}\rfloor$;

 2) $k+|\{i|L_i\neq\varnothing\}|\leq m$ and $d=2$;

 3) $k\leq 1+\min\limits_{\mathfrak{p}\subset\{1,...,m\}, |\mathfrak{p}|\leq \max(0;m-d-1)}\sum\limits_{i\in \mathfrak{p}}\lfloor \frac{d-\dim L_i}{2}\rfloor$ and there is a number $i_{tv}$ such that  $L_{i_{tv}}=\varnothing$.

\end{theorem}

\begin{remark}
 Under the conditions of Theorem~\ref{Roudneffs_theorem} take $A=X$, $k=0$. %Then $\left<A\right>\cap\left<B\right>=\varnothing$; $|B|=0$.
 Then $$|A|+\sum\limits_{i=1}^{m}(1+\dim L_i)\geq (m-1)(d+1)+1$$
 And $k\leq1$.
 Apply Theorem~\ref{main_result}, we obtain a partition $A_1, A_2, ... ,A_m$ of $A$ such that $\bigcap\limits_{i=1}^mC(L_i, A_i, \varnothing)\neq\varnothing$.
 Therefore Theorem~\ref{main_result} is a generalization of Theorem~\ref{Roudneffs_theorem}.
\end{remark}

Consider a some generalization of algebraic independence:

\begin{definition}\label{strong_independence}
 A set $S$ in $\R^d$ is said to be {\sl strongly independent} provided that each finite family $\{S_1, ..., S_r\}$ of pair-wise disjoint subsets of $S$ has the following property:

 $\bullet$ If $|S_i|\leq d+1$ for each $i$, then $$\dim(\bigcap_{i=1}^{r}\aff S_i)=\max(-1;-1-(d+1)(r-1)+\sum\limits_{i=1}^r|S_i|).$$
\end{definition}

In [Re68] it is mentioned that the algebraic independence of a set implies the strong independence.
 Let us generalize the definition~\ref{strong_independence} with the following definition.

\begin{definition}\label{strong_independence_with_flats}
 A set $S$ in $\R^d$ is said to be {\sl strongly independent with affine subspaces $L_1, ..., L_m$} provided that each finite family $\{S_1, ..., S_r\}$ of pair-wise disjoint subsets of $S$ and each subcollection $\{Q_1, ..., Q_r\}\subset\{L_1, ..., L_m\}$ have the following property:

 $\bullet$ If $|S_i|\leq d-\dim Q_i$ for each $i$, then $$\dim(\bigcap_{i=1}^{r}\aff (Q_i\cup S_i))=\max(-1;-d(r-1)+\sum\limits_{i=1}^r|S_i|+\sum\limits_{i=1}^r\dim Q_i).$$
\end{definition}

Clearly, the case $L_1=...=L_m=\varnothing$ in Definition~\ref{strong_independence_with_flats} is Definition~\ref{strong_independence}.

 \begin{theorem}\label{main_result_2}
  Take notations $A, L_1, ..., L_m, d$ from Theorem~\ref{main_result}.
  Suppose that

 $\bullet$ $|A|+\sum\limits_{i=1}^{m}(1+\dim L_i)= (m-1)(d+1)+1$;

$\bullet$ $|A|\geq k+|\{i|L_i=\varnothing\}|$.

Additionaly suppose that the set $A$  is  strongly independent with affine subspaces $L_1, ..., L_m$.
Then $A$ is divisible for a collection $L_1, L_2, ... ,L_m$ of affine subspaces of $\R^d$ with $k$ negative masses if at least one of the conditions holds:

 1) $k+|\{i|L_i\neq\varnothing\}|\leq m-1$;

 2) $k\leq \lfloor \frac {|A|}2\rfloor$ and $m=2$;

 3) $|A|$ is even and $k=\frac{1}{2}|A|$.
 \end{theorem}

\begin{lemma}\label{intersection_in_general_position}
 Let $A_1, ..., A_m, B_1, ..., B_m$ be a finite subsets of $\R^d$ and affine subspaces $L_1, ..., L_m$ such that:

 $\bullet$ the set $(\bigcup\limits_{i=1}^mA_i)\cup(\bigcup\limits_{i=1}^mB_i)$ is strongly independent with $L_1, ..., L_m$;

 $\bullet$ the intersection $\bigcap\limits_{i=1}^mC(L_i, A_i, B_i)$ is an unique point $\omega$.

 %Moreover, if the set $\{A\cup B\}$ has strong independence with subspaces $L_1, L_2, ... ,L_m$ then %the intersection $\bigcap\limits_{i=1}^mC(L_i, A_i, B_i)$ is unique point.
Then there is a family of numerical sets $\{\{\alpha_{x,i}\}_{x\in L_i\cup A_i\cup B_i}\}_{i=1}^m$
such that for any $i=1,...,m$:

$\bullet$ $\sum\limits_{x\in L_i\cup A_i\cup B_i}\alpha_{x,i}x=\omega$;

$\bullet$ $\sum\limits_{x\in L_i\cup A_i\cup B_i}\alpha_{x,i}=1$;

$\bullet$ $\alpha_{x,i}>0$ for $x\in A_i$;

$\bullet$ $\alpha_{x,i}<0$ for $x\in B_i$.
\end{lemma}

{\bf Proof.}
 For any $i=1, ..., m$ we have $\omega\in C(L_i, A_i, B_i)$.
 According to the definition of $C(L_i, A_i, B_i)$  there is a  numerical set $\{\alpha_{x,i}\}_{x\in L_i\cup A_i\cup B_i}$
such that:

$\bullet$ $\sum\limits_{x\in L_i\cup A_i\cup B_i}\alpha_{x,i}x=\omega$;

$\bullet$ $\sum\limits_{x\in L_i\cup A_i\cup B_i}\alpha_{x,i}=1$;

$\bullet$ $\alpha_{x,i}\geq 0$ for $x\in A_i$;

$\bullet$ $\alpha_{x,i}\leq 0$ for $x\in B_i$.

For any $i=1, ..., m$ denote by $A_i', B_i'$ a subsets of $A_i, B_i$ respectively such that $\alpha_{x,i}\neq 0$ for $x\in A_i'\cup B_i'$ .
Hence the point $\omega$ belong to relative interior of  $ C(L_i, A_i', B_i')$.
%Then \newline $\omega\in\aff((L_i\cup A_i\cup B_i)\setminus\{x\})$.
Therefore $$\{\omega\}=\bigcap\limits_{i=1}^m\aff(L_i\cup A_i'\cup B_i').$$
From strongly independence follows that:
$$d(m-1)=\sum\limits_{i=1}^m\dim\aff(L_i\cup A_i'\cup B_i').$$
Suppose that $\alpha_{y,j}= 0$ for some $j$ and $y\in A_j\cup B_j$.
Than
$$\dim\aff(L_j\cup A_j'\cup B_j'\cup\{y\})+\sum\limits_{i\neq j}^m\dim\aff(L_i\cup A_i'\cup B_i')=1+d(m-1).$$
Hence the intersection
$\aff(L_j\cup A_j'\cup B_j'\cup\{y\})\cap\bigcap_{i\neq j}^m\dim\aff(L_i\cup A_i'\cup B_i')$
js 1-dimensional.
Since the point $\omega$ belong to relative interior of  $ C(L_i, A_i', B_i')$ for any $i$, the both intersections $$C(L_j, A_j', B_j'\cup\{y\})\cap\bigcap_{i\neq j}^m\dim\aff(L_j\cup A_j'\cup B_j')$$
and $$C(L_j, A_j'\cup\{y\}, B_j')\cap\bigcap_{i\neq j}^m\dim\aff(L_j\cup A_j'\cup B_j')$$
contains a points besides $\omega$.
Therefore the intersection $\bigcap\limits_{i=1}^mC(L_i, A_i, B_i)$ contains a points besides $\omega$.
A contradiction.
Then $\alpha_{x,i}\neq 0$ for any $i$ and $x\in A_i\cup B_i$.

\begin{remark}\label{implications_of_Barany_Soberon_results}
 Under the conditions of Theorem~\ref{main_result} take $A$ as set of $(m-1)(d+1)+1$ algebraically independent point,  $L_1=...=L_m=\varnothing$.
 Hence $|\{i|L_i\neq\varnothing\}|=0$.
 Apply either Theorem~\ref{main_result} or~\ref{main_result_2}, we obtain a partition $A_1, A_2, ... ,A_m, B_1, B_2, ... ,B_m$ of $A$ such that  $\sum\limits_{i=1}^m|B_i|=k$ and
   $\bigcap\limits_{i=1}^mC(\varnothing, A_i, B_i)\neq\varnothing$.
 The set $A$ is strongly independent because it is algebraically independent.
 Then set $A$ is strongly independent with subspaces $L_1=...=L_m=\varnothing$.
   Hence the intersection $$\bigcap\limits_{i=1}^m\aff(A_i\cup B_i)$$ is the unique point.
   Denote it by $\omega$.
   Then $\{\omega\}=\bigcap\limits_{i=1}^mC(\varnothing, A_i, B_i)$ because $$\bigcap\limits_{i=1}^mC(\varnothing, A_i, B_i)\subset\bigcap\limits_{i=1}^m\aff(A_i\cup B_i).$$
 Apply Lemma~\ref{intersection_in_general_position} for the  point $\omega$, we obtain the a family $\{\{\alpha_{x,i}\}_{x\in  A_i\cup B_i}\}_{i=1}^m$.
 Take:

 $\bullet$ $r=m$;

  $\bullet$ $X=A$;

  $\bullet$ $X_1=A_1\cup B_1, ..., X_r=A_r\cup B_r$ - partition of $X$;

  $\bullet$ $\alpha(x)=\alpha_{x, i}$ for any $x\in X$ and $i$ such that $x\in X_i$.

Then $\alpha(x) <0$ for exactly $k$ elements $x \in X$.

Therefore Theorem~\ref{main_result} in conditions 2,3 implies Corollary~\ref{plane_intersections},  Corollary~\ref{maximum_of_negative_points} respectively.
Theorem~\ref{main_result_2} in conditions 1,2,3 implies Theorem~\ref{th:pm}, Corollary~\ref{Radon_generalization}, Corollary~\ref{half_is_negative} respectively.
\end{remark}

The Corollary~\ref{plane_intersections} for $r=3$ mentioned in [BS16].

\begin{corollary}\label{plane_intersections}
 Assume $k \in \{0,1,\ldots,r\}$. Under the conditions of Theorem~\ref{th:tverb} with $d=2$ there is a (proper) partition of $X$ into $r$ parts $X_1, ..., X_r$ so that in the unique solution to (\ref{aff}) $\alpha(x) <0$ for exactly $k$ elements $x \in X$.
\end{corollary}

According to both corollaries~\ref{plane_intersections} and~\ref{maximum_of_negative_points}  we can require exactly $r$ of the $\alpha(x)$ are negative in Theorem~\ref{aff} for $d\neq 1$, $r\geq 2d+1$.

\begin{corollary}\label{maximum_of_negative_points}
 Assume $k \in \{0,1,\ldots,1+\lfloor \frac {d+1}2\rfloor\max(0;r-d-1)\}$. Under the conditions of Theorem~\ref{th:tverb} there is a (proper) partition of $X$ into $r$ parts $X_1, ..., X_r$ so that in the unique solution to (\ref{aff}) $\alpha(x) <0$ for exactly $k$ elements $x \in X$.
\end{corollary}

Both corollaries~\ref{Radon_generalization} and~\ref{half_is_negative} proved in [BS16]:

\begin{corollary}\label{Radon_generalization}
 Assume $k \in \{0,1,\ldots,\lfloor \frac {d+2}2\rfloor\}$. Under the conditions of Theorem~\ref{th:tverb} with $r=2$ there is a (proper) partition of $X$ into $r$ parts $X_1, ..., X_r$ so that in the unique solution to (\ref{aff}) $\alpha(x) <0$ for exactly $k$ elements $x \in X$.
\end{corollary}

\begin{corollary}\label{half_is_negative}
 Assume $r,d$ are both even and positive integers, $A$ is a set of algebraically independent point, $|A|=(r-1)(d+1)+1$, and $k=\frac 12 [(r-1)(d+1)+1]$. Then $A$ has a proper partition into $r$ parts $X_1, ..., X_r$ such that in (\ref{aff}) exactly $k$ of the coefficients $\alpha(x)$ are negative.
\end{corollary}

According to Theorem~\ref{main_result_2} both corollaries~\ref{Radon_generalization} and~\ref{half_is_negative} holds for a set $A$ in strong general position.

 \section{Sufficiently general position}

\begin{definition}\label{sufficiently_general_position_def1}  A collection of affine subspaces $P_1, ..., P_m$ of $\R^{d+1}$ and subset $R\subset\R^{d+1}$ is in sufficiently general position
if for each pairwise disjoint subsets $R_1, ..., R_m\subset R$ and for each subcollection $Q_1, ..., Q_k$ of collection $\aff(P_1\cup R_1), ..., \aff(P_m\cup R_m)$ of affine subspaces the following proposition holds:

$\bullet$ if $\bigcap\limits_{i=1}^kQ_i\not\subset \R^d\times 0$ then $\dim \bigcap\limits_{i=1}^k Q_i=-(d+1)(k-1)+\sum\limits_{i=1}^k \dim Q_k$.
\end{definition}

\begin{definition}\label{sufficiently_general_position_def2}  A collection of affine subspaces $P_1, ..., P_m$ of $\R^{d+1}$ and subset $R\subset\R^{d+1}$ is in strong sufficiently general position
if it in sufficiently general position and for each pairwise disjoint subsets $R_1, ..., R_m\subset R$ and for each subcollection $Q_1, ..., Q_k$ of collection $\aff(P_1\cup R_1), ..., \aff(P_m\cup R_m)$ of affine subspaces the following proposition holds:

$\bullet$ if $\bigcap\limits_{i=1}^kQ_i\neq \varnothing$ then $\dim (\R^d\times 0)\cap\bigcap\limits_{i=1}^k Q_i=-d(k-1)+\sum\limits_{i=1}^k \dim Q_k\cap(\R^d\times 0)$.
\end{definition}

\begin{lemma}\label{sufficiently_general_position_lemma1}
 For any collection of affine subspaces $P_1, ..., P_m$ and a subset $R$ of $\R^{d+1}$ in sufficiently general
position and any open set $U\subset\R^{d+1}$ there is $p\in U$ such that the flats $P_1, ..., P_m$ and subset $R\cup\{p\}$ are in sufficiently general
position.
Moreover, if the affine subspaces $P_1, ..., P_m$ and the subset $R$ of $\R^{d+1}$ in strong sufficiently general position then the set $R\cup\{p\}$ in strong sufficiently general position with the affine subspaces $P_1, ..., P_m$.
\end{lemma}

{\it Proof.}
For any ordered family $(R_1, ..., R_m)$ of subsets of $R$
and any subcollection $Q_1, ..., Q_k$ of the collection $\aff(P_1\cup R_1), ..., \aff(P_m\cup R_m)$  let us declare the following subspace forbidden:

$$\aff((\bigcap\limits_{i=1}^{k-1}Q_i)\cup Q_k),$$

if this subspace is not equal to $\R^{d+1}$.
Note that the number of forbidden subspaces is finite.
Take a point $p\in U$ which does not belong to any of forbidden subspaces and prove that  $p$ is the required point.
%For proving this Lemma it sufficiently to proove that for any $p\notin U$ the flats $\aff(P_j\cup p), \{P_i\}_{i\neq j, 1\leq i\leq m} $ are in sufficiently general position.

 Suppose that the  flats $P_1, ..., P_m$ and the subset $R\cup\{p\}$ are not in sufficiently general
position.
  Then there is pairwise disjoint subsets $R_1, ..., R_m$ of $R$ and
  a subcollection $Q_1, ..., Q_k$ of collection  $\aff(P_1\cup R_1), ..., \aff(P_m\cup R_m)$ such that
  $$(\bigcap\limits_{i=1}^{k-1}Q_i)\cap \aff(Q_k\cup\{p\})\not\subset \R^d\times 0$$
  and
  $$\dim(\bigcap\limits_{i=1}^{k-1}Q_i)\cap \aff(Q_k\cup\{p\})\neq1-(d+1)(k-1)+\sum\limits_{i=1}^k \dim Q_i.$$

  Then $\bigcap\limits_{i=1}^{k-1}Q_i\not\subset \R^d\times 0$.
  Hence $$\dim(\bigcap\limits_{i=1}^{k-1}Q_i)=-(d+1)(k-2)+\sum\limits_{i=1}^{k-1}\dim Q_i.$$

   The space $\aff(Q_k\cup\{p\})$ is not equal to $Q_k$ because the collection $Q_1, ..., Q_k$ is in sufficiently general position.
   Hence $\dim \aff(Q_k\cup\{p\})=1+\dim Q_k$.

 Since $(\bigcap\limits_{i=1}^{k-1}Q_i)\cap \aff(Q_k\cup\{p\})\neq\varnothing$, we have

 \begin{multline*} \dim(\bigcap\limits_{i=1}^{k-1}Q_i)\cap \aff(Q_k\cup\{p\})\geq \dim(\bigcap\limits_{i=1}^{k-1}Q_i)+\dim \aff(Q_k\cup\{p\})-(d+1)=\\
  =1-(d+1)+\dim(\bigcap\limits_{i=1}^{k-1}Q_i)+\dim Q_k=1-(d+1)(k-1)+\sum\limits_{i=1}^{k}\dim Q_i.
 \end{multline*}

 Then $\dim(\bigcap\limits_{i=1}^{k-1}Q_i)\cap \aff(Q_k\cup\{p\})> \dim(\bigcap\limits_{i=1}^{k-1}Q_i)+\dim \aff(Q_k\cup\{p\})-(d+1)$ because $\dim(\bigcap\limits_{i=1}^{k-1}Q_i)\cap \aff(Q_k\cup\{p\})\neq1-(d+1)(k-1)+\sum\limits_{i=1}^k \dim Q_k.$
 Therefore $\aff((\bigcap\limits_{i=1}^{k-1}Q_i)\cup \aff(Q_k\cup p))\neq\R^{d+1}$.
 Hence $\aff((\bigcap\limits_{i=1}^{k-1}Q_i)\cup Q_k)\neq\R^{d+1}$.

  %{\it Case 1. $\aff((\bigcap\limits_{i=1}^{k-1}Q_i)\cup Q_k)\neq\R^{d+1}$.}

 Since the subspace $\aff((\bigcap\limits_{i=1}^{k-1}Q_i)\cup Q_k)$ is forbidden, the point $p$ does not belong to this subspace.
 Hence $(\bigcap\limits_{i=1}^{k-1}Q_i)\cap\aff(Q_k\cup p)=(\bigcap\limits_{i=1}^{k-1}Q_i)\cap Q_k=\bigcap\limits_{i=1}^{k}Q_i$.
 Therefore
 \begin{multline*}
  \dim(\bigcap\limits_{i=1}^{k}Q_i)=\dim(\bigcap\limits_{i=1}^{k-1}Q_i)\cap \aff(Q_k\cup\{p\})\geq\\
  \geq 1-(d+1)+\dim(\bigcap\limits_{i=1}^{k-1}Q_i)+\dim Q_k>-(d+1)+\dim(\bigcap\limits_{i=1}^{k-1}Q_i)+\dim Q_k=\\
  =-(d+1)(k-1)+\sum\limits_{i=1}^{k}\dim Q_i.
 \end{multline*}
 Then $\bigcap\limits_{i=1}^{k}Q_i\subset\R^d\times0$ because the subspaces $Q_1, ..., Q_k$ are in sufficiently general position.
 A contradiction because $$\bigcap\limits_{i=1}^{k}Q_i=(\bigcap\limits_{i=1}^{k-1}Q_i)\cap \aff(Q_k\cup\{p\})\not\subset \R^d\times 0.$$
 First part of Lemma proved.

 Next, we suppose that the  flats $P_1, ..., P_m$ and the subset $R$ are in strong sufficiently general
position.
For any ordered family $(R_1, ..., R_m)$ of subsets of $R$
and any subcollection $Q_1, ..., Q_k$ of the collection $$\aff(P_1\cup R_1), ..., \aff(P_m\cup R_m), \R^d\times0$$  let us additionaly declare as forbidden a hyperplane which contains the subspace $Q_k$ and parallel to the subspace $\bigcap\limits_{i=1}^{k-1}Q_i$.
Take a point $p\in U$ which does not belong to any of forbidden subspaces and prove that  $p$ is the required point.

Suppose that the flats $P_1, ..., P_m$ and the subset $R\cup\{p\}$ are not in strong sufficiently general
position.
Then there is pairwise disjoint subsets $R_1, ..., R_m$ of $R$ and
  a subcollection $Q_1, ..., Q_k$ of collection  $\aff(P_1\cup R_1), ..., \aff(P_m\cup R_m)$ such that
  $$(\bigcap\limits_{i=1}^{k-1}Q_i)\cap \aff(Q_k\cup\{p\})\neq \varnothing$$
  and
  $$\dim(\R^d\times0)\cap(\bigcap\limits_{i=1}^{k-1}Q_i)\cap \aff(Q_k\cup\{p\})\neq -d(k-1)+\sum\limits_{i=1}^{k-1} \dim Q_i\cap(\R^d\times0)+\dim\aff(Q_k\cup\{p\}).$$

 if $Q_k\subset\R^d\times0$ then $(\R^d\times0)\cap\aff(Q_k\cup\{p\})=(\R^d\times0)\cap Q_k$.
 Hence
 \begin{gather*}
  \dim(\R^d\times0)\cap(\bigcap\limits_{i=1}^{k-1}Q_i)\cap \aff(Q_k\cup\{p\})=\dim(\R^d\times0)\cap\bigcap\limits_{i=1}^{k}Q_i=\\
  =-d(k-1)+\sum\limits_{i=1}^k \dim Q_k\cap(\R^d\times 0)=\\
  =-d(k-1)+\sum\limits_{i=1}^{k-1} \dim Q_k\cap(\R^d\times 0)+\dim\aff(Q_k\cup\{p\}).
 \end{gather*}
A contradiction.
Then $Q_k\not\subset\R^d\times0$.

  According to first part of Lemma, the flats $P_1, ..., P_m$ and the subset $R\cup\{p\}$ are in sufficiently general position.
If $(\bigcap\limits_{i=1}^{k-1}Q_i)\cap\aff(Q_k\cup\{p\})\not\subset\R^d\times0$ then

$$\dim Q_i\cap(\R^d\times0)=-1+\dim Q_i$$ for $i=1,...,k-1$ and
\begin{gather*}
  \dim(\R^d\times0)\cap(\bigcap\limits_{i=1}^{k-1}Q_i)\cap \aff(Q_k\cup\{p\})=-1+\dim(\bigcap\limits_{i=1}^{k-1}Q_i)\cap \aff(Q_k\cup\{p\})=\\
  =-1-(k-1)(d+1)+\sum\limits_{i=1}^{k-1}\dim Q_i+\dim \aff(Q_k\cup\{p\})=\\
  =-(k-1)d+\sum\limits_{i=1}^{k-1}\dim Q_i\cap(\R^d\times0)+\dim (\R^d\times0)\cap\aff(Q_k\cup\{p\}).
\end{gather*}
A contradiction.
Therefore $(\bigcap\limits_{i=1}^{k-1}Q_i)\cap\aff(Q_k\cup\{p\})\subset\R^d\times0$.

%$\dim (\R^d\times0)\cap\aff(Q_k\cup\{p\})=-1+\dim \aff(Q_k\cup\{p\})$;

Let us consider a subcases:

{\it Case 1: $(\R^d\times0)\cap\bigcap\limits_{i=1}^{k}Q_i\neq\varnothing$}

%Then $\dim (\R^d\times0)\cap(\bigcap\limits_{i=1}^{k-1}Q_i)+\dim Q_k=$

Since $\dim\aff(Q_k\cup\{p\})=1+\dim Q_k$, we have
\begin{gather*}
  \dim(\R^d\times0)\cap(\bigcap\limits_{i=1}^{k-1}Q_i)\cap \aff(Q_k\cup\{p\})=1+\dim(\R^d\times0)\cap(\bigcap\limits_{i=1}^{k}Q_i)=\\
  =1-(k-1)d+\sum\limits_{i=1}^{k}\dim (\R^d\times0)\cap Q_i=\\
  =-(k-1)d+\sum\limits_{i=1}^{k-1}\dim Q_i\cap(\R^d\times0)+\dim (\R^d\times0)\cap\aff(Q_k\cup\{p\}).
\end{gather*}
A contradiction.

{\it Case 2: $(\R^d\times0)\cap\bigcap\limits_{i=1}^{k}Q_i=\varnothing$}

 If $\aff Q_k\cup((\R^d\times0)\cap\bigcap\limits_{i=1}^{k-1}Q_i)\neq\R^{d+1}$, then $p\notin\aff Q_k\cup((\R^d\times0)\cap\bigcap\limits_{i=1}^{k-1}Q_i)$.
 Hence $(\R^d\times0)\cap(\bigcap\limits_{i=1}^{k-1}Q_i)\cap\aff(Q_k\cup\{p\})=(\R^d\times0)\cap\bigcap\limits_{i=1}^{k}Q_i=\varnothing$.
Therefore $(\bigcap\limits_{i=1}^{k-1}Q_i)\cap\aff(Q_k\cup\{p\})=\varnothing$ because $(\bigcap\limits_{i=1}^{k-1}Q_i)\cap\aff(Q_k\cup\{p\})\subset\R^d\times0$.
A contradiction.
Then $\aff Q_k\cup((\R^d\times0)\cap\bigcap\limits_{i=1}^{k-1}Q_i)=\R^{d+1}$.

Hence $\dim Q_k+\dim (\R^d\times0)\cap\bigcap\limits_{i=1}^{k-1}Q_i=d$.
Since the point $\{p\}$ does not belong to forbidden subspace, the interestion $ (\R^d\times0)\cap(\bigcap\limits_{i=1}^{k-1}Q_i)\cap \aff(Q_k\cup\{p\})$ is unique point.
Therefore
\begin{gather*}
 \dim(\R^d\times0)\cap(\bigcap\limits_{i=1}^{k-1}Q_i)\cap \aff(Q_k\cup\{p\})=0=\\
 =-d+\dim Q_k+\dim (\R^d\times0)\cap\bigcap\limits_{i=1}^{k-1}Q_i=1-d(k-1)+\sum\limits_{i=1}^k\dim (\R^d\times0)\cap Q_i=\\
 =-d(k-1)+\sum\limits_{i=1}^{k-1}\dim (\R^d\times0)\cap Q_i+\dim (\R^d\times0)\cap \aff(Q_k\cup\{p\}).
\end{gather*}
A contradiction.

Lemma is proved

\begin{lemma}\label{taking_higher_points}
 Given affine subspaces $L_1, ..., L_m\subset\R^d=\R^d\times0$ in general position at infinity with $\bigcap\limits_{i=1}^mL_i=\varnothing$.
 Let $A\subset\R^d\times0$.
 For any number $$q\leq\min(|A|-|\{i|L_i=\varnothing\}|;(m-1)d-\sum\limits_{i=1}^m\dim L_i)$$ and any subset $\mathfrak{U}\subset A$ such that $|\mathfrak{U}|=q$ there is:

 $\bullet$ a finite subset $A_{higher}\subset \R^d\times(1;+\infty)$;

 $\bullet$ a partition $\mathfrak{U}_1,...,\mathfrak{U}_m$ of $\mathfrak{U}$;

 $\bullet$ a pairwise disjoint subsets $\mathfrak{U}_1',...,\mathfrak{U}_m'$ of $A_{higher}\cup (A\setminus\mathfrak{U})$

 with properties:

 $\bullet$ $|(\R^d\times(0;1))\cap\bigcap\limits_{i=1}^mC(L_i,\mathfrak{U}_i', \mathfrak{U}_i)|=1$;

 $\bullet$ the subspaces $L_1, ..., L_m$ and the subset $A\cup A_{higher}$ of $\R^{d+1}$ in  sufficiently general position.

 Moreover, if the set $A$ in strongly independent with the subspaces $L_1,...,L_m$ as subsets of $\R^d=\R^d\times0$ we can require that the subspaces $L_1, ..., L_m$ and the subset $A\cup A_{higher}$ of $\R^{d+1}$ in strong sufficiently general position.
\end{lemma}

{\it Proof.}

Denote by $m_1$ a maximal number such that $\bigcap\limits_{i=1}^{m_1}L_i\neq\varnothing$.
Then $1\leq m_1\leq m-1$.
Let us construct the sequence $W_{1}, ..., W_{m_1}$ as follows:
Suppose that  the points $W_{1}, ..., W_{i}$ are constructed for some $i=0,...,m_1-1$.
Apply Lemma~\ref{sufficiently_general_position_lemma1} and take $W_{i+1}\in \R^{d}\times(1;+\infty)$ such that the subspaces $L_1, ..., L_m$ and the set $A\cup\{W_{1},...,W_{i}\}$ are in sufficiently general position.
Moreover, if the set  $A\cup\{W_{1},...,W_{i-1}\}$ in strong sufficiently general position with the subspaces $L_1,...,L_m$, we can require that set in $A\cup\{W_{1},...,W_{i}\}$ in strong sufficiently general position with the subspaces $L_1,...,L_m$.
The sequence $W_{1}, ..., W_{m_1}$ is constructed.

Since $q\leq \min(|A|-|\{i|L_i=\varnothing\}|;(m-1)d-\sum\limits_{i=1}^m\dim L_i)$
there is a pairwise disjoint subsets $\mathfrak{A}_1, ..., \mathfrak{A}_m, \mathfrak{B}_1, ..., \mathfrak{B}_m$  of $A$ such that:

$\bullet$ $\sum\limits_{i=1}^m|\mathfrak{B}_i|=q$;

$\bullet$ $\mathfrak{A}_i\cup L_i\neq \varnothing$ for any $i=1, ..., m$;

$\bullet$ $\sum\limits_{i=1}^{i_0} (|\mathfrak{A}_i|+|\mathfrak{B}_i|+1+\dim L_i)= (i_0-1)(d+1)$ for any $i_0=m_1+1, ..., m$.

Then $|\mathfrak{A}_i|+|\mathfrak{B}_i|+1+\dim L_i\leq d+1$ for any $i=1, ..., m$.

Let us construct the sequence $W_{m_1+1}, ..., W_{m}$ as follows:
Suppose that for some $i\in\{m_1+1, ..., m\}$  the points $W_{m_1+1}, ..., W_{i-1}$ are constructed with the properties:

$\bullet$ $(\bigcap\limits_{j=1}^{i-1}C(L_j, \{W_{j}\}\cup \mathfrak{A}_j, \mathfrak{B}_j))\cap(\R^d\times(0;1))\neq\varnothing$;

$\bullet$ the subspaces $L_1, ..., L_m$ and the set $$A\cup\{W_1,...,W_{i-1}\}$$ are in sufficiently general position.

%Since the subspaces $L_1, ..., L_m$ do not crosses at infinity, the intersection

Since $\sum\limits_{j=1}^{i} (|\mathfrak{A}_j|+|\mathfrak{B}_j|+1+\dim L_j)= (i-1)(d+1)$,
there is an open set $U_{i}\subset \R^{d}\times(1;+\infty)$ such that for any point $E_{i}\in U_{i}$:
 $$ C(L_{i}, \{E_{i}\}\cup \mathfrak{A}_{i}, \mathfrak{B}_{i})\cap(\bigcap\limits_{j=1}^{i-1}C(L_j, \{W_{j}\}\cup \mathfrak{A}_j, \mathfrak{B}_j))\cap(\R^d\times(0;1))\neq\varnothing.$$
  From Lemma~\ref{sufficiently_general_position_lemma1} there is a point $W_{i}\in U_{i}$ such that  the subspaces $L_1, ..., L_m$ and the set $A\cup\{W_{1},...,W_{i}\}$ are in sufficiently general position.
  Moreover, if the set  $A\cup\{W_{1},...,W_{i-1}\}$ in strong sufficiently general position with the subspaces $L_1,...,L_m$, we can require that set in $A\cup\{W_{1},...,W_{i}\}$ in strong sufficiently general position with the subspaces $L_1,...,L_m$.
  The sequence $W_{m_1+1}, ..., W_{m}$ is constructed.
  Take:

  $\bullet$ $A_{higher}=\{W_1,...,W_m\}$;

  $\bullet$ $\mathfrak{U}_i=\mathfrak{B}_i$ for $i=1,...,m$;

  $\bullet$ $\mathfrak{U}_i'=\mathfrak{A}_i\cup\{W_i\}$ for $i=1,...,m$.

  Since $\sum\limits_{i=1}^{m} (|\mathfrak{A}_i|+|\mathfrak{B}_i|+1+\dim L_i)= (m-1)(d+1)$, we have $$\sum\limits_{i=1}^{i_0} (|\mathfrak{U}_i|+|\mathfrak{U}_i'|+\dim L_i)= (m-1)(d+1).$$
  From sufficiently general position follows that $|\bigcap\limits_{i=1}^mC(L_i,\mathfrak{U}_i', \mathfrak{U}_i)|=1$.
  Lemma is proved.

\smallskip

\section{Linear optimization}

For any polyhedron $Q$ denote by $\partial Q$ the union of its proper faces.
We will consider polyhedra of both types:  bounded and unbounded.
%Denote by $f$ the last coordinate functional on $\R^{d+1}$.

\smallskip

\begin{lemma}[on linear optimization]\label{linear_optimization} Let

$\bullet$ $L_1, L_2, ... ,L_m$ be a collection of affine subspaces of $\R^d\times0\subset\R^{d+1}$;

$\bullet$ $f$ the last coordinate functional on $\R^{d+1}$;

$\bullet$ $R$ be a finite subset of $\R^d\times[0;+\infty)$ such that the collection $L_1, L_2, ... ,L_m$ and set $R$ are  in sufficiently general position;

$\bullet$ $Y_1, Y_2, ... ,Y_m, Y_1', Y_2' ... ,Y_m'$ be pairwise disjoint subsets of $R$;

$\bullet$ $W\in R\setminus\{Y_1\cup...\cup Y_m\cup Y_1'\cup...\cup Y_m'\}$;

$\bullet$ $Z\in\R^d\times(0;+\infty)$ is the unique point of intersection $\bigcap\limits_{i=1}^mC(L_i, Y_i, Y_i')$.

Then there is $j_c$ such that:

1) if $f(W)<f(Z)$ then the point $Z$ is not minimal point of the functional $f$ on the intersection
$C(L_{j_c}, (Y_{j_c}\cup\{W\}),Y_{j_c}')\cap(\bigcap\limits_{i\neq j_c}C(L_i, Y_i,Y_i'))$;

2) if $f(W)>f(Z)$ then the point $Z$ is not minimal point of the functional $f$ on the intersection
$C(L_{j_c}, Y_{j_c},(Y_{j_c}'\cup\{W\}))\cap(\bigcap\limits_{i\neq j_c}C(L_i, Y_i,Y_i'))$.

\end{lemma}

\smallskip

{\it Proof}.

%Consider some $(d+1)$-dimensional cube $K$ with center in $Z$. %The intersection $K\cap(\bigcap\limits_{i=1}^m[L_i, (Y_i\cap\{W\})))=\bigcap\limits_{i=1}^m([L_i, Y_i\cap\{W\})\cap K)$ is a polyhedron.
Denote by:

$\bullet$ $\Delta_i$ a polyhedron $C(L_i, Y_i, Y_i')$ for $i=1, ..., m$;

$\bullet$ $k$ a number of non-full-dimensional polyhedra from $\Delta_1, ..., \Delta_m$.
%$\bullet$ $f$ the last coordinate functional on $\R^{d+1}$.

According to the conditions $f(Z)\neq f(W)$.
For $i=1, ..., m$ define a polyhedron
\begin{equation*}
\delta_i =
    \begin{cases}
        C(L_i, Y_i\cup\{W\}, Y_i') &\text{if } f(W)<f(Z)\\
        C(L_i, Y_i, Y_i'\cup\{W\}) & \text{if }  f(W)>f(Z).
    \end{cases}
\end{equation*}

Without loss of generality  the polyhedra $\Delta_{1}, ..., \Delta_k$ are not full-dimensional.
Then for $i=1, ..., k$ if the polyhedron $\Delta_i$ is a proper face of the polyhedron $\delta_i$.
From sufficiently general position of the collection $L_1, L_2, ... ,L_m$ and set $R$ it follows that $Z\notin \partial \Delta_i$ for $i=1, ..., m$.
Since $Z$ is interior point of polyhedra $\Delta_{k+1}, ..., \Delta_m$, we have $\bigcap\limits_{i=1}^k\Delta_i=\{Z\}$.
From Lemma~\ref{edge_of_interestion_of_polyhedron} follows that the point $Z$ is a vertex of $\bigcap\limits_{i=1}^k\delta_i$.
Denote $\mathfrak{M}=\bigcap\limits_{i=1}^k\delta_i$.

Prove that the point $Z$ is not a minimum point of the functional $f$ on $\mathfrak{M}$.
Let us consider two cases.
First case $f(W)<f(Z)$, then $$W\in\bigcap_{i=1}^kC(L_i, Y_i\cup\{W\}, Y_i')=\mathfrak{M}.$$
Second case  $f(W)>f(Z)$, then $$C(\varnothing, Z, W)\subset\bigcap_{i=1}^kC(L_i, Y_i, Y_i'\cup\{W\})=\mathfrak{M}.$$
The point $Z$ is not a minimum point of the functional $f$ on $C(\varnothing, Z, W)$ in this case.
Therefore the point $Z$ is not a minimum point of the functional $f$ on $\mathfrak{M}$ in both cases.

 %The last coordinate of $W\in \mathfrak{M}$ equals 0.
 %The last coordinate of $Z\in \mathfrak{M}$ is larger than 0.
 Since the set $\mathfrak{M}$ is convex, the point $Z$ is not a local minimum point of $f$ on $\mathfrak{M}$.
 Then there is an edge $E$ of $\mathfrak{M}$ such that $Z\in E$ and $f(Z_1)<f(Z)$ for any $Z_1\in E\setminus{Z}$.
 %Take some $Z_1\in E\setminus{Z}$.
From Lemma~\ref{edge_of_interestion_of_polyhedron} follows that $E=\delta_{j_c}\cap(\bigcap\limits_{i\neq j_c, i\leq k}\Delta_i)$ for some $j_c\leq k$.
Since the polyhedra $\Delta_{k+1}, ..., \Delta_m$ are full-dimensional,
the intersection $\delta_{j_c}\cap(\bigcap\limits_{i\neq j_c}\Delta_i)$ is edge of $\mathfrak{M}$.
Take a point $Z_2\in(\delta_{j_c}\cap(\bigcap\limits_{i\neq j_c}\Delta_i))\setminus\{Z\}$.
 Then  $f(Z_2)<f(Z)$.
 %and the last coordinate of $Z_1$ is smaller than the last coordinate of $Z$.
 Lemma~\ref{linear_optimization} is proved.

\smallskip

Lemma~\ref{polyhedron_border} below is simple and so could be known.

\begin{lemma}\label{polyhedron_border} Take convex polyhedra $A_1, ..., A_n$ in $\R^d$ which have more than one common point.
Then  $\partial(A_1\cap ... \cap A_n)\subset(\partial A_1)\cup ... \cup (\partial A_n)$.
\end{lemma}

 {\it Proof}. %ïðèìåíèì ïðèìåíèì îáîáùåíèå ýòîé ëåììû íà n ìíîãîãðàííèêîâ

 Suppose that there is a point $P\in\partial(A_1\cap ... \cap A_n)$ such that $P \notin (\partial A_1)\cup ... \cup (\partial A_n)$.
 Take $Q\in (A_1\cap ... \cap A_n)-\partial(A_1\cap ... \cap A_n)$.
 Denote by $l$ the line such that $P, Q\in l$.
 Denote by $I_k$ the segment  $I_k=A_k\cap l$ for each $k=1, ..., n$ and denote $I=l\cap A_1\cap ... \cap A_n$.
 Then $I=I_1\cap...\cap I_n$.
 Since $P\in\partial(A_1\cap ... \cap A_n)$, the point $P$ is an endpoint of $I$.
 Since $I=I_1\cap...\cap I_n$, there is $j$ such that $P$ is an endpoint of $I_j$.
 Since $I_j=A_j\cap l$, some face of $A_j$ contains a point $P$.
 Hence $P\in\partial A_j\subset(\partial A_1)\cup ... \cup (\partial A_n)$.
 
 A contradiction.
 
 \smallskip

 \begin{lemma}\label{edge_of_interestion_of_polyhedron}
  Let

$\bullet$ $L_1, L_2, ... ,L_k$ be a collection of flats of $\R^d\times0\subset\R^{d+1}$;

$\bullet$ $R$ be a finite subset $\R^d\times[0;+\infty)$ such that the collection $L_1, L_2, ... ,L_k$ and set $R$ are  in sufficiently general position;

$\bullet$ $A_1, A_2, ... ,A_k, A', B_1, B_2 ... ,B_k, B'$ be pairwise disjoint subsets of $R$;

%$\bullet$ $A_1', A_2', ... ,A_k', B_1', B_2' ... ,B_k'$ be  subsets of $R$ such that ;

$\bullet$ $C(L_i, A_i, B_i)$ is a proper face of polyhedron $C(L_i, A_i\cup A', B_i\cup B')$ for $i=1, ..., k$;

$\bullet$ $Z\in\R^d\times(0;+\infty)$ is the unique point of intersection $\bigcap\limits_{i=1}^kC(L_i, A_i, B_i)$.

Then:

$\bullet$ $Z$ is a vertex of polyhedron $\mathfrak{M}:=\bigcap\limits_{i=1}^kC(L_i, A_i\cup A', B_i\cup B')$;

$\bullet$ for any edge $E$ of polyhedron $\mathfrak{M}$ such that $Z\in E$ there is a number $j$ such that $E$ is an edge of $C(L_j, A_j\cup A', B_j\cup B')\cap(\bigcap\limits_{i\neq j}C(L_i, A_i, B_i))$.

 \end{lemma}

{\it Proof}.
Denote by:

$\bullet$ $\Delta_i$ a polyhedron $C(L_i, A_i, B_i)$ for $i=1, ..., k$;

$\bullet$ $\delta_i$ a polyhedron $C(L_i, A_i\cup A', B_i\cup B')$ for $i=1, ..., k$.

From sufficiently general position of the collection $L_1, L_2, ... ,L_k$ and set $R$ it follows that $Z\notin \partial \Delta_i$ for $i=1, ..., k$.

Prove that the point $Z$ is a vertex of $\mathfrak{M}$.
Assume to the contrary.
Then there is a non-degenerate segment $I\subset \mathfrak{M}$ with midpoint in $Z$.
Hence for any $i=1,..., k$ we have $I\subset\delta_i$.
Since $Z\in\Delta_i$ and $Z$ is midpoint of $I$, the segment $I$ contains in the set  $\Delta_i$ for $i=1,...,k$.
Therefore $I\subset \bigcap\limits_{i=1}^k\Delta_i$.
A contradiction, because $Z$ is the unique point of intersection $\bigcap\limits_{i=1}^k\Delta_i$.
Hence the point $Z$ is a vertex of $\mathfrak{M}$.

Take a minimal by inclusion subset $J\subset\{1, ..., k\}$ such that $E$ is an edge of $(\bigcap\limits_{j\in J} \delta_j)\cap(\bigcap\limits_{j\notin J} \Delta_j)$.
 Denote $\mathfrak{M}_J=(\bigcap\limits_{j\in J} \delta_j)\cap(\bigcap\limits_{j\notin J} \Delta_j)$.
 Suppose that $|J|\geq2$.
 For any $j\in J$ the intersection $\delta_j\cap(\bigcap\limits_{i\neq j}\Delta_i)$ at least 1-dimensional because the subspaces $L_1, L_2, ... ,L_m$ and the set $R$ are  in sufficiently general position.
 Then $\dim\mathfrak{M}_J\geq2$.
 %Take $Z_1\in E\setminus\{Z\}$.
 Hence $E\subset\partial\mathfrak{M}_J$ because $\dim\mathfrak{M}_J\geq2$.
 Since $Z\notin\partial\Delta_j$ for $1\leq j\leq k$, there is a neighbourhood $U$ of the point $Z$
 which does not intersect any face of polyhedra $\delta_1, ... , \delta_k$ besides $\Delta_1, ... , \Delta_k$.
 Take a point $Z_1\in E\cap (U\setminus\{Z\})$.
 Then $Z_1\in\partial\mathfrak{M}_J$.
 Apply Lemma~\ref{polyhedron_border}  to the polyhedron $\mathfrak{M}_J$.
 We get that $Z_1\in\partial\delta_{j_b}$ for some $j_b\in J$.
 From the choice of $U$ it follows that $Z_1\in\Delta_{j_b}$.
 Since $Z\in\Delta_{j_b}$, we have $E\subset\aff(Z,Z_1)\subset\aff(\Delta_{j_b})$.
 Hence $E\subset\aff(\Delta_{j_b})\cap\delta_{j_b}=\Delta_{j_b}$.
 Therefore $E\subset(\bigcap\limits_{j\in J-j_b} \delta_j)\cap(\bigcap\limits_{j\notin J-j_b} \Delta_j)$.
 Since $E$ is an edge of $\mathfrak{M}$ and $$E\subset(\bigcap\limits_{j\in J-j_b} \delta_j)\cap(\bigcap\limits_{j\notin J-j_b} \Delta_j)\subset\mathfrak{M},$$ the set $E$ is an edge of $(\bigcap\limits_{j\in J-j_b} \delta_j)\cap(\bigcap\limits_{j\notin J-j_b} \Delta_j)$.
 A contradiction with the definition of $J$.
 Hence $|J|=1$.
Take $j$ such that $J=\{j\}$.
Lemma is proved.

 \smallskip

 \begin{lemma}\label{add_negative_mass}
  In the conditions of Lemma~\ref{linear_optimization} if $f(W)=0<f(Z)$ and for some $j$:

  $\bullet$ $(\R^d\times0)\cap C(L_{j}, (Y_{j}\cup\{W\}),Y_{j}')\cap\aff(\bigcap\limits_{i \neq j} C(L_i, Y_i,Y_i'))\neq \varnothing$;

  $\bullet$ $L_j=Y_{j}'=\varnothing$;

  %$\bullet$ $|Y_{j}|\leq d+1$;

  then there is a number $j_{decr}$ such that the minimal point of the functional $f$ on the intersection
  $$C(L_{j_{decr}}, Y_{j_{decr}},Y_{j_{decr}}'\cup\{W\})\cap (\bigcap\limits_{i\neq j_{decr}}C(L_i, Y_i,Y_i'))$$
is not $Z$.

 \end{lemma}

{\it Proof}.
Since
\begin{gather*}
 (\R^d\times0)\cap C(L_{j}, (Y_{j}\cup\{W\}),Y_{j}')\cap\aff(\bigcap\limits_{i \neq j} C(L_i, Y_i,Y_i'))\neq \varnothing= (\R^d\times0)\cap\{Z\},
\end{gather*}
we have
$$C(L_{j}, (Y_{j}\cup\{W\}),Y_{j}')\cap\aff(\bigcap\limits_{i \neq j} C(L_i, Y_i,Y_i'))\neq\{Z\}.$$
Hence
$$C(L_{j}, (Y_{j}\cup\{W\}),Y_{j}')\neq C(L_{j}, Y_{j},Y_{j}').$$
Therefore the polyhedron $C(L_{j}, Y_{j},Y_{j}')=\left<Y_{j}\right>$ is not full-dimensional.
Then $|Y_{j}|\leq d+1$.

Denote by $k$ a number of non-full-dimensional polyhedra from $C(L_i, Y_i, Y_i')$ for $i=1, ..., m$.
Without loss of generality a polyhedra $C(L_i, Y_i, Y_i')$ for $i=1, ..., k$ are not full-dimensional and $j\leq k$.
Then:

$\bullet$ $Z$ is the unique point of intersection $\bigcap\limits_{i=1}^kC(L_i, Y_i, Y_i')$;

$\bullet$ a polyhedron $C(L_i, Y_i, Y_i')$ is a proper face of $C(L_i, Y_i,Y_i'\cap\{W\})$ for $i=1, ..., k$.

From Lemma~\ref{edge_of_interestion_of_polyhedron} follows that the point $Z$ is vertex of $\bigcap\limits_{i=1}^k C(L_i, Y_i,Y_i'\cap\{W\})$.
Denote $\mathfrak{M}=\bigcap\limits_{i=1}^k C(L_i, Y_i,Y_i'\cup\{W\})$.

Take $Z_2\in (C(L_{j}, (Y_{j}\cup\{W\}),Y_{j}')\cap(\bigcap\limits_{i \neq j} C(L_i, Y_i,Y_i')))\setminus\{Z\}$.
The line $\aff(Z, Z_2)$ intersects the set $$(\R^d\times0)\cap C(L_{j}, (Y_{j}\cup\{W\}),Y_{j}')=(\R^d\times0)\cap \left<Y_{j}\cup\{W\}\right>.$$
Denote by $Z_2'$ a point of this intersection.
Then $C(\varnothing, Z_2', W)\subset\R^d\times0$.
The intersection  $C(\varnothing, Z_2', W)\cap C(L_{j}, Y_{j},Y_{j}')=C(\varnothing, Z_2', W)\cap \left<Y_{j}\right>$ is the one-point set, because $Z_2'\in \left<Y_{j}\cup\{W\}\right>$.
Denote by $Z_2''$ a point such that $\{Z_2''\}=C(\varnothing, Z_2', W)\cap C(L_{j_c}, Y_{j_c},Y_{j_c}')$.
Then
$$C(\varnothing, \{Z, Z_2'\}, W)\cap C(L_{j}, Y_{j},Y_{j}')=\left<Z, Z_2''\right>.$$
Hence $$
 \left<Z, Z_2''\right>\subset C(L_{j}, Y_{j},Y_{j}')\cap(\bigcap_{i\neq j} C(L_i, Y_i,Y_i'\cup\{W\}))\subset\mathfrak{M}.$$

%Denote by $\delta_1, ..., \delta_k$ a  polyhedra from the collection $\{C(L_i, Y_i,Y_i'\cap\{W\})\}_{i=1}^k$.
 The point $Z$ is not a minimum point of the functional $f$ on $\left<Z, Z_2''\right>$.
Therefore the point $Z$ is not a minimum point of the functional $f$ on  $\mathfrak{M}$.
 Since the set $\mathfrak{M}$ is convex, the point $Z$ is not a local minimum point of $f$ on $\mathfrak{M}$.
Since $Z$ is a vertex of $\mathfrak{M}$, there is an edge $E_1$ of  $\mathfrak{M}$ such that $Z\in E_1$ and $f(Z)>f(Z_3)$ for any $Z_3\in E_1\setminus\{Z\}$. %$Z$ is a maximum point of $f$ on $E_1$.
From Lemma~\ref{edge_of_interestion_of_polyhedron} follows that there is a number $j_{decr}$ such that $$E_1=C(L_{j_{decr}}, Y_{j_{decr}},Y_{j_{decr}}'\cup\{W\})\cap (\bigcap\limits_{i\neq j_{decr}; i\leq k}C(L_i, Y_i,Y_i')).$$
Then the point $Z$ is not minimal point of the functional $f$ of the intersection $$C(L_{j_{decr}}, Y_{j_{decr}},Y_{j_{decr}}'\cup\{W\})\cap (\bigcap\limits_{i\neq j_{decr}}C(L_i, Y_i,Y_i'))$$
because $Z$ is the internal point of $C(L_i, Y_i,Y_i')$ for $i=(k+1), ..., m$.

Lemma is proved.

\smallskip

\begin{lemma}\label{add_negative_mass_to_full_dimensional_simplex_0}
 Let:

 $\bullet$ L is an affine subspace of $\R^d$;

 $\bullet$ $R$ be finite subset of $\R^d$;

 $\bullet$ $\aff(L\cup R)=\R^d$, $|R|=d-\dim L$;

 $\bullet$ $N$ is interior point of $C(L, A, B)$ for some $A,B\subset R$;

  $\bullet$ $\omega\in C(L, A', B')$ for some $A',B'\subset R$;

  $\bullet$ $\omega\notin \aff(L\cup R')$ for any subset $R'\subset R\cup\{N\}$ such that $|R'|\leq d-\dim L-1$;

 $\bullet$ $|A'|\geq |B'|+2$.

 Then there are pairwise disjoint subsets $A'', B''\subset R\cup\{N\}$ such that

 $\bullet$ $\omega$ is interior point of $C(L, A'', B'')$;

 $\bullet$ $|A''|+ |B''|=|R|$;

 $\bullet$ $|A''|=|A'|-1$, $|B''|=|B'|+1$.

 \end{lemma}

  {\it Proof}.

  For any $x\in L\cup R$ denote by $\alpha(x), \beta(x)$ the numbers such that:

  $$\sum\limits_{x\in L\cup R}\alpha(x)x=N$$

  and

 $$\sum\limits_{x\in L\cup R}\beta(x)x=\omega.$$

 Both sums have only finite number of non-zero terms.

 Hence
 $$tN+\sum\limits_{x\in L\cup R}(\beta(x)-t\alpha(x))x=\omega$$
 for any $t\in\R$.

 Since $N$ is interior point of $C(L, A, B)$, for $x\in R$ we have $\alpha(x)\neq 0$.
 Since $\omega$ is interior point of $C(L, A', B')$, for $x\in R$ we have $\beta(x)\neq 0$.

 Suppose that there are $x', x''\in R$ such that $\frac{\beta(x')}{\alpha(x')}=\frac{\beta(x'')}{\alpha(x'')}$.
 Then the line $\aff\{N, \omega\}$ intersects an affine subspace $\aff(L\cup (R\setminus\{x',x''\}))$.
 Hence $\omega\in \aff(L\cup(R\setminus\{x',x''\})\cup\{N\})$.
 A contradiction because $|(R\setminus\{x',x''\})\cup\{N\}|= |R|-1 = d-\dim L-1$.
 Hence the relations $\frac{\beta(x)}{\alpha(x)}$ for $x\in R$ are pairwise distinct.

 Denote:

 $$R_1=\{x|x\in R, \frac{\beta(x)}{\alpha(x)}>0\}$$

 and

 $$R_2=\{x|x\in R, \frac{\beta(x)}{\alpha(x)}<0\}.$$

 Since $|A'|\geq |B'|+2$, at least one of the following two inequalities holds:

 1) $|\{x|x\in R_1, \beta(x)>0\}|\geq |\{x|x\in R_1, \beta(x)<0\}|+2$;

 2) $|\{x|x\in R_2, \beta(x)>0\}|\geq |\{x|x\in R_2, \beta(x)<0\}|+1$.

 Let us consider two cases.

 {\it Case 1: $|\{x|x\in R_1, \beta(x)>0\}|\geq |\{x|x\in R_1, \beta(x)<0\}|+2$.}

 Let us enumerate the points of $R_1$: $x_1, ..., x_{|R_1|}$ such that

 $$0<\frac{\beta(x_1)}{\alpha(x_1)}<\frac{\beta(x_2)}{\alpha(x_2)}<...<\frac{\beta(x_{|R_1|})}{\alpha(x_{|R_1|})}.$$
 There is a number $j$ such that $\beta(x_j)>0$ and  $$|\{i|i<j,\beta(x_i)>0\}|=1+|\{i|i<j,\beta(x_i)<0\}|.$$
 Take $t=\frac{\beta(x_j)}{\alpha(x_j)}$. Then:

 $\bullet$ $\beta(x_j)-t\alpha(x_j)=0$;

 $\bullet$ $\beta(x_i)(\beta(x_i)-t\alpha(x_i))<0$ for $i=1, ..., j-1$;

 $\bullet$ $\beta(x)(\beta(x)-t\alpha(x))>0$ for $x\in R\setminus\{x_1,...,x_j\}$.

 In the set of numbers $\{t\}\cup\{\beta(x)-t\alpha(x)|x\in R\}$ exactly $|B'|+1$  negative elements, because in the set of numbers $\{\beta(x)|x\in R\}$ exactly $|B'|$  negative elements.

 Take $B''=\{x|x\in R,\beta(x)-t\alpha(x)<0\}$, $A''=\{N\}\cup\{x|x\in R,\beta(x)-t\alpha(x)>0\}$.
 Then $\omega \in C(L, A'', B'')$ because $tN+\sum\limits_{x\in L\cup R}(\beta(x)-t\alpha(x))x=\omega$.
 Lemma proved in this case.

 {\it Case 2: $|\{x|x\in R_2, \beta(x)>0\}|\geq |\{x|x\in R_2, \beta(x)<0\}|+1$.}

 Let us enumerate the points of $R_2$: $y_1, ..., y_{|R_2|}$ such that

 $$0>\frac{\beta(y_1)}{\alpha(y_1)}>\frac{\beta(y_2)}{\alpha(y_2)}>...>\frac{\beta(y_{|R_2|})}{\alpha(y_{|R_2|})}.$$

 There is a number $j$ such that $\beta(x_j)>0$ and  $$|\{i|i<j,\beta(x_i)>0\}|=|\{i|i<j,\beta(x_i)<0\}|.$$
 Take $t=\frac{\beta(x_j)}{\alpha(x_j)}$. Then:

 $\bullet$ $\beta(x_j)-t\alpha(x_j)=0$;

 $\bullet$ $\beta(x_i)(\beta(x_i)-t\alpha(x_i))<0$ for $i=1, ..., j-1$;

 $\bullet$ $\beta(x)(\beta(x)-t\alpha(x))>0$ for $x\in R\setminus\{x_1,...,x_j\}$.

 In the set of numbers $\{t\}\cup\{\beta(x)-t\alpha(x)|x\in R\}$ exactly $|B'|+1$  negative elements, because in the set of numbers $\{\beta(x)|x\in R\}$ exactly $|B'|$  negative elements.

 Take $B''=\{N\}\cup\{x|x\in R,\beta(x)-t\alpha(x)<0\}$, $A''=\{x|x\in R,\beta(x)-t\alpha(x)>0\}$.
 Then $\omega \in C(L, A'', B'')$ because $tN+\sum\limits_{x\in L\cup R}(\beta(x)-t\alpha(x))x=\omega$.
 Lemma proved in bots cases.

\smallskip

\begin{lemma}\label{add_negative_mass_to_full_dimensional_simplex}
In the conditions of Lemma~\ref{linear_optimization}
suppose that $f(W)=0$ and there is $j$ such that:

 $\bullet$  the polyhedron $C(L_j,Y_j,Y_j')$ is $d+1$-dimensional;

 $\bullet$   $|Y_j'|\leq \lfloor \frac{d-\dim L_j}{2}\rfloor$;

  $\bullet$   $Y_j'\subset\R^d\times0$;

 $\bullet$ $|Y_j\setminus(\R^d\times0)|=1$.

Then for any $n$ from $|Y_j'|$ to $\lfloor \frac{d-\dim L_j}{2}\rfloor$ there are pairwise disjoint subsets $P, Q\subset Y_j\cup Y_j'\cup \{W\}$ such that:

 $\bullet$ $Z\in C(L_j,P,Q)$;

 $\bullet$ $|Q|=n$;

 $\bullet$ $Q\subset \R^d\times0$;

 $\bullet$ $|P|+|Q|=d+1-\dim L_j$;

 $\bullet$ $(Y_j\cup Y_j'\cup \{W\})\setminus(P\cup Q)$ is the one point subset of $\R^d\times0$.

\end{lemma}

  {\it Proof}. We will prove by induction on $n$.

  Base case: For $n=|Y_j'|$ take $P=Y_j, Q=Y_j'$.

  Induction step:  Suppose that for some $n$ such that:
  $|Y_j'|<n\leq\lfloor \frac{d-\dim L_j}{2}\rfloor$ there are pairwise disjoint subsets $P', Q'\subset Y_j\cup Y_j'\cup \{W\}$ such that:

 $\bullet$ $Z\in C(L_j,P',Q')$;

 $\bullet$ $|Q'|=n-1$;

 $\bullet$ $Q'\subset \R^d\times0$;

 $\bullet$ $|P'|+|Q'|=d+1-\dim L_j$;

 $\bullet$ $(Y_j\cup Y_j'\cup \{W\})\setminus(P'\cup Q')$ is the one point subset of $\R^d\times0$.

 Since the polyhedron $C(L_j,Y_j,Y_j')$ is $d+1$-dimensional, from sufficiently general position follows that $\{Z\}=\bigcap\limits_{i\neq j} C(L_i,Y_i,Y_i')$.

  Denote by $D$ the unique point of $Y_j\setminus(\R^d\times0)$.
  Denote by $W'$ a point such that $\{W'\}=(Y_j\cup Y_j'\cup \{W\})\setminus(P'\cup Q')\subset\R^d\times0$.
  Since $Y_j', Q\subset\R^d\times0$ and $W\in\R^d\times0$, we have  $P'\setminus(\R^d\times0)=Y_j\setminus(\R^d\times0)=\{D\}$.
 From sufficiently general position follows that the point $Z$ does not belong to any affine subspace $ \aff(L_j\cup K)$ for any subset $K\subset Y_j\cup Y_j'\cup \{W\}$ such that $|K|\leq d-\dim L_j$ because $\{Z\}=\bigcap\limits_{i\neq j} C(L_i,Y_i,Y_i')$.
  Hence the ray $C(\varnothing, D, Z)$ does not intersect subspace $\aff(L_i\cup R')$ for any subset $R'\subset (Y_j\setminus \{D\})\cup Y_j'\cup \{W\}\subset\R^d\times0$ such that $|R'|\leq d-\dim L_j-1$.
  Denote by $\omega$ the unique point of intersection $C(\varnothing, D, Z)\cap(\R^d\times0)$.
 Then $\omega\in C(L_j, P', Q')\cap (\R^d\times0)=C(L_j,P'\cap(\R^d\times0),Q')$.
 Take $P_1=P'\cap(\R^d\times0)$.
 The polyhedron $C(L_j, P_1, Q')$ is $d$-dimensional and $C(L_j, P_1, Q')\subset \R^d\times0$.
 Hence $|P_1|+|Q'|=d-\dim L_j$.
 Therefore $|P_1|-|Q'|\geq2$
 because $|Q'|=n-1\leq\lfloor \frac{d-\dim L_j}{2}\rfloor-1$.
 Since $|P_1|+|Q'|=d-\dim L_j$, the point $W'$ belong to the interior of some $d$-dimensional polyhedron $C(L_j, P_2, Q_2)$ for some $P_2, Q_2\subset P_1, Q'$.
 From Lemma~\ref{add_negative_mass_to_full_dimensional_simplex_0} follows that there are pairwise disjoint subsets $P_3,Q_3\subset P_1\cup Q'\cup\{W'\}$ such that:

 $\bullet$ $\omega$ is interior point of $C(L_j, P_3, Q_3)$;

 $\bullet$ $|P_3|+|Q_3|=d-\dim L_j$;

 $\bullet$ $|P_3|=|P_1|-1$, $|Q_3|=|Q'|+1=n$.

 Take $P=P_3\cup\{D\}$, $Q=Q_3$.
 Then:

 $\bullet$ $Z\in \left<D,\omega\right>\subset C(L_j,P,Q)$;

 $\bullet$ $|Q|=n$;

 $\bullet$ $Q\subset \R^d\times0$;

 $\bullet$ $|P|+|Q|=d+1-\dim L_j=|Y_j|+|Y_j'|$.

 Hence $|Y_j\cup Y_j'\cup \{W\}|=1+|P\cup Q|$.
 Since $(Y_j\cup Y_j'\cup \{W\})\setminus(\R^d\times0)=\{D\}\subset P$, we have $(Y_j\cup Y_j'\cup \{W\})\setminus(P\cup Q)\subset \R^d\times0$.
 Therefore subsets $P,Q$ are required.

 The step of induction is completed.

\smallskip

\begin{lemma}\label{add_negative_mass_to_full_dimensional_simplexes}
In the conditions of Lemma~\ref{linear_optimization}
suppose  that:

$\bullet$ $f(W)=0$;

$\bullet$ $\sum\limits_{j\in J}|Y_j'|\leq \sum\limits_{j\in J}\lfloor \frac{d-\dim L_j}{2}\rfloor$;

and there is a subset $J\subset\{1,...,m\}$ such that for any $j\in J$:

 $\bullet$  the polyhedron $C(L_j,Y_j,Y_j')$ is $d+1$-dimensional;

  $\bullet$   $Y_j'\subset\R^d\times0$;

 $\bullet$ $|Y_j\setminus(\R^d\times0)|=1$.

Then for any $n$ from $\sum\limits_{j\in J}|Y_j'|$ to $\sum\limits_{j\in J}\lfloor \frac{d-\dim L_j}{2}\rfloor$ there is a family $\{P_j\}_{j\in J}\cup\{Q_j\}_{j\in J}$ pairwise disjoint subsets of $ \{W\}\cup\bigcup\limits_{j\in J}(Y_j\cup Y_j') $ such that
  $$\sum\limits_{j\in J}|Q_j|=n;$$
  $$\varnothing\neq (\{W\}\cup\bigcup\limits_{j\in J}(Y_j\cup Y_j') )\setminus (\bigcup\limits_{j\in J}(P_j\cup Q_j))\subset \R^d\times0$$
  and for any $j\in J$:

 $\bullet$ $Z\in C(L_j,P_j,Q_j)$;

 $\bullet$ $Q_j\subset \R^d\times0$;

 $\bullet$ $|P_j|+|Q_j|=d+1-\dim L_j$.

\end{lemma}

  {\it Proof}. We will prove by induction on $|J|$.

  Base case: it is Lemma~\ref{add_negative_mass_to_full_dimensional_simplex}.

  Induction step:

  Take a numbers $n_j$ for $j\in J$ such that:

  $\bullet$ $|Y_j'|\leq n_j\leq \max(\lfloor \frac{d-\dim L_j}{2}\rfloor; |Y_j'|)$ for $j\in J$;

  $\bullet$ $\sum\limits_{j\in J}n_j=n$.

  There is a number $j_1\in J$  such that
  $$\sum\limits_{j\in J\setminus{\{j_1\}}}|Y_j'|\leq \sum\limits_{j\in J\setminus{\{j_1\}}}n_j\leq \sum\limits_{j\in J\setminus{\{j_1\}}}\lfloor \frac{d-\dim L_j}{2}\rfloor.$$
  Apply induction hypothesis on the set $J\setminus\{j_1\}$ and the number $\sum\limits_{j\in J\setminus{\{j_1\}}}n_j$.
 There  is a family $\{P_j\}_{j\in  J\setminus{\{j_1\}}}\cup\{Q_j\}_{j\in  J\setminus{\{j_1\}}}$ pairwise disjoint subsets of $ \{W\}\cup\bigcup\limits_{j\in  J\setminus{\{j_1\}}}(Y_j\cup Y_j') $ such that
 $$\sum\limits_{j\in J\setminus{\{j_1\}}}|Q_j|=\sum\limits_{j\in J\setminus{\{j_1\}}}n_j$$
 $$\varnothing\neq (\{W\}\cup\bigcup\limits_{j\in J\setminus{\{j_1\}}}(Y_j\cup Y_j') )\setminus (\bigcup\limits_{j\in J\setminus{\{j_1\}}}(P_j\cup Q_j))\subset \R^d\times0$$
 and for any $j\in J\setminus{\{j_1\}}$:

 $\bullet$ $Z\in C(L_j,P_j,Q_j)$;

 $\bullet$ $Q_j\subset \R^d\times0$;

 $\bullet$ $|P_j|+|Q_j|=d+1-\dim L_j$.

 Take a point $W'\in (\{W\}\cup\bigcup\limits_{j\in J\setminus{\{j_1\}}}(Y_j\cup Y_j'))\setminus (\bigcup\limits_{j\in J\setminus{\{j_1\}}}(P_j\cup Q_j))$.
 Then $W'\in\R^d\times0$.

 Apply Lemma~\ref{add_negative_mass_to_full_dimensional_simplex} for sets $Y_{j_1}, Y_{j_1}'$, number $n_{j_1}$ and point $W'$.
 There are pairwise disjoint subsets $P_{j_1}, Q_{j_1}\subset Y_{j_1}\cup Y_{j_1}'\cup \{W'\}$ such that:

 $\bullet$ $Z\in C(L_{j_1},P_{j_1},Q_{j_1})$;

 $\bullet$ $|Q_{j_1}|=n_{j_1}$;

 $\bullet$ $Q_{j_1}\subset \R^d\times0$;

 $\bullet$ $|P_{j_1}|+|Q_{j_1}|=d+1-\dim L_{j_1}$.

 $\bullet$ $(Y_{j_1}\cup Y_{j_1}'\cup \{W'\})\setminus(P_{j_1}\cup Q_{j_1})$ is the one point subset of $\R^d\times0$.

Hence $$\varnothing\neq (\{W\}\cup\bigcup\limits_{j\in J}(Y_j\cup Y_j') )\setminus (\bigcup\limits_{j\in J}(P_j\cup Q_j))=(Y_{j_1}\cup Y_{j_1}'\cup \{W'\})\setminus(P_{j_1}\cup Q_{j_1})\subset\R^d\times0.$$
 Therefore the family $\{P_j\}_{j\in J}\cup\{Q_j\}_{j\in J}$ is required.

 The step of induction is completed.
\smallskip

\section{Proof by optimization}

{\it Proof of Theorem~\ref{main_result} in conditions 1,2,3}.

Note that if $\bigcap\limits_{i=1}^mL_i\neq\varnothing$ then statement of Theorem~\ref{main_result} holds.
Next, we suppose that  $\bigcap\limits_{i=1}^mL_i=\varnothing$.
From general position at infinity follows that $\sum\limits_{i=1}^m\dim L_i\leq d(m-1)-1$.
Therefore
$$\sum\limits_{i=1}^m(1+\dim L_i)\leq(d+1)(m-1)$$
and
$$|A|\geq 1.$$

Consider $\R^d$ as $\R^d\times0\subset\R^{d+1}$.

Note that inequality $k\leq(m-1)d-\sum\limits_{i=1}^m\dim L_i$ holds in any conditions 1,2,3 of Theorem~\ref{main_result}.

Take a subset $A_{higher}\subset\R^d\times(1;+\infty)$ from Lemma~\ref{taking_higher_points}.

  Denote by $f$ the last coordinate functional on $\R^{d+1}$.
  For $q=0,..., k$ take a pairwise disjoint subsets $Y_1^{(q)+}, ..., Y_m^{(q)+}, Y_1^{(q)-}, ..., Y_m^{(q)-}$ of $A\cup A_{higher}$, if it exists, such that:

$\bullet$ $Y_j^{(q)-}\subset A$ for $j=1, ..., m$;

$\bullet$ $\sum\limits_{j=1}^{m}|Y_j^{(q)-}|=q$;

$\bullet$ $|Y_j^{(q)+}\setminus A|=1$ for $j=1, ..., m$;

$\bullet$ $\bigcap\limits_{j=1}^{m}C(L_j, Y_j^{(q)+}, Y_j^{(q)-})$ is an unique point set; denote this point by $Z^{(q)}$;

$\bullet$ $f(Z^{(q)})>0$;

$\bullet$ $f(Z^{(q)})$ is minimal with respect to previous conditions.

According to the choice of $A_{higher}$ the subsets $Y_1^{(k)+}, ..., Y_m^{(k)+}, Y_1^{(k)-}, ..., Y_m^{(k)-}$ are exists.

Since the subspaces $L_1, ..., L_m$ and the set $A\cup A_{higher}$ are in sufficiently general position, we have
$$\sum\limits_{j=1}^{m}(\dim L_j +1+ | Y_j^{(q)+}|+| Y_j^{(q)-}|)-(m-1)(d+1)=0.$$
then the set $A\setminus(Y_1^{(q)+}\cup ...\cup Y_m^{(q)+}\cup Y_1^{(q)-}\cup ...\cup Y_m^{(q)-})$ is not empty.
Denote by $T^{(q)}$ a some point from $A\setminus(Y_1^{(q)+}\cup ...\cup Y_m^{(q)+}\cup Y_1^{(q)-}\cup ...\cup Y_m^{(q)-})$.

Let us prove some propositions.

\begin{proposition}\label{add_positive mass}
 If there is no a point $Z^{(k-1)}\in \R^d\times(0;f(Z^{(k)}))$ then Theorem~\ref{main_result} holds for this $k$.
\end{proposition}

{\it Proof of Proposition.}

From Lemma~\ref{linear_optimization} follows that there is a number $j_c(k)$ such that the point $Z^{(k)}$ is not a minimum point of the functional $f$ on the set
\begin{gather}\label{segment_from_Z}
C(L_{j_c(k)}, \{T^{(k)}\}\cup Y_{j_c(k)}^{(k)+}, Y_{j_c(k)}^{(k)-})\cap\bigcap\limits_{i\neq j_c(k)}C(L_i, Y_i^{(k)+}, Y_i^{(k)-})
\end{gather}
because $f(T^{(k)})=0<f(Z^{(k)})$.
%Take some $j_c(k)$ for any $k=0,...,k$.
 %In the case when
 Suppose that the set~\eqref{segment_from_Z} does not intersect the subspace $\R^d\times0$.
 Then the minimum of the functional $f$ on the set~\eqref{segment_from_Z} is bigger than $0$.
Denote by $Z_2^{(k)}$ a minimal point of the functional $f$ on the set~\eqref{segment_from_Z}.
Then the point $Z_2^{(k)}$ belong to some proper face of the polyhedron~\eqref{segment_from_Z}.
From Lemma~\ref{polyhedron_border} follows that
\begin{gather*}
Z_2^{(k)}\in \partial C(L_{j_c(k)}, \{T^{(k)}\}\cup Y_{j_c(k)}^{(k)+}, Y_{j_c(k)}^{(k)-})\cup\bigcup\limits_{j\neq j_c(k)}\partial C(L_j, Y_j^{(k)+}, Y_j^{(k)-}).
\end{gather*}
Then there are a subsets:
\begin{equation}\label{inclusions}
\begin{cases}
 Y_i^{(k)*+}\subset Y_i^{(k)+}, & \mbox{for }   i\neq j_c(k)\\
 Y_{ j_c(k)}^{(k)*+}\subset \{T^{(k)}\}\cup Y_{ j_c(k)}^{(k)+}\\
 Y_{i}^{(k)*-}\subset Y_i^{(k)-}, & \mbox{for }   i=1 , ..., m
\end{cases}
\end{equation}
such that $\{Z_2^{(k)}\}=\bigcap\limits_{i=1}^mC(L_j, Y_j^{(k)*+}, Y_j^{(k)*-})$.
%Take these subsets $Y_1^{(k)*+}, ..., Y_m^{(k)*+}, Y_1^{(k)*-}, ..., Y_m^{(k)*-}$ for this $k$.
Since the polyhedron~\eqref{segment_from_Z} is 1-dimensional, we have
\begin{multline*}
 \sum\limits_{i=1}^m\dim C(L_j, Y_j^{(k)*+}, Y_j^{(k)*-})=(m-1)(d+1)=\\
 =-1+\dim C(L_{j_c(k)}, \{T^{(k)}\}\cup Y_{j_c(k)}^{(k)+}, Y_{j_c(k)}^{(k)-})+\sum\limits_{j\neq j_c(k)}\dim C(L_j, Y_j^{(k)+}, Y_j^{(k)-}).
\end{multline*}
Therefore
\begin{multline*}
 \sum\limits_{i=1}^m(|Y_j^{(k)*+}|+|Y_j^{(k)*-}|)=\\
 =-1+| \{T^{(k)}\}\cup Y_{j_c(k)}^{(k)+}|+|Y_{j_c(k)}^{(k)-}|+\sum\limits_{j\neq j_c(k)}(|Y_j^{(k)+}|+|Y_j^{(k)-}|).
 \end{multline*}
 From the inclusions~\eqref{inclusions} follows either:
\begin{gather*}
 \sum\limits_{i=1}^m|Y_j^{(k)*-}| =-1+\sum\limits_{i=1}^m|Y_j^{(k)-}|=k-1
 \end{gather*}
or
\begin{gather*}
 \sum\limits_{i=1}^m|Y_j^{(k)*-}| =\sum\limits_{i=1}^m|Y_j^{(k)-}|=k.
 \end{gather*}
 If $\sum\limits_{i=1}^m|Y_j^{(k)*-}| =k$ we have a contradiction because $f(Z_2^{(k)})<f(Z^{(k)})$.
  If $\sum\limits_{i=1}^m|Y_j^{(k)*-}| =k-1$ then there is the point $Z^{(k-1)}$.
  From the definition of $Z^{(k-1)}$ follows that $f(Z^{(k-1)})\leq f(Z_2^{(k)})$.
A contradiction because $f(Z_2^{(k)})<f(Z^{(k)})\leq f(Z^{(k-1)})$.% and $\{Z_2^{(k)}\}=\bigcap\limits_{i=1}^mC(L_j, Y_j^{(k)*+}, Y_j^{(k)*-})$, we have

Hence the set~\eqref{segment_from_Z}  intersect the subspace $\R^d\times0$.
Take:

 $\bullet$ $A_i=Y_i^{(k)+}\cap(\R^d\times0)$ for $i\in\{1, ..., m\}\setminus\{j_c(k)\}$;

 $\bullet$ $B_i=Y_i^{(k)-}\cap(\R^d\times0)$ for $i\in\{1, ..., m\}$;

  $\bullet$  $A_{j_c(k)}=A\setminus((\bigcup\limits_{i\neq j_c(k)}A_i)\cup\bigcup\limits_{i=1}^mB_i)$.

  Then $(Y_{j_c(k)}^{(k)+}\cup\{T^{(k)}\})\cap(\R^d\times0)\subset A_{j_c(k)}$.
  Hence
  \begin{multline*}
\varnothing\neq(\R^d\times0)\cap C(L_{j_c(k)}, \{T^{(k)}\}\cup Y_{j_c(k)}^{(k)+}, Y_{j_c(k)}^{(k)-})\cap\bigcap\limits_{i\neq j_c(k)}C(L_i, Y_i^{(k)+}, Y_i^{(k)-})\subset\\ \subset\bigcap\limits_{i=1}^mC(L_i, A_i, B_i)
\end{multline*}

Therefore $A_1, ..., A_m, B_1, ..., B_m$ are required partition.

Proposition is proved.

\begin{remark}\label{Roudneffs_case}
 The point $Z^{(-1)}$ does not exists.
 %Then $Z^{(-1)}\notin\R^d\times(0;f(Z^{(0)}))$.
 From Proposition~\ref{add_positive mass} follows  Theorem~\ref{main_result} for $k=0$.
\end{remark}

\begin{proposition}\label{add_negative_mass_2}
 Suppose for some $k\geq 1$ there is a number $i_e$ such that $Z^{(k-1)}$ is not a minimum point of the functional $f$ on the set
\begin{gather}\label{segment_from_Z_k-1}
C(L_{i_e}, Y_{i_e}^{(k-1)+}, \{T^{(k-1)}\}\cup Y_{i_e}^{(k-1)-})\cap\bigcap\limits_{i\neq i_e}C(L_i, Y_i^{(k-1)+}, Y_i^{(k-1)-}).
\end{gather}
Then Theorem~\ref{main_result} holds for this $k$.
\end{proposition}

{\it Proof:}
If $f(Z^{(k-1)})\geq f(Z^{(k)})$ then Theorem~\ref{main_result} holds according to Proposition~\ref{add_positive mass}. Next we suppose that $f(Z^{(k-1)})< f(Z^{(k)})$.

Suppose that the set~\eqref{segment_from_Z_k-1} does not intersect the subspace $\R^d\times0$.
 Then the minimum of the functional $f$ on the set~\eqref{segment_from_Z_k-1} is bigger than $0$.
Denote by $Z_3^{(k-1)}$ a minimal point of the functional $f$ on the set~\eqref{segment_from_Z_k-1}.
Then the point $Z_3^{(k-1)}$ belong to some proper face of the polyhedron~\eqref{segment_from_Z_k-1}.
From Lemma~\ref{polyhedron_border} follows that
\begin{gather*}
Z_3^{(k-1)}\in \partial C(L_{i_e}, Y_{i_e}^{(k-1)+},\{T^{(k-1)}\}\cup  Y_{i_e}^{(k-1)-})\cup\bigcup\limits_{j\neq i_e}\partial C(L_j, Y_j^{(k-1)+}, Y_j^{(k-1)-}).
\end{gather*}
Then there are a subsets:
\begin{equation}\label{inclusions}
\begin{cases}
 Y_i^{(k-1)**-}\subset Y_i^{(k-1)-}, & \mbox{for }   i\neq i_e\\
 Y_{ i_e}^{(k-1)**-}\subset \{T^{(k-1)}\}\cup Y_{ i_e}^{(k-1)}\\
 Y_{i}^{(k-1)**+}\subset Y_i^{(k-1)+}, & \mbox{for }   i=1 , ..., m
\end{cases}
\end{equation}
such that $\{Z_3^{(k-1)}\}=\bigcap\limits_{i=1}^mC(L_j, Y_j^{(k-1)**+}, Y_j^{(k-1)**-})$.
%Take these subsets $Y_1^{(k-1)**+}, ..., Y_m^{(k-1)**+}, Y_1^{(k-1)**-}, ..., Y_m^{(k-1)**-}$ for this $k$.
Since the polyhedron~\eqref{segment_from_Z_k-1} is 1-dimensional, we have
\begin{multline*}
 \sum\limits_{i=1}^m\dim C(L_j, Y_j^{(k-1)**+}, Y_j^{(k-1)**-})=(m-1)(d+1)=\\
 =-1+\dim C(L_{i_e}, Y_{i_e}^{(k-1)+}, \{T^{(k-1)}\}\cup Y_{i_e}^{(k-1)-})+\sum\limits_{j\neq i_e}\dim C(L_j, Y_j^{(k-1)+}, Y_j^{(k-1)-}).
\end{multline*}
Therefore
\begin{multline*}
 \sum\limits_{i=1}^m(|Y_j^{(k-1)**+}|+|Y_j^{(k-1)**-}|)=\\
 =-1+|Y_{i_e}^{(k-1)+}|+| \{T^{(k-1)}\}\cup Y_{i_e}^{(k-1)-}|+\sum\limits_{j\neq i_e}(|Y_j^{(k-1)+}|+|Y_j^{(k-1)-}|).
 \end{multline*}
 From the inclusions~\eqref{inclusions} follows either:
\begin{gather*}
 \sum\limits_{i=1}^m|Y_j^{(k-1)**-}| =-1+\sum\limits_{i=1}^m|Y_j^{(k-1)-}|=k-1
 \end{gather*}
or
\begin{gather*}
 \sum\limits_{i=1}^m|Y_j^{(k-1)**-}| =\sum\limits_{i=1}^m|Y_j^{(k-1)-}|=k.
 \end{gather*}
A contradiction because $f(Z_3^{(k-1)})<f(Z^{(k-1)})< f(Z^{(k)})$.% and $\{Z_3^{(k-1)}\}=\bigcap\limits_{i=1}^mC(L_j, Y_j^{(k-1)**+}, Y_j^{(k-1)**-})$, we have

Hence the set~\eqref{segment_from_Z_k-1}  intersect the subspace $\R^d\times0$.
Take:

 $\bullet$ $A_i=Y_i^{(k-1)+}\cap(\R^d\times0)$ for $i\in\{1, ..., m\}$;

 $\bullet$ $B_i=Y_i^{(k-1)-}\cap(\R^d\times0)$ for $i\in\{1, ..., m\}\setminus\{i_e\}$;

  $\bullet$  $B_{i_e}=A\setminus((\bigcup\limits_{i\neq i_e}B_i)\cup\bigcup\limits_{i=1}^mA_i)$.

  Then $(Y_{i_e}^{(k-1)-}\cup\{T^{(k-1)}\})\cap(\R^d\times0)\subset B_{i_e}$.
  Hence
  \begin{multline*}
\varnothing\neq(\R^d\times0)\cap C(L_{i_e}, Y_{i_e}^{(k-1)+}, \{T^{(k-1)}\}\cup Y_{i_e}^{(k-1)-})\cap\bigcap\limits_{i\neq i_e}C(L_i, Y_i^{(k-1)+}, Y_i^{(k-1)-})\subset\\ \subset\bigcap\limits_{i=1}^mC(L_i, A_i, B_i)
\end{multline*}

Therefore $A_1, ..., A_m, B_1, ..., B_m$ are required partition.

Proposition is proved.

\begin{proposition}\label{optimization_on_line}
Let:

 $\bullet$ $\Gamma_1^+, ..., \Gamma_m^+, \Gamma_1^-, ..., \Gamma_m^-$ be pairwise disjoint subsets of  $A\cup A_{higher}$;

 $\bullet$ $\Gamma_1^-, ..., \Gamma_m^-\subset A$;

 $\bullet$ $|\Gamma_i^+\setminus A|=1$ for any $i=1, ..., m$;

 $\bullet$ $\dim\bigcap\limits_{i=1}^mC(L_i,\Gamma_i^+,\Gamma_i^-)=1$;

 $\bullet$ $|\Gamma_i^+|+|\Gamma_i^-|\leq d+1-\dim L_i$ for $i\in\{1,...,m\}$.

 Denote by $Z$ a minimal point of the functional $f$ on the set $\bigcap\limits_{i=1}^mC(L_i,\Gamma_i^+,\Gamma_i^-)$.
 Suppose that $\{Z\}=C(L_j,\Gamma_j^+,\Gamma_j^{-*})\cap\bigcap\limits_{i\neq j}C(L_i,\Gamma_i^+,\Gamma_i^-)$ for some $j$ and proper subset $\Gamma_j^{-*}\subset\Gamma_j^-$.

 Suppose  $\dim\bigcap\limits_{i\in \mathfrak{p}_1}C(L_i,\Gamma_i^+,\Gamma_i^-)=1$ for some subset $\mathfrak{p}_1\subset\{1,...,m\}$.
  Denote by $l$ the line: $l=\aff \bigcap\limits_{i=1}^mC(L_i,\Gamma_i^+,\Gamma_i^-)$.
  For any $i\in \mathfrak{p}_1$ take a partition $\Xi_i^+,\Xi_i^-$ of the set $ \Gamma_i^+\cup\Gamma_i^-$ such that $f(Z)\geq\max\limits_{\bigcap\limits_{i\in\mathfrak{p}_1}C(L_i,\Xi_i^+,\Xi_i^-)}f>0$.
  Denote $\mathfrak{p}_2=\{1, ..., m\}\setminus \mathfrak{p}_1$.
  Suppose that $$\sum\limits_{i=1}^m|\Gamma_i^-|\leq (\sum\limits_{i\in \mathfrak{p}_2}\lfloor\frac{d-\dim L_i}{2}\rfloor)+\sum\limits_{i\in \mathfrak{p}_1}|\Xi_i^-|.$$

  Then for any $t$ from $\sum\limits_{i=1}^m|\Gamma_i^-|-1$ to $(\sum\limits_{i\in \mathfrak{p}_2}\lfloor\frac{d-\dim L_i}{2}\rfloor)-1+\sum\limits_{i\in \mathfrak{p}_1}|\Xi_i^-|$

  there are a pairwise disjoint subsets $$\Phi_1^+, ..., \Phi_m^+, \Phi_1^-, ..., \Phi_m^-\subset \bigcup\limits_{i=1}^m(\Gamma_i^+\cup\Gamma_i^-)$$ such that:

%$\bullet$ $\bigcap\limits_{i=1}^mC(L_i,\Phi_i^+,\Phi_i^-)\not\subset\R^d\times0$;

$\bullet$ $\Phi_1^-,...,\Phi_m^-\subset A$;

$\bullet$ $\max\limits_{\bigcap\limits_{i=1}^mC(L_i,\Phi_i^+,\Phi_i^-)}f >0$;

$\bullet$ $\min\limits_{\bigcap\limits_{i=1}^mC(L_i,\Phi_i^+,\Phi_i^-)}f < f(Z)$;

$\bullet$  $\sum\limits_{i=1}^m|\Phi_i^-|\in\{t,t+1\}$;

$\bullet$ either   1) $\sum\limits_{i=1}^m|\Phi_i^-|=t+1$ or 2) $\dim\bigcap\limits_{i=1}^mC(L_i,\Phi_i^+,\Phi_i^-) =0$.

%Moreover if $l\cap(\R^d\times0)\subset C(L_i,\Xi_i^+,\Xi_i^-)$
\end{proposition}

{\it Proof:} %Denote by $Z$ the minimal point of $f$ in the intersection $\bigcap\limits_{i=1}^mC(L_i,\Gamma_i^+,\Gamma_i^-)$.

  Take $Z_0=Z$.
Denote by $T_0$ a point such that $\{T_0\}=\Gamma_j^{-}\setminus \Gamma_j^{-*}$.
Take:

 $\bullet$ $\Phi_{1,j}^+=\Gamma_j^+\cup\{T_0\}$;

 $\bullet$ $\Psi_{0,j}^+=\Gamma_j^+$;

 $\bullet$ $\Psi_{0,j}^-= \Phi_{1,j}^-=\Gamma_j^{-*}$;

 $\bullet$ $\Psi_{0,i}^+=\Phi_{1,i}^+=\Gamma_i^+, \Psi_{0,i}^-=\Phi_{1,i}^-=\Gamma_i^{-}$ for $i\neq j$.

Then:

$\bullet$ $Z$ is maximal point of functional $f$ on the set $\bigcap\limits_{i=1}^mC(L_i,\Phi_{1,i}^+,\Phi_{1,i}^-)$;

$\bullet$ $\{Z\}=\bigcap\limits_{i=1}^mC(L_i,\Psi_{0,i}^+,\Psi_{0,i}^-)$;

$\bullet$  $\sum\limits_{i=1}^m|\Phi_{1,i}^-|= \sum\limits_{i=1}^m|\Gamma_i^-|-1$.

Let us consider a subcases.

{\it Case 1: $ \max\limits_{\bigcap\limits_{i\in\mathfrak{p}_1}C(L_i,\Xi_i^+,\Xi_i^-)}f =f(Z)$}

 Since $\max\limits_{\bigcap\limits_{i\in\mathfrak{p}_1}C(L_i,\Xi_i^+,\Xi_i^-)}f =f(Z)=\max\limits_{\bigcap\limits_{i\in\mathfrak{p}_1}C(L_i,\Phi_{1,i}^+,\Phi_{1,i}^-)}f$, both subsets $\bigcap\limits_{i\in\mathfrak{p}_1}C(L_i,\Xi_i^+,\Xi_i^-)$ and $\bigcap\limits_{i=1}^mC(L_i,\Phi_{1,i}^+,\Phi_{1,i}^-)$ of line $l$ have common relative interior point.
  Then for $i\in\mathfrak{p}_1$ the polyhedra $C(L_i,\Phi_{1,i}^+,\Phi_{1,i}^-)$ and $C(L_i,\Xi_{i}^{+},\Xi_{i}^{-})$ have common relative interior point.
  Hence $\Phi_{1,i}^+=\Xi_{i}^{+}$, $\Phi_{1,i}^-=\Xi_{i}^{-}$ for any $i\in\mathfrak{p}_1$ because $\Xi_{i}^{+}\cup\Xi_{i}^{-}=\Phi_{1,i}^+\cup\Phi_{1,i}^-$ and $|\Phi_{1,i}^+\cup\Phi_{1,i}^-|\leq d+1-\dim L_i$.
 Therefore $\sum\limits_{i\in \mathfrak{p}_1}|\Psi_{0,i}^-|=\sum\limits_{i\in \mathfrak{p}_1}|\Phi_{1,i}^-|=\sum\limits_{i\in\mathfrak{p}_1}|\Xi_{i}^{-}|$.
 Then
 \begin{multline*}
  (\sum\limits_{i\in\mathfrak{p}_2}\lfloor\frac{d-\dim L_i}{2}\rfloor) -1+\sum\limits_{i\in\mathfrak{p}_1}|\Xi_{i}^{-}| \geq t\geq\\
   \geq\sum\limits_{i=1}^m|\Phi_{1,i}^-|  = \sum\limits_{i\in\mathfrak{p}_2}|\Phi_{1,i}^-|+\sum\limits_{i\in\mathfrak{p}_1}|\Psi_{0,i}^{-}|\geq \sum\limits_{i\in\mathfrak{p}_2}|\Phi_{1,i}^-| -1+ \sum\limits_{i\in\mathfrak{p}_1}|\Xi_{i}^{-}|.
 \end{multline*}

From Lemma~\ref{add_negative_mass_to_full_dimensional_simplexes} follows that there is a family $\{\Phi_i^-\}_{i\in \mathfrak{p}_2}\cup\{\Phi_i^+\}_{i\in \mathfrak{p}_2}$ pairwise disjoint subsets of $ \{T_{0}\}\cup\bigcup\limits_{i\in \mathfrak{p}_2}(\Gamma_{i}^{+}\cup \Gamma_{i}^{-}) $ such that
  $$\sum\limits_{i\in \mathfrak{p}_2}|\Phi_i^-|=t+1- \sum\limits_{i\in\mathfrak{p}_1}|\Xi_{i}^{-}|$$ and for any $i\in \mathfrak{p}_2$:

 $\bullet$ $Z\in C(L_i,\Phi_i^+,\Phi_i^-)$;

 $\bullet$ $\Phi_i^-\subset \R^d\times0$;

 $\bullet$ $|\Phi_i^+|+|\Phi_i^-|=d+1-\dim L_i$.

 From Lemma~\ref{linear_optimization} follows that there is a number $j_0$ such that the point $Z$ is not minimal point of the functional $f$ on the intersection
$C(L_{j_0}, \Psi_{0, j_0}^+\cup\{T'\},\Psi_{0, j_0}^-)\cap(\bigcap\limits_{i\in \mathfrak{p}_1\setminus \{j_0\}}C(L_i, \Psi_{0, i}^+,\Psi_{0, i}^-))$.

Take $\Phi_i^+=\Psi_{0, i}^+$, $\Phi_i^-=\Psi_{0, i}^-$ for $i\in \mathfrak{p}_1\setminus\{j_0\}$, $\Phi_{j_0}^+=\Psi_{0, j_0}^+\cup\{T'\}$, $\Phi_{j_0}^-=\Psi_{0, j_0}^-$.

 Then

$\bullet$ $\Phi_1^-,...,\Phi_m^-\subset A$;

$\bullet$ $\max\limits_{\bigcap\limits_{i=1}^mC(L_i,\Phi_i^+,\Phi_i^-)}f = f(Z) >0$;

$\bullet$ $\min\limits_{\bigcap\limits_{i=1}^mC(L_i,\Phi_i^+,\Phi_i^-)}f < f(Z)$;

$\bullet$ $\sum\limits_{i=1}^m|\Phi_i^-|=t+1$.

Then sets $\Phi_1^+, ..., \Phi_m^+, \Phi_1^-, ..., \Phi_m^-$ are required.

{\it Case 2: $ \max\limits_{\bigcap\limits_{i\in\mathfrak{p}_1}C(L_i,\Xi_i^+,\Xi_i^-)}f <f(Z)$}

For some integer $r\geq 1$ let us construct the system of pairwise disjoint subsets $\{\Phi_{r,1}^+, ..., \Phi_{r,m}^+, \Phi_{r,1}^-, ..., \Phi_{r,m}^-\}$ of $A\cup A_{higher}$ and point $Z_{r-1}\in l$ with properties:

$\bullet$ $\Phi_{r,1}^-, ..., \Phi_{r,m}^-\subset A$;

 $\bullet$ $\dim\bigcap\limits_{i=1}^mC(L_i,\Phi_{r,i}^+,\Phi_{r,i}^-)=1$;

 $\bullet$ $\Phi_{r,i}^+\cup \Phi_{r,i}^-=\Gamma_{i}^+\cup \Gamma_{i}^-$ for any $i\in \mathfrak{p}_1$;
 %$\bullet$ $\bigcap\limits_{i=1}^mC(L_i,\Phi_{p,i}^+,\Phi_{p,i}^-)\subset l$;;

  $\bullet$ $\max\limits_{\bigcap\limits_{i=1}^mC(L_i,\Phi_{r,i}^+,\Phi_{r,i}^-)}f=f(Z_{r-1})$

  $\bullet$ $\sum\limits_{i=1}^m|\Phi_{r,i}^-|\leq t$;

 $\bullet$ $\max\limits_{\bigcap\limits_{i\in\mathfrak{r}_1}C(L_i,\Xi_i^+,\Xi_i^-)}f \leq f(Z_{r-1})\leq f(Z)$.

Since $\bigcap\limits_{i=1}^m \aff(L_i\cup\Gamma_{i}^+\cup\Gamma_{i}^-)=l=\bigcap\limits_{i\in\mathfrak{p}_1} \aff(L_i\cup\Gamma_{i}^+\cup\Gamma_{i}^-)$, we have
\begin{multline*}
 -(d+1)(m-1)+\sum\limits_{i=1}^m \dim \aff(L_i\cup\Gamma_{i}^+\cup\Gamma_{i}^-)=\dim \bigcap\limits_{i=1}^m \aff(L_i\cup\Gamma_{i}^+\cup\Gamma_{i}^-)=\\
=\dim \bigcap\limits_{i\in \mathfrak{p}_1} \aff(L_i\cup\Gamma_{i}^+\cup\Gamma_{i}^-) =-(d+1)(|\mathfrak{p}_1|-1)+\sum\limits_{i\in \mathfrak{p}_1} \aff(L_i\cup\Gamma_{i}^+\cup\Gamma_{i}^-).
\end{multline*}
Hence $\sum\limits_{i\in \mathfrak{p}_2} \dim \aff(L_i\cup\Gamma_{i}^+\cup\Gamma_{i}^-)=(d+1)|\mathfrak{p}_2|$.
Therefore for $i\in \mathfrak{p}_2$ the polyhedron $C(L_i,\Gamma_{i}^+,\Gamma_{i}^-)$ is $(d+1)$-dimensional.
From sufficiently general position follows that the point $Z$ is interior point of polyhedron $C(L_i,\Gamma_{i}^+,\Gamma_{i}^-)$ for $i\in \mathfrak{p}_2$.

For $r=1$:

The subsets $\{\Phi_{1,1}^+, ..., \Phi_{1,m}^+, \Phi_{1,1}^-, ..., \Phi_{1,m}^-\}$ are constructed with properties:

$\bullet$ $\Phi_{1,1}^-, ..., \Phi_{1,m}^-\subset A$;

 $\bullet$ $\dim\bigcap\limits_{i=1}^mC(L_i,\Phi_{1,i}^+,\Phi_{1,i}^-)=1$;

 $\bullet$ $\Phi_{1,i}^+\cup \Phi_{1,i}^-=\Gamma_{i}^+\cup \Gamma_{i}^-$ for any $i\in \mathfrak{p}_1$;

$\bullet$  $\max\limits_{\bigcap\limits_{i\in\mathfrak{r}_1}C(L_i,\Xi_i^+,\Xi_i^-)}f <f(Z)=f(Z_0)=\max\limits_{\bigcap\limits_{i=1}^mC(L_i,\Phi_{1,i}^+,\Phi_{1,i}^-)}f$;

$\bullet$  $\sum\limits_{i=1}^m|\Phi_{1,i}^-|= \sum\limits_{i=1}^m|\Gamma_i^-|-1\leq t$.

For $r\geq 2$:

Suppose that for $r-1=p\geq 1$ the pairwise disjoint subsets $\{\Phi_{p,1}^+, ..., \Phi_{p,m}^+, \Phi_{p,1}^-, ..., \Phi_{p,m}^-\}$ of $A\cup A_{higher}$ are constructed.% for some $p$ with properties:

  Since $\Phi_{p,i}^+\cup \Phi_{p,i}^-=\Gamma_{i}^+\cup \Gamma_{i}^-$ for any $i\in \mathfrak{p}_1$, we have $$\bigcap\limits_{i\in \mathfrak{p}_1}\aff(L_i\cup\Phi_{p,i}^+\cup\Phi_{p,i}^-)=\bigcap\limits_{i\in\mathfrak{p}_1}\aff(L_i\cup\Gamma_{i}^+\cup \Gamma_{i}^-)=l.$$
  Then $$\dim \bigcap\limits_{i\in \mathfrak{p}_1}\aff(L_i\cup\Phi_{p,i}^+\cup\Phi_{p,i}^-)=1=\dim \bigcap\limits_{i=1}^m\aff(L_i\cup\Phi_{p,i}^+\cup\Phi_{p,i}^-).$$
  Hence $$-(d+1)(|\mathfrak{p}_1|-1)+ \sum\limits_{i\in \mathfrak{p}_1}\dim\aff(L_i\cup\Phi_{p,i}^+\cup\Phi_{p,i}^-)=-(d+1)(m-1)+ \sum\limits_{i=1}^m\dim\aff(L_i\cup\Phi_{p,i}^+\cup\Phi_{p,i}^-).$$
  Therefore $$(d+1)|\mathfrak{p}_2|= \sum\limits_{i\in\mathfrak{p}_2}\dim\aff(L_i\cup\Phi_{p,i}^+\cup\Phi_{p,i}^-).$$
  Then $\dim C(L_i,\Phi_{p,i}^+,\Phi_{p,i}^-)=d+1$ for $i\in\mathfrak{p}_2$.

Denote by $Z_p$ the minimal point of $f$ on the intersection $\bigcap\limits_{i=1}^mC(L_i,\Phi_{p,i}^+,\Phi_{p,i}^-)$.
According to Lemma~\ref{polyhedron_border} there are subsets
$\Psi_{p,1}^+\subset\Phi_{p,1}^+, ..., \Psi_{p,m}^+\subset\Phi_{p,m}^+, \Psi_{p,1}^-\subset\Phi_{p,1}^-, ..., \Psi_{p,m}^-\subset\Phi_{p,m}^-$ such that $\bigcap\limits_{i=1}^mC(L_i,\Psi_{p,i}^+,\Psi_{p,i}^-)=\{Z_p\}$. Then
\begin{gather*}
 \dim \bigcap\limits_{i=1}^m C(L_i,\Psi_{p,i}^+,\Psi_{p,i}^-)=\dim\{Z_p\}=0=-1+\dim\bigcap\limits_{i=1}^mC(L_i,\Phi_{p,i}^+,\Phi_{p,i}^-).
\end{gather*}
Hence
\begin{gather*}
 \sum\limits_{i=1}^m(|\Psi_{p,i}^+|+|\Psi_{p,i}^-|)=-1+\sum\limits_{i=1}^m(|\Phi_{p,i}^+|+|\Phi_{p,i}^-|).
\end{gather*}
Then there is a number $j_1$ such that $|\Psi_{p,j_1}^+|+|\Psi_{p,j_1}^-|=-1+|\Phi_{p,j_1}^+|+|\Phi_{p,j_1}^-|$.

Consider a point $T_p\in(\Phi_{p,j_1}^+\cup\Phi_{p,j_1}^-)\setminus(\Psi_{p,j_1}^+\cup\Psi_{p,j_1}^-)$.
Since $Z_p\notin \R^d\times0$, the polyhedron $C(L_{j_1},\Psi_{p,j_1}^+,\Psi_{p,j_1}^-)$ is not subset of $\R^d\times0$.
Therefore $\Psi_{j_1}^+\not\subset\R^d\times0$.
Since $\Psi_{p,j_1}^+\subset\Phi_{p,j_1}^+$ and $|\Phi_{p,j_1}^+\setminus A|=1$, we have $\Psi_{p,j_1}^+\setminus A=\Phi_{p,j_1}^+\setminus A$.
Then $T_p\in A$.

From $\sum\limits_{i=1}^m(|\Psi_{p,i}^+|+|\Psi_{p,i}^-|)=-1+\sum\limits_{i=1}^m(|\Phi_{p,i}^+|+|\Phi_{p,i}^-|)$ follows that the sum $\sum\limits_{i=1}^m|\Psi_{p,i}^-|$ equals  either $\sum\limits_{i=1}^m|\Phi_{p,i}^-|$
or
$-1+\sum\limits_{i=1}^m|\Phi_{p,i}^-|$.
Since $\dim\bigcap\limits_{i=1}^mC(L_i,\Phi_{p,i}^+,\Phi_{p,i}^-)=1$, we have $f(Z_p)<f(Z)$.

Denote by $Z'$ the point of maximum of the functional $f$ on the set $\bigcap\limits_{i\in\mathfrak{p}_1}C(L_i,\Xi_i^+,\Xi_i^-)$. %Else, denote by $Z'$ the point of minimum of the functional $f$ on the set $\bigcap\limits_{i\in\mathfrak{p}_1}C(L_i,\Xi_i^+,\Xi_i^-)$.
We have $f(Z)>f(Z')>0$ and $\{Z'\}$ is a proper face of the polyhedron $\bigcap\limits_{i\in\mathfrak{p}_1}C(L_i,\Xi_i^+,\Xi_i^-)$.
According to Lemma~\ref{polyhedron_border} there are subsets
$\Xi_i^{+*}\subset\Xi_i^+$, $\Xi_i^{-*}\subset\Xi_i^-$ for $i\in\mathfrak{p}_1$ such that:
$$\{Z'\}=\bigcap\limits_{i\in\mathfrak{p}_1}C(L_i,\Xi_i^{+*},\Xi_i^{-*}).$$
From sufficiently general position follows that the point $Z'$ is relative interior point of $C(L_i,\Xi_i^{+*},\Xi_i^{-*})$ for any $i\in\mathfrak{p}_1$.

Let us consider a subcases:

{\it Case 2.1: $f(Z') < f(Z_{p-1})$.}

Then $f(Z') < f(Z_{p-1}) = \max\limits_{\bigcap\limits_{i=1}^mC(L_i,\Phi_{p,i}^+,\Phi_{p,i}^-)}f \leq \max\limits_{\bigcap\limits_{i\in \mathfrak{p}_1}C(L_i,\Phi_{p,i}^+,\Phi_{p,i}^-)}f$.
%Therefore either: $$\bigcap\limits_{i\in\mathfrak{p}_1}C(L_i,\Xi_i^+,\Xi_i^-)\subset\bigcap\limits_{i=1}^mC(L_i,\Phi_{p,i}^+,\Phi_{p,i}^-)$$ or $$(\bigcap\limits_{i\in\mathfrak{p}_1}C(L_i,\Xi_i^+,\Xi_i^-))\cap(\bigcap\limits_{i=1}^mC(L_i,\Phi_{p,i}^+,\Phi_{p,i}^-))=\{Z'\}.$$

Suppose that $Z'$ is relative interior point of $\bigcap\limits_{i\in \mathfrak{p}_1}C(L_i,\Phi_{p,i}^+,\Phi_{p,i}^-)$.
Then both subsets $\bigcap\limits_{i\in\mathfrak{p}_1}C(L_i,\Xi_i^+,\Xi_i^-) $ and $ \bigcap\limits_{i\in \mathfrak{p}_1}C(L_i,\Phi_{p,i}^+,\Phi_{p,i}^-)$  have common relative interior point. Hence all polyhedra:
$C(L_i,\Xi_i^+,\Xi_i^-)$
and
$C(L_i,\Phi_{p,i}^+,\Phi_{p,i}^-)$ for $i\in\mathfrak{p}_1$
have common relative interior point.
Therefore
$C(L_i,\Xi_i^+,\Xi_i^-)=C(L_i,\Phi_{p,i}^+,\Phi_{p,i}^-)$ for $i\in\mathfrak{p}_1$ because for any $i\in\mathfrak{p}_1$ we have $\Phi_{p,i}^+\cup\Phi_{p,i}^-=\Xi_i^+\cup\Xi_i^-$ and $|\Phi_{p,i}^+\cup\Phi_{p,i}^-|=|\Xi_i^+\cup\Xi_i^-|\leq d+1-\dim L_i$.
Then $\bigcap\limits_{i\in\mathfrak{p}_1}C(L_i,\Xi_i^+,\Xi_i^-) = \bigcap\limits_{i\in \mathfrak{p}_1}C(L_i,\Phi_{p,i}^+,\Phi_{p,i}^-)$.
A contradiction because $Z'$ is not relative interior point of $\bigcap\limits_{i\in\mathfrak{p}_1}C(L_i,\Xi_i^+,\Xi_i^-)$.
Hence $f(Z') \leq \min\limits_{\bigcap\limits_{i\in \mathfrak{p}_1}C(L_i,\Phi_{p,i}^+,\Phi_{p,i}^-)}f \leq \min\limits_{\bigcap\limits_{i=1}^mC(L_i,\Phi_{p,i}^+,\Phi_{p,i}^-)}f=f(Z_p)$.

Let us consider a subcases:

{\it Case 2.1.1: $\sum\limits_{i=1}^m|\Psi_{p,i}^-|=t$.}
Then
$$0<\min\bigcap\limits_{i\in\mathfrak{p}_1}C(L_i,\Psi_{p,i}^+,\Psi_{p,i}^-)=\max\bigcap\limits_{i\in\mathfrak{p}_1}C(L_i,\Psi_{p,i}^+,\Psi_{p,i}^-)<f(Z)$$
because $\bigcap\limits_{i\in\mathfrak{p}_1}C(L_i,\Psi_{p,i}^+,\Psi_{p,i}^-)=\{Z_p\}$.
The dimension of intersection $\bigcap\limits_{i\in\mathfrak{p}_1}C(L_i,\Psi_{p,i}^+,\Psi_{p,i}^-)=\{Z_p\}$ equals 0.
Since $\Psi_{p,1}^-, ..., \Psi_{p,m}^-\subset A$, take $\Phi_1^+=\Psi_{p,1}^+, ..., \Phi_m^+=\Psi_{p,m}^+, \Phi_1^-=\Psi_{p,1}^-, ..., \Phi_m^-=\Psi_{p,m}^-$ as required subsets.

{\it Case 2.1.2: $\sum\limits_{i=1}^m|\Psi_{p,i}^-|<t$.}
Take $\Phi_{p+1,i}^+=\Phi_{p,i}^+$, $\Phi_{p+1,i}^-=\Phi_{p,i}^-$ for $i\neq j_1$.
If $T_p\in\Phi_{p,j_1}^+$ take $\Phi_{p+1,j_1}^+=\Psi_{p,j_1}^+$ and $\Phi_{p+1,j_1}^-=\Psi_{p,j_1}^-\cup\{T_p\}$.
If $T_p\in\Phi_{p,j_1}^-$ take $\Phi_{p+1,j_1}^+=\Psi_{p,j_1}^+\cup\{T_p\}$ and $\Phi_{p+1,j_1}^-=\Psi_{p,j_1}^-$.
Then the polyhedra $C(L_{j_1},\Phi_{p,j_1}^+,\Phi_{p,j_1}^-)$ and $C(L_{j_1},\Phi_{p+1,j_1}^+,\Phi_{p+1,j_1}^-)$ have common face:  $C(L_{j_1},\Psi_{p,j_1}^+,\Psi_{p,j_1}^-)$.
 Therefore the interestions
 $$\bigcap\limits_{i\in\mathfrak{p}_1}C(L_i,\Phi_{p,i}^+,\Phi_{p,i}^-)\subset l$$
 and
 $$\bigcap\limits_{i\in\mathfrak{p}_1}C(L_i,\Phi_{p+1,i}^+,\Phi_{p+1,i}^-)\subset l$$
 have common endpoint $Z_p$.
 We have

$\bullet$ $\Phi_{p+1,1}^-, ..., \Phi_{p+1,m}^-\subset A$;

 $\bullet$ $\dim\bigcap\limits_{i=1}^mC(L_i,\Phi_{p+1,i}^+,\Phi_{p+1,i}^-)=1$;

 $\bullet$ $\Phi_{p+1,i}^+\cup \Phi_{p+1,i}^-=\Gamma_{i}^+\cup \Gamma_{i}^-$ for any $i\in \mathfrak{p}_1$;

 $\bullet$ $\max\limits_{\bigcap\limits_{i=1}^mC(L_i,\Phi_{p+1,i}^+,\Phi_{p+1,i}^-)}f=f(Z_p)$;

  $\bullet$ $\sum\limits_{i=1}^m|\Phi_{p+1,i}^-|\leq 1+\sum\limits_{i=1}^m|\Psi_{p,i}^-|\leq t$

 $\bullet$ $\max\limits_{\bigcap\limits_{i\in\mathfrak{p}_1}C(L_i,\Xi_i^+,\Xi_i^-)}f \leq  \max\limits_{\bigcap\limits_{i=1}^mC(L_i,\Phi_{p+1,i}^+,\Phi_{p+1,i}^-)}f <\max\limits_{\bigcap\limits_{i=1}^mC(L_i,\Phi_{p,i}^+,\Phi_{p,i}^-)}f\leq f(Z)$.

  The system of pairwise disjoint subsets $\{\Phi_{p+1,1}^+, ..., \Phi_{p+1,m}^+, \Phi_{p+1,1}^-, ..., \Phi_{p+1,m}^-\}$ of $A\cup A_{higher}$ is constructed.

{\it Case 2.2: $f(Z') = f(Z_{p-1})$.}

 Then $Z'=Z_{p-1}$.
 %Since $f(Z)>f(Z')$, the point $Z_{p-1}$ is not $Z=Z_0$.
 %Hence $p\geq 2$.
 %$\bullet$ $Z'$ belong to proper face of $C(L_{j_2},\Phi_{p,j_2}^+,\Phi_{p,j_2}^-)$;
  The point $Z_{p-1}$ is relative interior point of $C(L_i,\Psi_{p-1,i}^+,\Psi_{p-1,i}^-)$ for $j\in\{1, ..., n\}$.
 The point $Z'$ is relative interior point of $C(L_i,\Xi_{i}^{+*},\Xi_{i}^{-*})$ for $j\in\mathfrak{p}_1$.

  Then for $i\in\mathfrak{p}_1$ the polyhedra $C(L_i,\Psi_{p-1,i}^+,\Psi_{p-1,i}^-)$ and $C(L_i,\Xi_{i}^{+*},\Xi_{i}^{-*})$ have common relative interior point.
  Hence $\Psi_{p-1,i}^+=\Xi_{i}^{+*}$, $\Psi_{p-1,i}^-=\Xi_{i}^{-*}$ for any $i\in\mathfrak{p}_1$ because $\Xi_{i}^{+}\cup\Xi_{i}^{-}=\Phi_{p-1,i}^+\cup\Phi_{p-1,i}^-$ and $|\Phi_{p,i}^+\cup\Phi_{p,i}^-|\leq d+1-\dim L_i$.

 Therefore
 \begin{multline*}
  (\sum\limits_{i\in\mathfrak{p}_2}\lfloor\frac{d-\dim L_i}{2}\rfloor) -1+\sum\limits_{i\in\mathfrak{p}_1}|\Xi_{i}^{-}|\geq t\geq\\
  \geq\sum\limits_{i=1}^m|\Phi_{p,i}^-|\geq \sum\limits_{i\in\mathfrak{p}_2}|\Phi_{p,i}^-|+\sum\limits_{i\in\mathfrak{p}_1}|\Xi_{i}^{-*}|\geq \sum\limits_{i\in\mathfrak{p}_2}|\Phi_{p,i}^-|-1+ \sum\limits_{i\in\mathfrak{p}_1}|\Xi_{i}^{-}|.
 \end{multline*}

  From Lemma~\ref{add_negative_mass_to_full_dimensional_simplexes} follows that there is a family $\{\Phi_i^-\}_{i\in \mathfrak{p}_2}\cup\{\Phi_i^+\}_{i\in \mathfrak{p}_2}$ pairwise disjoint subsets of $ \{T_{p-1}\}\cup\bigcup\limits_{i\in \mathfrak{p}_2}(\Phi_{p,i}^{+}\cup \Phi_{p,i}^{-}) $ such that
  $$\sum\limits_{i\in \mathfrak{p}_2}|\Phi_i^-|=t+1- \sum\limits_{i\in\mathfrak{p}_1}|\Xi_{i}^{-}|$$ and for any $i\in \mathfrak{p}_2$:

 $\bullet$ $Z'=Z_p\in C(L_i,\Phi_i^+,\Phi_i^-)$;

 $\bullet$ $\Phi_i^-\subset \R^d\times0$;

 $\bullet$ $|\Phi_i^+|+|\Phi_i^-|=d+1-\dim L_i$.

 Take for $i\in\mathfrak{p}_1$:
 $\Phi_i^+=\Xi_i^{+*}$ and $\Phi_i^-=\Xi_i^{-*}$.
Then:

$\bullet$ $\Phi_1^-,...,\Phi_m^-\subset A$;

$\bullet$ $f(Z)> \min\limits_{\bigcap\limits_{i=1}^mC(L_i,\Phi_i^+,\Phi_i^-)}f = f(Z') =\max\limits_{\bigcap\limits_{i=1}^mC(L_i,\Phi_i^+,\Phi_i^-)}f >0$;

$\bullet$  $\sum\limits_{i=1}^m|\Phi_i^-|=t+1- \sum\limits_{i\in\mathfrak{p}_1}|\Xi_{i}^{-}|+\sum\limits_{i\in\mathfrak{p}_1}|\Xi_{i}^{-*}|\in \{t;t+1\}$;

$\bullet$  $\dim\bigcap\limits_{i=1}^mC(L_i,\Phi_i^+,\Phi_i^-)=\dim\{Z'\}=0$.

Therefore $\Phi_1^+, ..., \Phi_m^+, \Phi_1^-, ...,  \Phi_m^-$ is required subsets.

\smallskip

In all cases we construct either: required subsets $\Phi_1^+, ..., \Phi_m^+, \Phi_1^-, ...,  \Phi_m^-$ or the subsets
$\Phi_{r,1}^+, ..., \Phi_{r,m}^+, \Phi_{r,1}^-, ..., \Phi_{r,m}^-$.
There are only finite ways to take $2m$ pairwise disjoint subsets of the set $A\cup A_{higher}$.

Since $\max\limits_{\bigcap\limits_{i=1}^mC(L_i,\Phi_{r,i}^+,\Phi_{r,i}^-)}f<\max\limits_{\bigcap\limits_{i=1}^mC(L_i,\Phi_{r-1,i}^+,\Phi_{r-1,i}^-)}f$, we can not construct the system of pairwise disjoint subsets $\{\Phi_{r,1}^+, ..., \Phi_{r,m}^+, \Phi_{r,1}^-, ..., \Phi_{r,m}^-\}$ for any positive integer $r$.
Hence required subsets $\Phi_1^+, ..., \Phi_m^+, \Phi_1^-, ...,  \Phi_m^-$ are exists.

Proposition is proved.

\begin{proposition}\label{optimization_on_line_2}
 In the conditions of Proposition~\ref{optimization_on_line} suppose that $f(Z)\leq\min(f(Z^{(t)});f(Z^{(t+1)}))$.
 Then $(\R^d\times0)\cap\bigcap\limits_{i=1}^mC(L_i,\Phi_i^+,\Phi_i^-)\neq\varnothing$ and $\sum\limits_{i=1}^m|\Phi_i^-|=t+1$.
\end{proposition}

{\it Proof.}
If the intersection $\bigcap\limits_{i=1}^mC(L_i,\Phi_i^+,\Phi_i^-)$ is unique point set then the value of the functional $f$ in this point belong to both intervals $(0;f(Z^{(t)}))$ and $(0;f(Z^{(t+1)}))$.
A contradiction with choice of points $Z^{(t-1)}$ and $Z^{(t)}$ because $\sum\limits_{i=1}^m|\Phi_i^-|\in\{t,t+1\}$.
Then $\dim\bigcap\limits_{i=1}^mC(L_i,\Phi_i^+,\Phi_i^-)\neq 0$.
Hence $\sum\limits_{i=1}^m|\Phi_i^-|=t+1$.

If $\min\limits_{\bigcap\limits_{i=1}^mC(L_i,\Phi_i^+,\Phi_i^-)}f>0$ denote by $Z''$ a point of minimum of the functional $f$ on the intersection $\bigcap\limits_{i=1}^mC(L_i,\Phi_i^+,\Phi_i^-)$.
According Lemma~\ref{polyhedron_border} there are a subsets $\Phi_1^{+''}\subset \Phi_1^+, ..., \Phi_m^{+''}\subset\Phi_m^+, \Phi_1^{-''}\subset\Phi_1^-, ...,  \Phi_m^{-''}\subset\Phi_m^-$ such that $\{Z''\}=\bigcap\limits_{i=1}^mC(L_i;\Phi_i^{+''};\Phi_i^{-''})$.
From sufficiently general position follows that
\begin{gather*}
 \sum\limits_{i=1}^m(\dim L_i+|\Phi_i^{+''}|+|\Phi_i^{-''}|)=(d+1)(m-1)=\\
 =-1+\sum\limits_{i=1}^m(\dim L_i+|\Phi_i^{+}|+|\Phi_i^{-}|).
\end{gather*}
Then
\begin{gather*}
 \sum\limits_{i=1}^m(|\Phi_i^{+''}|+|\Phi_i^{-''}|) =-1+\sum\limits_{i=1}^m(|\Phi_i^{+}|+|\Phi_i^{-}|).
\end{gather*}
Hence $\sum\limits_{i=1}^m|\Phi_i^{-''}|\in\{t;t+1\}$.
A contradiction with choice of points $Z^{(t-1)}$ and $Z^{(t)}$ because $f(Z'')<f(Z)\leq\min(f(Z^{(t)});f(Z^{(t+1)}))$.
Then  $\min\limits_{\bigcap\limits_{i=1}^mC(L_i,\Phi_i^+,\Phi_i^-)}f=0$.
Therefore $(\R^d\times0)\cap\bigcap\limits_{i=1}^mC(L_i,\Phi_i^+,\Phi_i^-)\neq\varnothing$.

Proposition is proved.

{\it Proof of Theorem~\ref{main_result} in conditions 1:}

For $k=0$ Theorem proved according to Remark~\ref{Roudneffs_case}.
Next we suppose that $k\geq 1$.
Then the both points $Z^{(k)}$ and $Z^{(k-1)}$ are exists.

Suppose that $f(Z^{(k)})> f(Z^{(k-1)})$.
Take a subset $\mathfrak{p}_0\subset \{1, ..., m\}$:

$$\mathfrak{p}_0=\{i|\dim C(L_i, Y_i^{(k-1)+}, Y_i^{(k-1)-})=d+1\}.$$
Then $|\mathfrak{p}_0|\geq \max(0;m-d-1)$.
Hence
\begin{gather*}
\sum\limits_{i\in \mathfrak{p}_0}\lfloor \frac{d-\dim L_i}{2}\rfloor
\geq \min\limits_{\mathfrak{p}\subset\{1,...,m\}, |\mathfrak{p}|\leq \max(0;m-d-1)}\sum\limits_{i\in \mathfrak{p}}\lfloor \frac{d-\dim L_i}{2}\rfloor\geq k\geq\\
\geq k-\sum\limits_{i\notin \mathfrak{p}_0}|Y_i^{(k-1)-}|>k-1-\sum\limits_{i\notin \mathfrak{p}_0}|Y_i^{(k-1)-}|=\sum\limits_{i\in \mathfrak{p}_0}|Y_i^{(k-1)-}|
\end{gather*}

Apply Lemma~\ref{add_negative_mass_to_full_dimensional_simplexes} to the point $Z^{(k-1)}\in\bigcap\limits_{i\in\mathfrak{p}_0}C(L_i, Y_i^{(k-1)+}, Y_i^{(k-1)-})$.
There is a family $\{P_j\}_{j\in \mathfrak{p}_0}\cup\{Q_j\}_{j\in \mathfrak{p}_0}$ pairwise disjoint subsets of $ \{T^{(k-1)}\}\cup\bigcup\limits_{i\in \mathfrak{p}_0}(Y_i^{(k-1)+}\cup Y_i^{(k-1)-}) $ such that
  $$\sum\limits_{i\in \mathfrak{p}_0}|Q_j|=k-\sum\limits_{i\notin \mathfrak{p}_0}|Y_i^{(k-1)-}|;$$
  $$\varnothing\neq ( \{T^{(k-1)}\}\cup\bigcup\limits_{i\in \mathfrak{p}_0}(Y_i^{(k-1)+}\cup Y_i^{(k-1)-}) )\setminus (\bigcup\limits_{i\in \mathfrak{p}_0}(P_i\cup Q_i))\subset \R^d\times0$$
  and for any $i\in \mathfrak{p}_0$:

 $\bullet$ $Z^{(k-1)}\in C(L_i,P_i,Q_i)$;

 $\bullet$ $Q_i\subset \R^d\times0$;

 $\bullet$ $|P_i|+|Q_i|=d+1-\dim L_i$.

Then $Z^{(k-1)}\in (\bigcap\limits_{i\in \mathfrak{p}_0}C(L_i,P_i,Q_i))\cap\bigcap\limits_{i\notin \mathfrak{p}_0}C(L_i, Y_i^{(k-1)+}, Y_i^{(k-1)-})$ and \newline $\sum\limits_{i\in \mathfrak{p}_0}|Q_j|+\sum\limits_{i\notin \mathfrak{p}_0}|Y_i^{(k-1)-}|=k$.
A contradiction with definition of $Z^{(k)}$ because $f(Z^{(k-1)})<f(Z^{(k)})$.

Hence $f(Z^{(k)})\leq f(Z^{(k-1)})$.
Theorem~\ref{main_result} follows from Proposition~\ref{add_positive mass}.

  {\it Proof of Theorem~\ref{main_result} in conditions 2}

For $k=0$ Theorem proved according to Remark~\ref{Roudneffs_case}.
Next, we suppose that $k\geq 1$.
Then the both points $Z^{(k)}$ and $Z^{(k-1)}$ are exists.

If $f(Z^{(k-1)})\notin (0;f(Z^{(k)}))$ Theorem follows from Proposition~\ref{add_positive mass}.
Next, we suppose that  $f(Z^{(k-1)})\in (0;f(Z^{(k)}))$.

  We have $|\{i|L_i\neq\varnothing\}|+\sum\limits_{i=1}^m|Y_i^{(k-1)-}|= |\{i|L_i\neq\varnothing\}|+k-1\leq m-1$.
  Therefore there is a number $j$ such that $L_j=Y_j^{(k-1)-}=\varnothing$.
  Then $C(L_j; Y_j^{(k-1)+}; Y_j^{(k-1)-})=\left<Y_j^{(k-1)+}\right>$ and $\dim C(L_j; Y_j^{(k-1)+}; Y_j^{(k-1)-})\leq d+1=3$.
  Take $j$ such that $\dim C(L_j; Y_j^{(k-1)+}; Y_j^{(k-1)-})$ is maximal.

  If the point $Z^{(k-1)}$ is not minimal point of the functional $f$ of the polyhedron
   $C(L_j; Y_j^{(k-1)+}; \{T^{(k-1)}\}\cup Y_j^{(k-1)-})\cap\bigcap\limits_{i\in\{1,...,m\}\setminus\{j\}}C(L_i; Y_i^{(k-1)+}; Y_i^{(k-1)-})$
 then Theorem~\ref{main_result} follows from Proposition~\ref{add_negative_mass_2}.
  Next, we suppose that the point $Z^{(k-1)}$ is minimal point of the functional $f$ of this polyhedron.

  Let us consider a subcases.
  
  {\it Case 1: The polyhedron $C(L_j; Y_j^{(k-1)+}; Y_j^{(k-1)-})$ is 3-dimensional.}

 From Lemma~\ref{add_negative_mass_to_full_dimensional_simplex} follows that there is a subsets $\mathfrak{Y}_j^{(k-1)+}, \mathfrak{Y}_j^{(k-1)-}$ of $\{T^{(k-1)}\}\cup Y_j^{(k-1)+}\cup Y_j^{(k-1)-}$ such that:

 $\bullet$ $|\mathfrak{Y}_j^{(k-1)-}|=1$;

 $\bullet$ $\mathfrak{Y}_j^{(k-1)-}\subset A$;

 $\bullet$ $Z^{(k-1)}\in C(L_j;\mathfrak{Y}_j^{(k-1)+}; \mathfrak{Y}_j^{(k-1)-})$.

 Then

 $\bullet$ $|\mathfrak{Y}_j^{(k-1)-}|+\sum\limits_{i\neq j} |Y_j^{(k-1)-}|=k$;

 $\bullet$ $Z^{(k-1)}=C(L_i;\mathfrak{Y}_j^{(k-1)+};\mathfrak{Y}_j^{(k-1)-})\cap\bigcap\limits_{i\in\{1,...,m\}\setminus\{j\}}C(L_i;Y_i^{(k-1)+};Y_i^{(k-1)-})$.

 A contradiction with choice of $Z^{(k)}$ because $f(Z^{(k)})>f(Z^{(k-1)})>0$.

  {\it Case 2: The polyhedron $C(L_j; Y_j^{(k-1)+}; Y_j^{(k-1)-})$ is 2-dimensional.}

Then the polyhedron $C(L_j; \{T^{(k-1)}\}\cup Y_j^{(k-1)+}; Y_j^{(k-1)-})=\left<\{T^{(k-1)}\}\cup Y_j^{(k-1)+}\right>$ is 3-dimensional.
From sufficiently general position follows that the intersection $\bigcap\limits_{i\in\{1,...,m\}\setminus\{j\}}C(L_i; Y_i^{(k-1)+}; Y_i^{(k-1)-})$ is 1-dimensional.
  Denote by $Z_4^{(k-1)}$ the minimal point of the functional $f$ of the polyhedron\newline
   $C(L_j; \{T^{(k-1)}\}\cup Y_j^{(k-1)+}; Y_j^{(k-1)-})\cap\bigcap\limits_{i\in\{1,...,m\}\setminus\{j\}}C(L_i; Y_i^{(k-1)+}; Y_i^{(k-1)-})$.

   Let us consider a subcases.

   {\it Case 2.1: $Z_4^{(k-1)}\in \R^d\times0$.}

   Then:

  $$(\R^d\times0)\cap\left<\{T^{(k-1)}\}\cup Y_j^{(k-1)+}\right>\cap\bigcap\limits_{i\neq j}\aff(L_i\cup Y_i^{(k-1)+}\cup Y_i^{(k-1)-})\neq \varnothing.$$

From Lemma~\ref{add_negative_mass} follows that there is a number $j_{decr}$ such that $Z^{(k-1)}$ is not minimal point of the functional $f$ on the polyhedron
$$C(L_{j_{decr}}; Y_{j_{decr}}^{(k-1)+};\{T^{(k-1)}\}\cup  Y_{j_{decr}}^{(k-1)-})\cap\bigcap\limits_{i\neq j_{decr}}C(L_i; Y_i^{(k-1)+}; Y_i^{(k-1)-}).$$
Apply Proposition~\ref{add_negative_mass_2}.
Theorem~\ref{main_result} proved in this case.

  {\it Case 2.2: $Z_4^{(k-1)}$ is interior point of $\left<\{T^{(k-1)}\}\cup Y_j^{(k-1)+}\right>$.}

  From Lemma~\ref{polyhedron_border} follows that there are subsets $Y_i^{(k-1)+''}\subset Y_i^{(k-1)+}$, $Y_i^{(k-1)-''}\subset Y_i^{(k-1)-}$ such that $\{Z_4^{(k-1)}\}=\bigcap\limits_{i\neq j}C(L_i; Y_i^{(k-1)+''}; Y_i^{(k-1)-''})$.
 Then
 \begin{gather*}
  \sum\limits_{i\neq j}(|Y_i^{(k-1)+''}|+|Y_i^{(k-1)-''}|)=(d+1)(m-2)-\sum\limits_{i\neq j}\dim L_i=\\
  =-1+\sum\limits_{i\neq j}(|Y_i^{(k-1)+}|+|Y_i^{(k-1)-}|)
 \end{gather*}
because $\dim \bigcap\limits_{i\neq j}C(L_i; Y_i^{(k-1)+}; Y_i^{(k-1)-})=1$.
Therefore either:
$$\sum\limits_{i\neq j}|Y_i^{(k-1)-''}|=\sum\limits_{i\neq j}|Y_i^{(k-1)-}|=k-1$$
or
$$\sum\limits_{i\neq j}|Y_i^{(k-1)-''}|=-1+\sum\limits_{i\neq j}|Y_i^{(k-1)-}|=k-2.$$

  According to Lemma~\ref{add_negative_mass_to_full_dimensional_simplex} there are pairwise disjoint subsets $\mathfrak{Y}_j^+, \mathfrak{Y}_j^-\subset \{T^{(k-1)}\}\cup Y_j^{(k-1)+}$ such that:

  $\bullet$ $Z_4^{(k-1)}\in C(L_j, \mathfrak{Y}_j^+, \mathfrak{Y}_j^-)$;

  $\bullet$ $|\mathfrak{Y}_j^-|=1$;

  $\bullet$ $\mathfrak{Y}_j^-\subset \R^d\times0$;

  $\bullet$ $|\mathfrak{Y}_j^+|+|\mathfrak{Y}_j^-|=d+1-\dim L_j=3$.

  Hence
  $\{Z_4^{(k-1)}\}=C(L_j, \mathfrak{Y}_j^+, \mathfrak{Y}_j^-)\cap\bigcap\limits_{i\neq j}C(L_i; Y_i^{(k-1)+''}; Y_i^{(k-1)-''})$ and
$|\mathfrak{Y}_j^-|+\sum\limits_{i\neq j}|Y_i^{(k-1)-''}|\in \{k-1,k\}$.

A contradiction because $f(Z_4^{(k-1)})<f(Z^{(k-1)})<f(Z^{(k)})$.

{\it Case 2.3: $Z_4^{(k-1)}$ is not interior point of $\left<\{T^{(k-1)}\}\cup Y_j^{(k-1)+}\right>$ and  $f(Z_4^{(k-1)})>0$.}

Then the point $Z_4^{(k-1)}$ belong to some proper face of $\left<\{T^{(k-1)}\}\cup Y_j^{(k-1)+}\right>$.
There is a subset $Y_j^{(k-1)+''}\subset \{T^{(k-1)}\}\cup Y_j^{(k-1)+}$ such that the point $Z_4^{(k-1)}$ belong to relative interior of $\left<Y_j^{(k-1)+''}\right>$.
From sufficiently general position follows $\dim \left<Y_j^{(k-1)+''}\right>=2$.
Then $(\{T^{(k-1)}\}\cup Y_j^{(k-1)+})\setminus Y_j^{(k-1)+''}$ is an unique point set.
Denote this point by $T^{(k-1)''}$

The point $Z_4^{(k-1)}$ is minimum point of functional $f$ on the intersection
\begin{gather*}
 \left<\{T^{(k-1)}\}\cup Y_j^{(k-1)+}\right>\cap\bigcap\limits_{i\neq j}C(L_i; Y_i^{(k-1)+}; Y_i^{(k-1)-})=\\
 =\left<\{T^{(k-1)''}\}\cup Y_j^{(k-1)+''}\right>\cap\bigcap\limits_{i\neq j}C(L_i; Y_i^{(k-1)+}; Y_i^{(k-1)-}).
\end{gather*}
Then the point $Z_4^{(k-1)}$ is not a minimum point of functional $f$ on the intersection
\begin{gather*}
C(\varnothing; Y_j^{(k-1)+''};\{T^{(k-1)''}\})\cap\bigcap\limits_{i\neq j}C(L_i; Y_i^{(k-1)+}; Y_i^{(k-1)-}).
\end{gather*}
Therefore Theorem~\ref{main_result} implies from Proposition~\ref{add_negative_mass_2}.

  {\it Case 3: The polyhedron $C(L_j; Y_j^{(k-1)+}; Y_j^{(k-1)-})$ is 1-dimensional.}

 From sufficiently general position follows that the intersection $\bigcap\limits_{i\neq j}C(L_i; Y_i^{(k-1)+}; Y_i^{(k-1)-})$ is 2-dimensional.
 There is a number $j_1$ such that the dimension of polyhedron $C(L_{j_1}; Y_{j_1}^{(k-1)+}; Y_{j_1}^{(k-1)-})$ is equal 2.
 For any $j_2\notin\{j,j_1\}$ the dimension of polyhedron $C(L_{j_2}; Y_{j_2}^{(k-1)+}; Y_{j_2}^{(k-1)-})$ is equal 3.
 According to the choice of number $j$ we have $L_i\cup Y_i^{(k-1)-}\neq \varnothing$ for $i\neq j$.
 Hence
 \begin{gather*}
m-1\leq|\{i|L_i\neq\varnothing\}|+\sum\limits_{i\neq j}|Y_i^{(k-1)-}|=\\
=|\{i|L_i\neq\varnothing\}|+\sum\limits_{i=1}^m|Y_i^{(k-1)-}|
=|\{i|L_i\neq\varnothing\}|+k-1\leq m-1.
 \end{gather*}
 Therefore for any $i\in \{1,...,m\}\setminus\{j\}$:

 $\bullet$ either set $L_i$ or $Y_i^{(k-1)-}$ is empty;

 $\bullet$ $|Y_i^{(k-1)-}|\leq1$.

 Similarly to Case 2 denote by $Z_4^{(k-1)}$ the minimal point of the functional $f$ of the polyhedron\newline
   $C(L_j; \{T^{(k-1)}\}\cup Y_j^{(k-1)+}; Y_j^{(k-1)-})\cap\bigcap\limits_{i\in\{1,...,m\}\setminus\{j\}}C(L_i; Y_i^{(k-1)+}; Y_i^{(k-1)-})$.

   Let us consider a subcases.

  {\it Case 3.1: $$(\R^d\times0)\cap\left<\{T^{(k-1)}\}\cup Y_j^{(k-1)+}\right>\cap\bigcap\limits_{i\neq j}\aff(L_i\cup Y_i^{(k-1)+}\cup Y_i^{(k-1)-})\neq \varnothing.$$}

From Lemma~\ref{add_negative_mass} follows that there is a number $j_{decr}$ such that $Z^{(k-1)}$ is not minimal point of the functional $f$ on the polyhedron
$$C(L_{j_{decr}}; Y_{j_{decr}}^{(k-1)+};\{T^{(k-1)}\}\cup  Y_{j_{decr}}^{(k-1)-})\cap\bigcap\limits_{i\neq j_{decr}}C(L_i; Y_i^{(k-1)+}; Y_i^{(k-1)-}).$$
Apply Proposition~\ref{add_negative_mass_2}.
Theorem~\ref{main_result} proved in this case.

{\it Case 3.2: $$(\R^d\times0)\cap\left<\{T^{(k-1)}\}\cup Y_j^{(k-1)+}\right>\cap\bigcap\limits_{i\neq j}\aff(L_i\cup Y_i^{(k-1)+}\cup Y_i^{(k-1)-})= \varnothing.$$}

Then $f(Z_4^{(k-1)})>0$.
Denote by

$\bullet$ $T^{(k-1)''}$ the unique point of intersection $Y_j^{(k-1)+}\cap(\R^d\times0)$;

$\bullet$ $Y_j^{(k-1)+''}=(Y_j^{(k-1)+}\setminus\{T^{(k-1)''}\})\cup\{T^{(k-1)}\}$.

Then
$$(\R^d\times0)\cap C(L_j;Y_j^{(k-1)+''};T^{(k-1)''})\cap\bigcap\limits_{i\neq j}\aff(L_i\cup Y_i^{(k-1)+}\cup Y_i^{(k-1)-})\neq \varnothing.$$
Hence
$$(\R^d\times0)\cap C(L_j;Y_j^{(k-1)+''};T^{(k-1)''})\cap \aff(L_{j_1}\cup Y_{j_1}^{(k-1)+}\cup Y_{j_1}^{(k-1)-})\neq \varnothing.$$
Therefore there are pairwise disjoint subsets $\Xi_{j_1}^{+''},\Xi_{j_1}^{-''}\subset Y_{j_1}^{(k-1)+}\cup Y_{j_1}^{(k-1)-}$ such that $C(L_{j_1};\Xi_{j_1}^{+''};\Xi_{j_1}^{-''})\cap C(L_j;Y_j^{(k-1)+''};T^{(k-1)''})\neq \varnothing$.

Let us consider a subcases.

{\it Case 3.2.1: $|\Xi_{j_1}^{-''}|\geq|Y_{j_1}^{(k-1)-}|$.}

Take

$\bullet$ $\Xi_{j_1}^{+}=\Xi_{j_1}^{+''}$, $\Xi_{j_1}^{-}=\Xi_{j_1}^{-''}$, $\Xi_{j}^{+}=Y_j^{(k-1)+''}$, $\Xi_{j}^{-}=\{T^{(k-1)''}\}$;

$\bullet$ $\Gamma_j^+= Y_j^{(k-1)+}$,  $\Gamma_j^-=\{T^{(k-1)}\}$;

$\bullet$ $\Gamma_i^+= Y_i^{(k-1)+}$,  $\Gamma_i^-= Y_i^{(k-1)-}$ for $i\neq j$;

$\bullet$ $\mathfrak{p}_1=\{j,j_1\}$; $\mathfrak{p}_2=\{1,...,m\}\setminus\{j,j_1\}$.

Then
\begin{gather*}
 \sum\limits_{i=1}^m|\Gamma_i^-|= k = 1+\sum\limits_{i=1}^m|Y_i^{(k-1)-}|\leq \sum\limits_{i\in\mathfrak{p}_1}|\Xi_i^-|+\sum\limits_{i\in\mathfrak{p}_2}|Y_i^{(k-1)-}|\leq\\
 \leq \sum\limits_{i\in\mathfrak{p}_1}|\Xi_i^-|+|\{i|i\in\mathfrak{p}_2, L_i = \varnothing\}|= \sum\limits_{i\in\mathfrak{p}_1}|\Xi_i^-|+\sum\limits_{i\in\mathfrak{p}_2}\lfloor \frac{d-\dim L_i}{2} \rfloor
\end{gather*}

Apply Proposition~\ref{optimization_on_line} for these $\Gamma_1^+, ..., \Gamma_m^+, \Gamma_1^-, ..., \Gamma_m^-$ and $t=k-1$.
We get a pairwise disjoint subsets $$\Phi_1^+, ..., \Phi_m^+, \Phi_1^-, ..., \Phi_m^-\subset \bigcup\limits_{i=1}^m(\Gamma_i^+\cup\Gamma_i^-)$$ such that:

$\bullet$ $\Phi_1^-,...,\Phi_m^-\subset A$;

$\bullet$ $\max\limits_{\bigcap\limits_{i=1}^mC(L_i,\Phi_i^+,\Phi_i^-)}f >0$;

$\bullet$ $\min\limits_{\bigcap\limits_{i=1}^mC(L_i,\Phi_i^+,\Phi_i^-)}f < f(Z^{(k-1)})$;

$\bullet$  $\sum\limits_{i=1}^m|\Phi_i^-|\in\{k-1,k\}$;

$\bullet$ either   1) $\sum\limits_{i=1}^m|\Phi_i^-|=k$ or 2) $\dim\bigcap\limits_{i=1}^mC(L_i,\Phi_i^+,\Phi_i^-) =0$.

Since $f(Z^{(k-1)})\leq \min(f(Z^{(k-1)});f(Z^{(k)}))$, we can apply Proposition~\ref{optimization_on_line_2}.
Then $\sum\limits_{i=1}^m|\Phi_i^-|=k$ and $(\R^d\times0)\cap\bigcap\limits_{i=1}^mC(L_i,\Phi_i^+,\Phi_i^-)\neq\varnothing$.

Take:

$\bullet$ $A_i=\Phi_i^+\cap A$ for $i=1,...,m-1$;

$\bullet$ $B_i=\Phi_i^-$ for $i=1,...,m$;

$\bullet$ $A_m=A\setminus((\bigcup\limits_{i=1}^{m-1}A_i)\cup\bigcup\limits_{i=1}^mB_i)$.

Then

 $\bullet$ $\sum\limits_{i=1}^m|B_i|=k$;

  $\bullet$ $\bigcap\limits_{i=1}^mC(L_i, A_i, B_i)\neq\varnothing$.

  Theorem proved in this case.

{\it Case 3.2.2: $|\Xi_{j_1}^{-''}|<|Y_{j_1}^{(k-1)-}|$.}

Then $|Y_{j_1}^{(k-1)-}|\geq 1$.
Hence  $Y_{j_1}^{(k-1)-}$ is unique point set and $L_{j_1}=\varnothing$.
Since $|\Xi_{j_1}^{-''}|<|Y_{j_1}^{(k-1)-}|$, the set $\Xi_{j_1}^{-''}$ is empty and
$$(\R^d\times0)\cap C(L_j;Y_j^{(k-1)+''};T^{(k-1)''})\cap \left< Y_{j_1}^{(k-1)+}\cup Y_{j_1}^{(k-1)-}\right>\neq \varnothing.$$

Denote:

$\bullet$ $D_{j_1}=Y_{j_1}^{(k-1)+}\setminus(\R^d\times0)$;

$\bullet$ $y_{j_1}=Y_{j_1}^{(k-1)+}\cap(\R^d\times0)$.

The set $D_{j_1}$ is not empty because $Z^{(k-1)}\in C(L_j;Y_{j_1}^{(k-1)+};Y_{j_1}^{(k-1)-})$.
The set $y_{j_1}$ is not empty because $(\R^d\times (0;1) )\cap C(L_j;Y_{j_1}^{(k-1)+};Y_{j_1}^{(k-1)-})\neq\varnothing$.
Since $\dim C(L_j;Y_{j_1}^{(k-1)+};Y_{j_1}^{(k-1)-})=2$, we have   $|Y_{j_1}^{(k-1)+}|=2$.
Hence $|D_{j_1}|=|y_{j_1}|=1$.

The segment $\left< y_{j_1} \cup Y_{j_1}^{(k-1)-}\right>$ equal to intersection $\left< Y_{j_1}^{(k-1)+}\cup Y_{j_1}^{(k-1)-}\right>\cap(\R^d\times0)$.
The ray $C(\varnothing; T^{(k-1)}; T^{(k-1)''})$ equal to the intersection $C(L_j;Y_j^{(k-1)+''};T^{(k-1)''})\cap(\R^d\times0)$.

Therefore the ray $C(\varnothing; T^{(k-1)}; T^{(k-1)''})$ intersect the segment $\left< y_{j_1} \cup Y_{j_1}^{(k-1)-}\right>$.
Then $T^{(k-1)''}\in C(\varnothing; T^{(k-1)}; y_{j_1} \cup Y_{j_1}^{(k-1)-})$.
Hence $$\left<Y_{j}^{(k-1)+}\right>\cap C(\varnothing; T^{(k-1)} \cup D_{j_1}; y_{j_1} \cup Y_{j_1}^{(k-1)-})$$ is a segment with endpoint $T^{(k-1)''}\in\R^d\times0$.
Take:

$\bullet$ $\Gamma_{j_1}^+=Y_{j_1}^{(k-1)+}$, $\Gamma_{j_1}^-=T^{(k-1)} \cup Y_{j_1}^{(k-1)-}$, $\Gamma_{j_1}^{-*}= Y_{j_1}^{(k-1)-}$;

$\bullet$ $\Gamma_{i}^+=Y_{i}^{(k-1)+}$, $\Gamma_{i}^-=Y_{i}^{(k-1)-}$ for $i\neq j_1$;

$\bullet$ $\Xi_{j_1}^+=T^{(k-1)} \cup D_{j_1}$, $\Xi_{j_1}^-=y_{j_1} \cup Y_{j_1}^{(k-1)-}$;

$\bullet$ $\Xi_{j}^+=Y_j^{(k-1)+}$, $\Xi_{j}^-=Y_j^{(k-1)-}=\varnothing$;

$\bullet$ $\mathfrak{p}_1=\{j,j_1\}$; $\mathfrak{p}_2=\{1,...,m\}\setminus\{j,j_1\}$.

Then $$\dim\bigcap\limits_{i=1}^m C(L_i;\Gamma_{i}^+;\Gamma_{i}^-)=1=\dim\bigcap\limits_{i\in\mathfrak{p}_1} C(L_i;\Gamma_{i}^+;\Gamma_{i}^-)$$
and $\Gamma_i^+\cup \Gamma_i^-=\Xi_i^+\cup \Xi_i^-$ for $i\in\mathfrak{p}_1$.
The point $Z^{(k-1)}$ is minimal point of the functional $f$ on the polyhedron $\bigcap\limits_{i=1}^m C(L_i;\Gamma_{i}^+;\Gamma_{i}^-)$.
The maximum of the functional $f$ on the intersection $\bigcap\limits_{i\in\mathfrak{p}_1}C(L_i;\Xi_i^+;\Xi_i^-)$ is not bigger that $f(Z^{(k-1)})$
because $$C(\varnothing; T^{(k-1)} \cup D_{j_1}; y_{j_1} \cup Y_{j_1}^{(k-1)-})\neq C(L_{j_1};Y_{j_1}^{(k-1)+};T^{(k-1)} \cup Y_{j_1}^{(k-1)-}).$$
Since $|\Xi_j^-|+|\Xi_{j_1}^-|=0+2=2$, we have
\begin{gather*}
(\sum\limits_{i\in\mathfrak{p}_2}\lfloor \frac{d-\dim L_i}{2}\rfloor)+ \sum\limits_{i\in\mathfrak{p}_1}|\Xi_i^-|\geq |\{i|i\in\mathfrak{p}_2, L_i=\varnothing\}|+2\geq\\
\geq  |\{i| L_i=\varnothing\}|\geq k = \sum\limits_{i=1}^m|\Gamma_i^-|.
\end{gather*}

Apply Proposition~\ref{optimization_on_line} for these $\Gamma_1^+, ..., \Gamma_m^+, \Gamma_1^-, ..., \Gamma_m^-$ and $t=k-1$.
We get a pairwise disjoint subsets $$\Phi_1^+, ..., \Phi_m^+, \Phi_1^-, ..., \Phi_m^-\subset \bigcup\limits_{i=1}^m(\Gamma_i^+\cup\Gamma_i^-)$$ such that:

$\bullet$ $\Phi_1^-,...,\Phi_m^-\subset A$;

$\bullet$ $\max\limits_{\bigcap\limits_{i=1}^mC(L_i,\Phi_i^+,\Phi_i^-)}f >0$;

$\bullet$ $\min\limits_{\bigcap\limits_{i=1}^mC(L_i,\Phi_i^+,\Phi_i^-)}f < f(Z^{(k-1)})$;

$\bullet$  $\sum\limits_{i=1}^m|\Phi_i^-|\in\{k-1,k\}$;

$\bullet$ either   1) $\sum\limits_{i=1}^m|\Phi_i^-|=k$ or 2) $\dim\bigcap\limits_{i=1}^mC(L_i,\Phi_i^+,\Phi_i^-) =0$.

Since $f(Z^{(k-1)})\leq \min(f(Z^{(k-1)});f(Z^{(k)}))$, we can apply Proposition~\ref{optimization_on_line_2}.
Then $\sum\limits_{i=1}^m|\Phi_i^-|=k$ and $(\R^d\times0)\cap\bigcap\limits_{i=1}^mC(L_i,\Phi_i^+,\Phi_i^-)\neq\varnothing$.

Take:

$\bullet$ $A_i=\Phi_i^+\cap A$ for $i=1,...,m-1$;

$\bullet$ $B_i=\Phi_i^-$ for $i=1,...,m$;

$\bullet$ $A_m=A\setminus((\bigcup\limits_{i=1}^{m-1}A_i)\cup\bigcup\limits_{i=1}^mB_i)$.

Then

 $\bullet$ $\sum\limits_{i=1}^m|B_i|=k$;

  $\bullet$ $\bigcap\limits_{i=1}^mC(L_i, A_i, B_i)\neq\varnothing$.

  Theorem proved in this case.

  {\it Case 4: The polyhedron $C(L_j; Y_j^{(k-1)+}; Y_j^{(k-1)-})$ is 0-dimensional.}
  Then the set $Y_j^{(k-1)+}$ is unique point.
  Hence
  $$Y_j^{(k-1)+}=C(L_j; Y_j^{(k-1)+}; Y_j^{(k-1)-})=\bigcap\limits_{i=1}^mC(L_i; Y_i^{(k-1)+}; Y_i^{(k-1)-})=\{Z^{(k-1)}\}.$$
  Since $Y_j^{(k-1)+}\subset A\cup\{W_1^{(0)}, ..., W_m^{(0)}, ..., W_1^{(k)}, ..., W_m^{(k)}\}$ and $Z^{(k-1)}\notin \R^d\times0$, we have $Z^{(k-1)}\in \{W_1^{(0)}, ..., W_m^{(0)}, ..., W_1^{(k)}, ..., W_m^{(k)}\}$.
  A contradiction because $f(Z^{(k-1)})<1<f(W_i^{(q)})$ for any $i=1,... , m$; $q=0, ... ,k$.

Theorem~\ref{main_result} proved in conditions 2.

  {\it Proof of Theorem~\ref{main_result} in conditions 3}

If $k\leq \min\limits_{\mathfrak{p}\subset\{1,...,m\}, |\mathfrak{p}|\leq \max(0;m-d-1)}\sum\limits_{i\in \mathfrak{p}}\lfloor \frac{d-\dim L_i}{2}\rfloor$ then Theorem proved above.
Next, we suppose that $k=1+ \min\limits_{\mathfrak{p}\subset\{1,...,m\}, |\mathfrak{p}|\leq \max(0;m-d-1)}\sum\limits_{i\in \mathfrak{p}}\lfloor \frac{d-\dim L_i}{2}\rfloor$.
Then the both points $Z^{(k)}$ and $Z^{(k-1)}$ are exists.

If $f(Z^{(k-1)})\notin (0;f(Z^{(k)}))$ Theorem follows from Proposition~\ref{add_positive mass}.
Next, we suppose that  $f(Z^{(k-1)})\in (0;f(Z^{(k)}))$.

If $Z^{(k-1)}$ is not a minimum point of the functional $f$ on the intersection $C(L_{t_v};Y_{t_v}^{(k-1)+};\{T^{(k-1)}\}\cup Y_{t_v}^{(k-1)-})\cap\bigcap\limits_{i\neq t_v}C(L_{i};Y_{i}^{(k-1)+}; Y_{i}^{(k-1)-})$ then Theorem follows from Proposition~\ref{add_negative_mass_2}.
Next, we suppose that the point $Z^{(k-1)}$ is minimum point of the functional $f$ on the intersection $$C(L_{t_v};Y_{t_v}^{(k-1)+};\{T^{(k-1)}\}\cup Y_{t_v}^{(k-1)-})\cap\bigcap\limits_{i\neq t_v}C(L_{i};Y_{i}^{(k-1)+}; Y_{i}^{(k-1)-}).$$
Then the point $Z^{(k-1)}$ is not a minimum point of the functional $f$ on the intersection $$C(L_{t_v};\{T^{(k-1)}\}\cup Y_{t_v}^{(k-1)+};Y_{t_v}^{(k-1)-})\cap\bigcap\limits_{i\neq t_v}C(L_{i};Y_{i}^{(k-1)+}; Y_{i}^{(k-1)-}).$$

Let us consider a subcases.

{\it Case 1:
$$(\R^d\times0)\cap C(L_{t_v};\{T^{(k-1)}\}\cup Y_{t_v}^{(k-1)+};Y_{t_v}^{(k-1)-})\cap\bigcap\limits_{i\neq t_v}\aff(L_{i}\cup Y_{i}^{(k-1)+}\cup Y_{i}^{(k-1)-})\neq\varnothing.$$ }

Apply Lemma~\ref{add_negative_mass}.
There is a number $j_{decr}$ such that the minimal point of the functional $f$ on the intersection
  $$C(L_{j_{decr}}, Y_{j_{decr}},Y_{j_{decr}}^{(k-1)}\cup\{T^{(k-1)}\})\cap (\bigcap\limits_{i\neq j_{decr}}C(L_i, Y_i^{(k-1)+},Y_i^{(k-1)-}))$$
is not $Z^{(k-1)}$.
Apply Proposition~\ref{add_negative_mass_2}.
Theorem~\ref{main_result} proved in this case.

{\it Case 2:
$$(\R^d\times0)\cap C(L_{t_v};\{T^{(k-1)}\}\cup Y_{t_v}^{(k-1)+};Y_{t_v}^{(k-1)-})\cap\bigcap\limits_{i\neq t_v}\aff(L_{i}\cup Y_{i}^{(k-1)+}\cup Y_{i}^{(k-1)-})=\varnothing.$$ }

Let us consider a subcases.

{\it Case 2.1: $\dim C(L_{t_v}, Y_{t_v}^{(k-1)+}, Y_{t_v}^{(k-1)-})=d+1$.}

Apply Lemma~\ref{add_negative_mass_to_full_dimensional_simplex}. There are pairwise disjoint subsets $P, Q\subset Y_{t_v}^{(k-1)+}\cup Y_{t_v}^{(k-1)-}\cup \{T^{(k-1)}\}$ such that:

 $\bullet$ $Z^{(k-1)}\in C(L_{t_v},P,Q)$;

 $\bullet$ $|Q|=1$;

 $\bullet$ $Q\subset \R^d\times0$;

 $\bullet$ $|P|+|Q|=d+1-\dim L_{t_v}$;

 $\bullet$ $(Y_{t_v}^{(k-1)+}\cup Y_{t_v}^{(k-1)-}\cup \{T^{(k-1)}\})\setminus(P\cup Q)$ is the one point subset of $\R^d\times0$.
 Then
 $$\{Z^{(k-1)}\}= C(L_{t_v},P,Q) \cap \bigcap\limits_{i\neq t_v}C(L_{i}, Y_{i}^{(k-1)+}, Y_{i}^{(k-1)-})$$
 and
 $$|Q|+\sum\limits_{i\neq t_v}|Y_{i}^{(k-1)-}|=k.$$
 A contradiction because $f(Z^{(k-1)})<f(Z^{(k)})$.

{\it Case 2.2: $\dim C(L_{t_v}, Y_{t_v}^{(k-1)+}, Y_{t_v}^{(k-1)-})<d+1$.}

Then the intersection $$ C(L_{t_v};\{T^{(k-1)}\}\cup Y_{t_v}^{(k-1)+};Y_{t_v}^{(k-1)-})\cap\bigcap\limits_{i\neq t_v}\aff(L_{i}\cup Y_{i}^{(k-1)+}\cup Y_{i}^{(k-1)-})$$ is 1-dimensional.

Denote by $\mathfrak{p}_1$ a subset of $\{1, ..., m\}$:
 $$\mathfrak{p}_1=\{i|\dim C(L_{i}, Y_{i}^{(k-1)+}, Y_{i}^{(k-1)-})<d+1\}.$$
 Then $|\mathfrak{p}_1|\leq d+1$ because $\dim\bigcap\limits_{i\in\mathfrak{p}_1}C(L_{i}, Y_{i}^{(k-1)+}, Y_{i}^{(k-1)-})=\dim\{Z^{(k-1)}\}=1$.
Denote $\mathfrak{p}_2=\{1,...,m\}\setminus\mathfrak{p}_1$.
Then $|\mathfrak{p}_2|\leq \max(0; m-d-1)$.
Hence
\begin{gather*}
 \sum\limits_{i\in \mathfrak{p}_2}\lfloor \frac{d-\dim L_i}{2}\rfloor\geq \min\limits_{\mathfrak{p}\subset\{1,...,m\}, |\mathfrak{p}|\leq \max(0;m-d-1)}\sum\limits_{i\in \mathfrak{p}}\lfloor \frac{d-\dim L_i}{2}\rfloor=k-1.
\end{gather*}

Let us consider a subcases.

{\it Case 2.2.1: $\sum\limits_{i\in\mathfrak{p}_1}|Y_i^{(k-1)-}|\geq 1$.}

Then
\begin{gather*}
 \sum\limits_{i\in\mathfrak{p}_2}|Y_i^{(k-1)-}|=k-1-\sum\limits_{i\in\mathfrak{p}_1}|Y_i^{(k-1)-}|<\\
 <k-\sum\limits_{i\in\mathfrak{p}_1}|Y_i^{(k-1)-}|\leq k-1\leq \sum\limits_{i\in \mathfrak{p}_2}\lfloor \frac{d-\dim L_i}{2}\rfloor.
\end{gather*}

The polyhedra $C(L_i, Y_i^{(k-1)+},Y_i^{(k-1)-})$ for $i\in\mathfrak{p}_2$ are full-dimensional and $Z^{(k-1)}\in\bigcap\limits_{i\in\mathfrak{p}_2}C(L_i, Y_i^{(k-1)+},Y_i^{(k-1)-})$.
%From definition of $Z^{(k-1)}$ follows that
%$$\sum\limits_{i\in\mathfrak{p}_2}|Y_i^{(k-1)-}|=k-1-\sum\limits_{i\in\mathfrak{p}_1}|Y_i^{(k-1)-}|.$$

Apply Lemma~\ref{add_negative_mass_to_full_dimensional_simplexes}.
There is a family $\{P_i\}_{i\in\mathfrak{p}_2}\cup\{Q_i\}_{i\in\mathfrak{p}_2}$ pairwise disjoint subsets of  $ \{T^{(k-1)}\}\cup\bigcup\limits_{i\in\mathfrak{p}_2}(Y_i^{(k-1)+}\cup Y_i^{(k-1)-}) $ such that
  $$\sum\limits_{i\in\mathfrak{p}_2}|Q_i|=k-\sum\limits_{i\in\mathfrak{p}_1}|Y_i^{(k-1)-}|;$$
  $$\varnothing\neq (\{T^{(k-1)}\}\cup\bigcup\limits_{i\in\mathfrak{p}_2}(Y_i^{(k-1)+}\cup Y_i^{(k-1)-}) )\setminus (\bigcup\limits_{i\in\mathfrak{p}_2}(P_i\cup Q_i))\subset \R^d\times0$$
  and for any $i\in\mathfrak{p}_2$:

 $\bullet$ $Z^{(k-1)}\in C(L_i,P_i,Q_i)$;

 $\bullet$ $Q_i\subset \R^d\times0$;

 $\bullet$ $|P_i|+|Q_i|=d+1-\dim L_i$.

 Then

 $$\{Z^{(k-1)}\}=(\bigcap\limits_{i\in\mathfrak{p}_1}C(L_i,Y_i^{(k-1)+},Y_i^{(k-1)-}))\cap \bigcap\limits_{i\in\mathfrak{p}_2}C(L_i,P_i,Q_i)$$
 and
 $$\sum\limits_{i\in\mathfrak{p}_2}|Q_i|+\sum\limits_{i\in\mathfrak{p}_1}|Y_i^{(k-1)-}|=k.$$
 A contradiction because $f(Z^{(k-1)})<f(Z^{(k)})$.

{\it Case 2.2.2: $\sum\limits_{i\in\mathfrak{p}_1}|Y_i^{(k-1)-}|=0$.}

Then

 $\bullet$ $\sum\limits_{i\in\mathfrak{p}_2}|Y_i^{(k-1)-}|=k-1$;

 $\bullet$ $Y_i^{(k-1)-}=\varnothing$ for $i\in\mathfrak{p}_1$.

 Hence $C(L_{t_v},\{T^{(k-1)}\}\cup Y_{t_v}^{(k-1)+},Y_{t_v}^{(k-1)-})=\left<\{T^{(k-1)}\}\cup Y_{t_v}^{(k-1)+}\right>$ because $t_v\in\mathfrak{p}_1$.
 Therefore the intersection
 \begin{gather}\label{segment_from_Z_k-1_2}
  C(L_{t_v},\{T^{(k-1)}\}\cup Y_{t_v}^{(k-1)+},Y_{t_v}^{(k-1)-})\cap\bigcap\limits_{i\neq t_v}\aff(L_{i}\cup Y_{i}^{(k-1)+}\cup Y_{i}^{(k-1)-})
 \end{gather}
 is a segment.
 The  point of maximum of the functional $f$ on the segment~\eqref{segment_from_Z_k-1_2} is $Z^{(k-1)}$.

Denote by:

$\bullet$ $Z^{(k-1)''}$ the point of minimum of the functional $f$ on the segment~\eqref{segment_from_Z_k-1_2};

$\bullet$ $\Xi_{t_v}^{(k-1)+}$ a proper subset of $\{T^{(k-1)}\}\cup Y_{t_v}^{(k-1)+}$ such that  $Z^{(k-1)''}\in\left<\Xi_{t_v}^{(k-1)+}\right>$.

 Since
 $$(\R^d\times0)\cap C(L_{t_v},\{T^{(k-1)}\}\cup Y_{t_v}^{(k-1)+},Y_{t_v}^{(k-1)-})\cap\bigcap\limits_{i\neq t_v}\aff(L_{i}\cup Y_{i}^{(k-1)+}\cup Y_{i}^{(k-1)-})=\varnothing,$$
 the set $\Xi_{t_v}^{(k-1)+}$ is not subset of $\R^d\times0$.
 Then $(\{T^{(k-1)}\}\cup Y_{t_v}^{(k-1)+})\setminus\Xi_{t_v}^{(k-1)+}\subset\R^d\times0$.
 Take $\Xi_{t_v}^{(k-1)-}=(\{T^{(k-1)}\}\cup Y_{t_v}^{(k-1)+})\setminus\Xi_{t_v}^{(k-1)+}$.
 The intersection $\left<\Xi_{t_v}^{(k-1)+}\right>\cap\bigcap\limits_{i\neq t_v}\aff(L_{i}\cup Y_{i}^{(k-1)+}\cup Y_{i}^{(k-1)-})$ is unique point set.
 From sufficiently general position follows that $|\Xi_{t_v}^{(k-1)+}|=|Y_{t_v}^{(k-1)+}|$.
 Hence $|\Xi_{t_v}^{(k-1)-}|=1$.

Denote:

%$\bullet$ $\Xi_{i}^+=Y_i^{(k-1)+}$, $\Xi_{i}^-=Y_i^{(k-1)-}=\varnothing$ for $i\in\mathfrak{p}_1\setminus\{t_v\}$;

$\bullet$ $\Gamma_i^+=Y_i^{(k-1)+}$, $\Gamma_{i}^-=Y_i^{(k-1)-}$ for $i\neq t_v$;

$\bullet$ $\Gamma_{t_v}^+=Y_{t_v}^{(k-1)+}$, $\Gamma_{t_v}^-=\{T^{(k-1)}\}\cup Y_{t_v}^{(k-1)-}$.

For $i\in\mathfrak{p}_1\setminus\{t_v\}$ take a partition $\Xi_{i}^+$, $\Xi_{i}^-$ of $\Gamma_i^+\cup \Gamma_i^-$ such that   $Z^{(k-1)''}\in C(L_i,\Xi_{i}^+,\Xi_{i}^-)$.

Then :

$\bullet$ $\sum\limits_{i=1}^m|\Gamma_i^-|=k$;

$\bullet$ $Z^{(k-1)}$ is minimal point of $f$ on $\bigcap\limits_{i=1}^mC(L_i,\Gamma_i^+,\Gamma_i^-)$;

$\bullet$ $\Xi_i^+\cup \Xi_i^-=\Gamma_i^+\cup \Gamma_i^-$ for $i\in\mathfrak{p}_1$;

$\bullet$ $\max\limits_{\bigcap\limits_{i=1}^mC(L_i,\Xi_{i}^+,\Xi_{i}^-)}f=f(Z^{(k-1)''})<f(Z^{(k-1)})$;

$\bullet$ $\sum\limits_{i\in\mathfrak{p}_1}|\Xi_{i}^-|+\sum\limits_{i\in\mathfrak{p}_2}\lfloor\frac{d-\dim L_i}{2}\rfloor\geq 1+(k-1)=k$.

Apply Proposition~\ref{optimization_on_line} for these $\Gamma_1^+, ..., \Gamma_m^+, \Gamma_1^-, ..., \Gamma_m^-$ and $t=k-1$.
We get a pairwise disjoint subsets $$\Phi_1^+, ..., \Phi_m^+, \Phi_1^-, ..., \Phi_m^-\subset \bigcup\limits_{i=1}^m(\Gamma_i^+\cup\Gamma_i^-)$$ such that:

$\bullet$ $\Phi_1^-,...,\Phi_m^-\subset A$;

$\bullet$ $\max\limits_{\bigcap\limits_{i=1}^mC(L_i,\Phi_i^+,\Phi_i^-)}f >0$;

$\bullet$ $\min\limits_{\bigcap\limits_{i=1}^mC(L_i,\Phi_i^+,\Phi_i^-)}f < f(Z^{(k-1)})$;

$\bullet$  $\sum\limits_{i=1}^m|\Phi_i^-|\in\{k-1,k\}$;

$\bullet$ either   1) $\sum\limits_{i=1}^m|\Phi_i^-|=k$ or 2) $\dim\bigcap\limits_{i=1}^mC(L_i,\Phi_i^+,\Phi_i^-) =0$.

Since $f(Z^{(k-1)})\leq \min(f(Z^{(k-1)});f(Z^{(k)}))$, we can apply Proposition~\ref{optimization_on_line_2}.
Then $\sum\limits_{i=1}^m|\Phi_i^-|=k$ and $(\R^d\times0)\cap\bigcap\limits_{i=1}^mC(L_i,\Phi_i^+,\Phi_i^-)\neq\varnothing$.

Take:

$\bullet$ $A_i=\Phi_i^+\cap A$ for $i=1,...,m-1$;

$\bullet$ $B_i=\Phi_i^-$ for $i=1,...,m$;

$\bullet$ $A_m=A\setminus((\bigcup\limits_{i=1}^{m-1}A_i)\cup\bigcup\limits_{i=1}^mB_i)$.

Then

 $\bullet$ $\sum\limits_{i=1}^m|B_i|=k$;

  $\bullet$ $\bigcap\limits_{i=1}^mC(L_i, A_i, B_i)\neq\varnothing$.

  Theorem proved in this case.

  Therefore Theorem~\ref{main_result} proved in conditions 3.

\section{Moving the point}

In this section we proove some lemmas for using method of moving the points.
Let $f$ the last coordinate functional on $\R^{d+1}$.
We consider a collection $L_1, L_2, ... ,L_m$ of affine subspaces of $\R^d\times0\subset\R^{d+1}$.
Let $R$ be a finite subset of $\R^d\times[0;+\infty)$ such that the collection $L_1, L_2, ... ,L_m$ and set $R$ are  in strong sufficiently general position.

The main idea is similarly to [Tv66]: start with a family of subsets of $R$ which has the required conditions. Then, as moves one point continuously, we construct a family of sets with required conditions.

Take a hyperplane $H$ of $\R^d\times0$ which does not intersect the projection of $R$ to $\R^d\times0$.
Denote by $H^+,H^-$ a half-spaces into which $H$ divides the space $\R^d\times0$.

For a subset $Y\subset\R^{d+1}$ denote by $D^+(Y)$ and $D^-(Y)$ the polyhedra $$C(H\cap Y,(H^+\times\R)\cap Y,(H^-\times\R)\cap Y)$$
and
$$C(H\cap Y,(H^-\times\R)\cap Y,(H^+\times\R)\cap Y)$$
respectively.

For a subset $Y\subset\R^{d+1}$ and a number $i=1, ...,m$ denote $D_i^+(Y):=D^+(Y\cup L_i)$ and $D_i^-(Y):=D^-(Y\cup L_i)$.

\begin{definition}\label{sufficiently_divisible}
 Let $A$ is a subset of $\R^{d+1}$, $|A|=(m-1)(d+1)+1-\sum\limits_{i=1}^m\dim L_i$.
 The set $A$ is sufficiently divisible (with subspaces $L_1, ..., L_m$ and separating hyperplane $H$) if there is a partition $A_1, ..., A_m$ of $A$ such that

$\bullet$ $|A_i|\leq d+1-\dim L_i$ for $i=1,...,m$;

$\bullet$ either of the polyhedra $\bigcap\limits_{i=1}^m D_i^+(A_i)$ or $\bigcap\limits_{i=1}^m D_i^-(A_i)$
 does not contains in the hyperplane $\R^d\times0$ but intersect it.
\end{definition}

The two following lemmas describe the construction of family of sets with required conditions.

\begin{lemma}\label{moving_the_point}Let $X_0, X_1\in R$ be a points which does not separated by $H$.
Exactly one point of $X_0, X_1$ belong to $\R^d\times0$.
 Denote by $X_t$ a point $tX_1+(1-t)X_0$ for $t\in[0,1]$.
 Take a subset $Y\subset R\setminus\{X_0,X_1\}$ such that $|Y|=(m-1)(d+1)-\sum\limits_{i=1}^m\dim L_i$.
 Consider a set $T$ of numbers $t\in[0,1]$ such that set $ Y\cup\{X_t\}$ is sufficiently divisible.

Suppose that there is no sufficiently divisible subset $Y^*\subset R$ such that\newline $|Y^*\cap(\R^d\times(0;+\infty))|\leq |Y\cap(\R^d\times(0;+\infty))|$.

Then $T$ is open subset of $[0,1]$.

\end{lemma}

{\it Proof.} If $T=\varnothing$ then $T$ is open.

 Next, we suppose that $T\neq\varnothing$.
 For $t\in T$ take a partition $Y_1(t), ..., Y_m(t)$ of the set $ Y\cup\{X_t\}$ according to Definition~\ref{sufficiently_divisible}.
 Take $t_1\in T$.
 Without loss of generality the polyhedron $\bigcap\limits_{i=1}^m D_i^+(Y_i(t_1))$ does not contains in the hyperplane $\R^d\times0$ but intersect it.
Denote by $i_1$ a number such that $X_{t_1}\in Y_{i_1}(t_1)$ and:

$\bullet$ $\bar{Y}_{i_1}(t)=(Y_{i_1}(t_1)\setminus\{X_{t_1}\})\cup\{X_t\}$;

$\bullet$ $\bar{Y}_i(t)=Y_i(t_1)$ for $i\neq i_1$.

Let us prove some propositions.

\begin{proposition}\label{polyhedron_intersect_interior}
Suppose that for some partition $J_1,J_2$ of the set $\{1,...,m\}$ the polyhedron $\bigcap\limits_{i\in J_1} D_i^+(Y_i(t_1))\cap\bigcap\limits_{i\in J_2} \aff(Y_i(t_1))$:

$\bullet$ intersects a some open set $U$;

$\bullet$ does not contains in $\partial D_i^+(Y_i(t_1))$ for any $i\in J_1$.

Then there is an $\varepsilon>0$ such that the  polyhedron
$\bigcap\limits_{i\in J_1} D_i^+(\bar{Y}_i(t))\cap\bigcap\limits_{i\in J_2} \aff(\bar{Y}_i(t))$
intersects the set $U$ for any $t\in(t_1-\varepsilon;t_1+\varepsilon)$.
\end{proposition}

{\it Proof the Proposition.}

Take a set $R_1\subset U$ such that:

$\bullet$ the collection $L_1, L_2, ... ,L_m$ and set $R\cup R_1$ are  in sufficiently general position;

$\bullet$ the intersection of $\left< R_1\right>$ and the polyhedron $\bigcap\limits_{i\in J_1} D_i^+(Y_i(t_1))\cap\bigcap\limits_{i\in J_2} \aff(Y_i(t_1))$ is relative interior point of both of it. Denote this point by $\omega(t_1)$.

Take a subsets $$Y_1'(t_1)\subset Y_1(t_1)\cup L_1, ..., Y_m'(t_1)\subset Y_{m}(t_1)\cup L_m$$
such that

$\bullet$  $|Y_i'(t_1)|=1+\dim\aff(Y_i(t_1)\cup L_i)$;

$\bullet$ $\aff (Y_i'(t_1))=\aff(Y_i(t_1)\cup L_i)$.

Then  $\left<R_1\right>\cap\bigcap\limits_{i=1}^m\aff(Y_i'(t_1))=\{\omega(t_1)\}$.
For $x\in R_1\cup\bigcup\limits_{i=1}^mY_i'(t_1)$ there are a numbers $\alpha_x(t_1)$ such that:

$\bullet$ $\sum\limits_{x\in Y_i'(t_1)}\alpha_x(t_1)x=\omega(t_1)$ for $i=1,...,m$;

$\bullet$ $\sum\limits_{x\in Y_i'(t_1)}\alpha_x(t_1)=1$ for $i=1,...,m$;

$\bullet$ $\sum\limits_{x\in R_1}\alpha_x(t_1)x=\omega(t_1)$;

$\bullet$ $\sum\limits_{x\in R_1}\alpha_x(t_1)=1$.

The numbers $\{\alpha_x(t_1)|x\in R_1\cup\bigcup\limits_{i=1}^mY_i'(t_1)\}$ and coordinates of $\omega(t_1)$ is the unique solution of the following system of linear equations for $t=t_1$:

\begin{equation}\label{system_of_equations}
\begin{cases}
  \sum\limits_{x\in Y_i'(t)}\alpha_xx=\omega,       & \text{for } i=1,...,m \\
  \sum\limits_{x\in Y_i'(t)}\alpha_x=1,        & \text{for } i=1,...,m \\
  \sum\limits_{x\in R_1}\alpha_xx=\omega \\
  \sum\limits_{x\in R_1}\alpha_x=1
\end{cases}
\end{equation}
with variables $\alpha_x$ and coordinates of $\omega$.
This system has $(d+2)(m+1)$ variables and the same number of equations.
The matrix of the system~\eqref{system_of_equations} is non-degenerate for $t=t_1$.
The coefficients of this system are depend on $t$ continiously.
Hence the matrix of the system~\eqref{system_of_equations} is non-degenerate for $t$ in a some neighbourhood of $t_1$.
Therefore the solution of the system~\eqref{system_of_equations} depend on $t$ continiously.
Denote by $\omega(t)$ a point of the intersection $\left<R_1\right>\cap\bigcap\limits_{i=1}^m\aff(Y_i'(t))$.
Clearly, this notation has no conflict with notation of $\omega(t_1)$.
The point $\omega(t)$ depend on $t$ continiously because its coordinates is the solution of the system~\eqref{system_of_equations}.
Since $\omega(t_1)\in U$, there is an $\varepsilon>0$ such that $\omega(t)\in U$ for any $t\in (t_1-\varepsilon;t_1+\varepsilon)$.
Then the polyhedron $\bigcap\limits_{i\in J_1} D_i^+(\bar{Y}_i(t))\cap\bigcap\limits_{i\in J_2} \aff(\bar{Y}_i(t))$ intersect the set $U$ for any $t\in (t_1-\varepsilon;t_1+\varepsilon)$.

Proposition is proved.

\begin{proposition}\label{polyhedron_in_border}
 Suppose that the polyhedron $\bigcap\limits_{i=1}^mD_i^+(Y_i(t_1))$ contains in some proper face of the polyhedra $D_i^+(Y_{i}(t_1))$, $i=1,...,m$.

Take a minimal on inclusion subsets $Q_1\subset Y_1(t_1), ..., Q_m\subset Y_m(t_1)$ such that $\bigcap\limits_{i=1}^mD_i^+(Y_i(t_1))=\bigcap\limits_{i=1}^mD_i^+(Q_i)$.

Then the set $(\bigcup\limits_{i=1}^m Y_i(t_1))\setminus(\bigcup\limits_{i=1}^m Q_{i})$ is the unique point set.
\end{proposition}

{\it Proof the Proposition.}
Note that $$\sum\limits_{i=1}^m|Q_i|\leq-1+\sum\limits_{i=1}^m|Y_i(t_1)|=(m-1)(d+1)-\sum\limits_{i=1}^m\dim L_i.$$
The dimension of polyhedron $\bigcap\limits_{i=1}^mD_i^+(Q_i)$ at least 1 because this polyhedron does not contains in the hyperplane $\R^d\times0$ but intersect it.
Hence the collection of affine subspaces $L_1, ..., L_m$ and subset $Q_1\cup ...\cup Q_m$ does not in sufficiently general position.
Therefore $X_{t_1}\in Q_{i_1}$ and $X_{t_1}\notin\{X_0;X_1\}$.

If the polyhedron $D_i^+(Q_{i_1})$ is full-dimensional, the intersection $\bigcap\limits_{i\neq i_1}D_i^+(Q_i)$ is 1-dimensional.
According to sufficiently general position $$\sum\limits_{i\neq i_1}|Q_i|=1+(m-2)(d+1)-\sum\limits_{i\neq i_1}\dim L_i.$$
Then $\sum\limits_{i=1}^m|Q_i|=1+(m-1)(d+1)-\sum\limits_{i=1}^m\dim L_i$.
A contradiction.
Hence the polyhedron $D_i^+(Q_{i_1})$ is not full-dimensional.

Then
\begin{gather*}
 |(Q_{i_1}\setminus\{X_{t_1}\})\cup\{X_0,X_1\}|+\sum\limits_{i\neq i_1}|Q_i|=1+\sum\limits_{i=1}^m|Q_i|\leq \\ \leq\sum\limits_{i=1}^m|Y_i(t_1)|=(m-1)(d+1)+1-\sum\limits_{i=1}^m\dim L_i.
\end{gather*}
Since $D_{i_1}^+(Q_{i_1})\subset D_{i_1}^+((Q_{i_1}\setminus\{X_{t_1}\})\cup\{X_0,X_1\})$, the intersection $$D_{i_1}^+((Q_{i_1}\setminus\{X_{t_1}\})\cup\{X_0,X_1\})\cap\bigcap\limits_{i\neq i_1}D_i^+(Q_{i})$$
does not contains in the hyperplane $\R^d\times0$ but intersect it.
Therefore the dimension of this intersection at least 1.
According to sufficiently general position
\begin{gather*}
 |(Q_{i_1}\setminus\{X_{t_1}\})\cup\{X_0,X_1\}|+\sum\limits_{i\neq i_1}|Q_i|\geq(m-1)(d+1)+1-\sum\limits_{i=1}^m\dim L_i.
\end{gather*}
Then
\begin{gather*}
 |(Q_{i_1}\setminus\{X_{t_1}\})\cup\{X_0,X_1\}|+\sum\limits_{i\neq i_1}|Q_i|=(m-1)(d+1)+1-\sum\limits_{i=1}^m\dim L_i.
\end{gather*}

Hence the intersection $$D_{i_1}^+((Q_{i_1}\setminus\{X_{t_1}\})\cup\{X_0,X_1\})\cap\bigcap\limits_{i\neq i_1}D_i^+(Q_{i})$$ is 1-dimensional and does not contains in a proper face of $D_{i_1}^+((Q_{i_1}\setminus\{X_{t_1}\})\cup\{X_0,X_1\})$.
Since
\begin{gather*}
 \sum\limits_{i=1}^m|Q_i|=-1+|(Q_{i_1}\setminus\{X_{t_1}\})\cup\{X_0,X_1\}|+\sum\limits_{i\neq i_1}|Q_i|=\\=(m-1)(d+1)-\sum\limits_{i=1}^m\dim L_i=-1+\sum\limits_{i=1}^m|Y_i(t_1)|,
\end{gather*}
the set $(\bigcup\limits_{i=1}^m Y_i(t_1))\setminus(\bigcup\limits_{i=1}^m Q_{i})$ is the unique point set.

Proposition is proved.

\begin{proposition}\label{polyhedron_in_border_2}
In the conditions of Proposition~\ref{polyhedron_in_border} denote:

$\bullet$ by $T_1$ the unique point of the set $(\bigcup\limits_{i=1}^m Y_i(t_1))\setminus(\bigcup\limits_{i=1}^m Q_{i})$;

$\bullet$ $Q_{i_1}(t)=(Q_{i_1}\setminus\{X_{t_1}\})\cup\{X_{t}\}$ for $t\in [0;1]$;

$\bullet$ $Q_i(t)=Q_i$ for $i\neq i_1$ and $t\in [0;1]$.

Then $T_1\in \R^d\times0$ and for any open subset $U\subset\R^{d+1}$ there is a neighbourhood of $t_1$ such that for any $t$ from this neighbourhood there is a number $i_2$ such that the polyhedron:
$$D_{i_2}^+(Q_{i_2}(t)\cup\{T_1\})\cap\bigcap\limits_{i\neq i_2}D_i^+(Q_i(t))$$
intersect the set $U$.
\end{proposition}

{\it Proof the Proposition.}

If $T_1\in\R^d\times(0;+\infty)$ then
\begin{gather*}
 1+|Y\cap(\R^d\times(0;+\infty))|=|(\bigcup\limits_{i=1}^m Y_i(t_1))\cap(\R^d\times(0;+\infty))|>\\>|(\bigcup\limits_{i=1}^m Q_{i})\cap(\R^d\times(0;+\infty))|=|((Q_{i_1}\setminus\{X_{t_1}\})\cup\{X_0,X_1\}\cup\bigcup\limits_{i\neq i_1} Q_{i})\cap(\R^d\times(0;+\infty))|
\end{gather*}
because
$|((Q_{i_1}\setminus\{X_{t_1}\})\cup\{X_0,X_1\})\cap(\R^d\times(0;+\infty))|=|Q_{i_1}\cap(\R^d\times(0;+\infty))|$.
Hence
\begin{gather*}
 |Y\cap(\R^d\times(0;+\infty))|\geq|((Q_{i_1}\setminus\{X_{t_1}\})\cup\{X_0,X_1\}\cup\bigcup\limits_{i\neq i_1} Q_{i})\cap(\R^d\times(0;+\infty))|.
\end{gather*}
A contradiction with conditions of Lemma because the set $(Q_{i_1}\setminus\{X_{t_1}\})\cup\{X_0,X_1\}\cup\bigcup\limits_{i\neq i_1} Q_{i}$ is sufficiently divisible.
Therefore $T_1\notin\R^d\times(0;+\infty)$.
Then $T_1\in\R^d\times0$.

If the polyhedron $D_{i_1}^+((Q_{i_1}\setminus\{X_{t_1}\})\cup\{X_0,X_1\})$ is full-dimensional then the subspace $\aff D_{i_1}^+(Q_{i_1})$ is hyperplane which contains the polyhedron $\bigcup\limits_{i=1}^mD_i^+(Q_i)$.
Hence there are an intervals $I_1, I_2$ with $t_1$ as border point such that:

$\bullet$ $I_1\cap I_2=\varnothing$;

$\bullet$ the polyhedron $D_{i_1}^+(Q_{i_1}(t)\cup\{T_1\})\cap\bigcap\limits_{i\neq i_1} D_i^+(Q_i)$ intersect $U$ for any $t\in I_1$;

$\bullet$ the polyhedron $D_{i_1}^+(Q_{i_1}(t))\cap D^+(\{T_1\}\cup\bigcap\limits_{i\neq i_1} D_i^+(Q_i))$ intersect $U$ for any $t\in I_2$.

Then there is an interval $I_3\subset I_2$ with $t_1$ as border point such that for any $t\in I_3$ the set $D_{i_1}^+(Q_{i_1}(t))\cap U$ intersect a 2-dimensional face of the polyhedron $\bigcap\limits_{i\neq i_1} D_i^+(Q_i\cup\{T_1\})$ which contains the intersection $\bigcap\limits_{i\neq i_1} D_i^+(Q_i)$.
This face equals to $D_{i_2}^+(Q_{i_2}\cup\{T\})\cap\bigcap\limits_{i\neq i_1, i_2} D_i^+(Q_i)$ for some $i_2$.
Therefore the polyhedron
$$ D_{i_2}^+(Q_{i_2}(t)\cup\{T\})\cap\bigcap\limits_{i\neq  i_2} D_i^+(Q_i(t))$$ intersect $U$ for any $t\in I_3$.
The union $I_1\cup\{t_1\}\cup I_3$ is required neighbourhood of $t_1$.

Next, we consider remaining case when the polyhedron\newline $D_{i_1}^+((Q_{i_1}\setminus\{X_{t_1}\})\cup\{X_0,X_1\})$ is not full-dimensional.

There are an intervals $I_4, I_5$ with $t_1$ as border point such that:

$\bullet$ $I_4\cap I_5=\varnothing$;

$\bullet$ $\left<Q_{i_1}(t)\cup\{T_1\}\right>\cap\bigcap\limits_{i\neq i_1}D_{i}^+(Q_{i})$ intersect $U$ for any $t\in I_4$;

$\bullet$ $D_{i_1}^+(Q_{i_1}(t))\cap\left<(\bigcap\limits_{i\neq i_1}D_{i}^+(Q_{i}))\cup\{T_1\}\right>$ intersect $U$ for any $t\in I_5$.

Then for some permutation $\{I_4',I_5'\}$ of the set of intervals $\{I_4,I_5\}$:

$\bullet$ $D_{i_1}^+(Q_{i_1}(t)\cup\{T_1\})\cap\bigcap\limits_{i\neq i_1}D_{i}^+(Q_{i})$ intersect $U$ for any $t\in I_4'$;

$\bullet$ $D_{i_1}^+(Q_{i_1}(t))\cap D^+((\bigcap\limits_{i\neq i_1}D_{i}^+(Q_{i}))\cup\{T_1\})$ intersect $U$ for any $t\in I_5'$.

Hence there is an interval $I_6\subset I_5'$ with $t_1$ as border point such that the set $D_{i_1}^+(Q_{i_1}(t))\cap U$ intersect a some $1+\dim \bigcap\limits_{i\neq i_1} D_i^+(Q_i)$-dimensional face of the polyhedron $\bigcap\limits_{i\neq i_1} D_i^+(Q_i\cup\{T_1\})$ which contains the intersection $\bigcap\limits_{i\neq i_1} D_i^+(Q_i)$.
This face equals to $D_{i_2}^+(Q_{i_2}\cup\{T\})\cap\bigcap\limits_{i\neq i_1, i_2} D_i^+(Q_i)$ for some $i_2$.
Therefore the polyhedron
$$ D_{i_2}^+(Q_{i_2}(t)\cup\{T\})\cap\bigcap\limits_{i\neq  i_2} D_i^+(Q_i(t))$$ intersect $U$ for any $t\in I_6$.
The union $I_4'\cup\{t_1\}\cup I_6$ is required neighbourhood of $t_1$.

Proposition is proved.

For proving Lemma let us consider a subcases.

{\it Case 1: On the polyhedron $\bigcap\limits_{i=1}^m D_i^+(Y_i(t_1))$ both $min f<0$ and $max f>0$.}

Then this polyhedron intersect both open sets: $\R^d\times(-\infty;0)$ and $\R^d\times(0;+\infty)$.

Let us consider a subcases.

{\it Case 1.1: The polyhedron $\bigcap\limits_{i=1}^mD_i^+(Y_i(t_1))$ does not contains in any proper face of the polyhedra $D_i^+(Y_{i}(t_1))$, $i=1,...,m$.}

Hence there are numbers $\varepsilon_1, \varepsilon_2>0$ such that:

$\bullet$ the polyhedron $\bigcap\limits_{i=1}^m D_i^+(\bar{Y}_i(t))$ intersect the set $\R^d\times(-\infty;0)$ for any $t\in(t_1-\varepsilon_1;t_1+\varepsilon_1)$;

$\bullet$ the polyhedron $\bigcap\limits_{i=1}^m D_i^+(\bar{Y}_i(t))$ intersect the set $\R^d\times(0;+\infty)$ for any $t\in(t_1-\varepsilon_2;t_1+\varepsilon_2)$;

Hence for any $t$ such that $|t-t_1|<\min(\varepsilon_1;\varepsilon_2)$  the polyhedron $\bigcap\limits_{i=1}^m D_i^+(\bar{Y}_i(t))$ does not contains in  the set $\R^d\times(-\infty;0)$ but intersect it.
Therefore $t_1$ is internal point of $T$.

{\it Case 1.2: The polyhedron $\bigcap\limits_{i=1}^mD_i^+(Y_i(t_1))$ contains in some proper face of the polyhedra $D_i^+(Y_{i}(t_1))$, $i=1,...,m$.}

According to Proposition~\ref{polyhedron_in_border_2} there is a neighbourhood of $t_1$ such that for any $t$ from this neighbourhood there is a number $i_3$ such that the polyhedron:
$$D_{i_3}^+(Q_{i_3}(t)\cup\{T_1\})\cap\bigcap\limits_{i\neq i_3}D_i^+(Q_i(t))$$
intersect the open set $\R^d\times(0;+\infty)$.
If $t_1$ is not internal point of $T$ then there is an interval $I_7$ with $t_1$ as border point such that for any $t\in I_7$ the polyhedron
$$D_{i_3}^+(Q_{i_3}(t)\cup\{T_1\})\cap\bigcap\limits_{i\neq i_3}D_i^+(Q_i(t))$$
does not intersect the hyperplane $\R^d\times0$.
Take a point $t_2$ in the intersection of this neighbourhood and the interval $I_7$.
Consider a some segment which connect the polyhedron
$$D_{i_3}^+(Q_{i_3}(t_2)\cup\{T_1\})\cap\bigcap\limits_{i\neq i_3}D_i^+(Q_i(t_2))$$
with the set $(\R^d\times(-\infty;0))\cap\bigcap\limits_{i=1}^m D_i^+(Y_i(t_1))$.
This segment intersect the hyperplane $\R^d\times0$ in a some point $Z$.
For some point $t_3$ between the points $t_1$ and $t_2$ the point $Z$ belong to the polyhedron
$$D_{i_3}^+(Q_{i_3}(t_3)\cup\{T_1\})\cap\bigcap\limits_{i\neq i_3}D_i^+(Q_i(t_3)).$$

A contradiction with choice of the interval $I_7$.

{\it Case 2: On the polyhedron $\bigcap\limits_{i=1}^m D_i^+(Y_i(t_1))$ either $min f=0$ or $max f=0$.}

Take a point $Z_1$ in the relative interior of $(\R^d\times0)\cap \bigcap\limits_{i=1}^m D_i^+(Y_i(t_1))$.
If $Z_1\in D_{i_1}^+(Y_{i_1}(t_1)\setminus\{X_{t_1}\})$ then $Z_1\in \bigcap\limits_{i=1}^m D_i^+(\bar{Y}_i(t))$ for any $t\in [0;1]$.
Hence $T=[0;1]$.
Next, we consider the case where $Z_1\notin D_{i_1}^+(Y_{i_1}(t_1)\setminus\{X_{t_1}\})$.

Suppose that $t_1$ is not internal point of $T$.
There is an interval $I_8\subset[0,1]$ with $t_1$ as border point such that $I_8\cap T=\varnothing$.
Then $\bigcap\limits_{i=1}^m D_i^+(\bar{Y}_i(t))$ does not intersect the hyperplane $\R^d\times0$ for any $t\in I_8$.

 Let us call a family of pairwise disjoint subsets $P_1,..., P_m$ of $Y$ good if:

$\bullet$ $Z_1$ is the relative internal point of the polyhedron $D_{i_1}^+(P_{i_1}\cup X_{t_1})$;

$\bullet$ $\{Z_1\}= D_{i_1}^+(P_{i_1}\cup X_{t_1})\cap\bigcap\limits_{i\neq i_1}D_i^+(P_i)$;

$\bullet$ $\dim\aff (P_{i_1}\cup L_{i_1}\cup X_{t_1})+\sum\limits_{i\neq i_1}\dim\aff (P_i\cup L_i)=(m-1)(d+1)$;

$\bullet$ $\aff (P_{i_1}\cup L_{i_1}\cup X_{t})\cap\bigcap\limits_{i\neq i_1}D_i^+ (P_i)\neq \varnothing$ for some $t\in I_8$;

$\bullet$ $Y\setminus\{P_1\cup ...\cup P_m\}\subset\R^d\times0$.

The existence of good family will be proved later.
For a good family $P_1,..., P_m$ define an interval $I_{P_1,...,P_m}\subset I_8$ as follows:

The set $Y\setminus\{P_1\cup ...\cup P_m\}$ is unique point set because
\begin{gather*}
 \sum\limits_{i=1}^m|P_i|=\sum\limits_{i=1}^m\dim(P_i\cup L_i)-\sum\limits_{i=1}^m\dim L_i=\\
 =-1+\dim\aff (P_{i_1}\cup L_{i_1}\cup X_{t_1})+\sum\limits_{i\neq i_1}\dim(P_i\cup L_i)-\sum\limits_{i=1}^m\dim L_i=\\
 =-1+(m-1)(d+1)-\sum\limits_{i=1}^m\dim L_i=|Y|-1.
\end{gather*}
Denote:

$\bullet$ by $T_2$ the unique point of the set $Y\setminus\{P_1\cup ...\cup P_m\}$;

$\bullet$ $\dot{Y}_{i_1}(t)=P_{i_1}\cup \{X_t\}$ for any $t\in [0;1]$;

$\bullet$ $\dot{Y}_{i}(t)=P_i$ for $i\neq i_1$, $t\in[0;1]$.

Since $\bigcap\limits_{i=1}^mD_i^+(\dot{Y}_{i}(t))$ does not intersect the hyperplane $\R^d\times0$ for $t\in I_8$, the union $\bigcup\limits_{t\in I_8}(\bigcap\limits_{i=1}^mD_i^+(\dot{Y}_{i}(t)))$ contains in either half-spaces $\R^d\times(0;+\infty)$ or \newline$\R^d\times(-\infty;0)$.
Denote by $U_1$ one of open half-spaces $\R^d\times(-\infty;0)$, $\R^d\times(0;+\infty)$ which do not contains the union $\bigcup\limits_{t\in I_8}(\bigcap\limits_{i=1}^mD_i^+(\dot{Y}_{i}(t)))$.
There is a number $i_4$ such that the interestion $$U_1\cap D^+(\{T_2\}\cup \aff(\dot{Y}_{i_4}(t_1)\cup L_{i_4}))\cap\bigcap\limits_{i\neq i_4} \aff(\dot{Y}_{i}(t_1)\cup L_i)$$ is not empty because $T_2\in \R^d\times0$.
Denote:

$\bullet$ $\ddot{Y}_{i_4}(t)=\dot{Y}_{i_4}(t)\cup\{T_2\}$ for $t\in[0;1]$;

$\bullet$ $\ddot{Y}_{i}(t)=\dot{Y}_{i}(t)$ for $i\neq i_4$, $t\in[0;1]$.

Then the interestion $$U_1\cap D_{i_1}^+(\ddot{Y}_{i_1}(t_1))\cap\bigcap\limits_{i\neq i_1} \aff(\ddot{Y}_{i}(t_1)\cup L_i)$$
is not empty because the point $Z_1\in\bigcap\limits_{i=1}^m\aff(\ddot{Y}_i\cup L_i)$ is relative interior point of $D_{i_1}^+(\dot{Y}_{i_1}(t_1))$.
According to Proposition~\ref{polyhedron_intersect_interior} there is an $\varepsilon_5>0$ such that $$U_1\cap D_{i_1}^+(\ddot{Y}_{i_1}(t))\cap\bigcap\limits_{i\neq i_1} \aff(\ddot{Y}_{i}(t_1)\cup L_i)\neq\varnothing$$
for any $t\in(t_1-\varepsilon_5;t_1+\varepsilon_5)$.
Hence $$(\R^d\times0)\cap D_{i_1}^+(\ddot{Y}_{i_1}(t))\cap\bigcap\limits_{i\neq i_1} \aff(\ddot{Y}_{i}(t_1)\cup L_i)\neq\varnothing$$
for $t\in(t_1-\varepsilon_5;t_1+\varepsilon_5)$ because the sets $U_1$ and $\bigcup\limits_{t\in I_8}(\bigcap\limits_{i=1}^mD_i^+(\dot{Y}_{i}(t)))$ separated by the hyperplane $\R^d\times0$.

Since $Z_1$ is the relative internal point of the polyhedron $D_{i_1}^+(P_{i_1}\cup X_{t_1})$, there is an interval $I_{P_1,...,P_m}'\subset I_8$ with $t_1$ as border point such that $$D_{i_1}^+(P_{i_1}\cup X_{t})\cap\bigcap\limits_{i\neq i_1}\aff D_i^+(P_i)\neq\varnothing$$ for any $t\in I_{P_1,...,P_m}'$.
Then $\bigcap\limits_{i=1}^m\aff\dot{Y}_i(t)\neq\varnothing$ for any $t\in I_{P_1,...,P_m}'$.
%Hence $(\R^d\times0)\cap\bigcap\limits_{i=1}^m\aff\ddot{Y}_i(t)\neq\varnothing$ for any $t\in I_{P_1,...,P_m}'$.

There is an $\varepsilon_6>0$ such that $D_{i_1}^+(P_{i_1}\cup X_{t})\cap\bigcap\limits_{i\neq i_1}\aff(P_i\cup L_i)\neq\varnothing$ for any $t\in(t_1-\varepsilon_6;t_1+\varepsilon_6)$ because $ D_{i_1}^+(P_{i_1}\cup X_{t_1})\cap\bigcap\limits_{i\neq i_1}D_i^+(P_i)\neq\varnothing$ and $Z_1$ is relative interior point of $D_{i_1}^+(P_{i_1}\cup X_{t_1})$.
Then there is an interval $I_{P_1,...,P_m}''\subset I_8$ with $t_1$ as border point such that $\bigcap\limits_{i=1}^mD_i^+(\dot{Y}_i(t))= D_{i_1}^+(P_{i_1}\cup X_{t})\cap\bigcap\limits_{i\neq i_1}D_i^+(P_i)\neq\varnothing$ because there is a point $t\in I_8$ such that $\aff (P_{i_1}\cup L_{i_1}\cup X_{t})\cap\bigcap\limits_{i\neq i_1}D_i^+ (P_i)\neq \varnothing$.
Moreover, the intersection $\bigcap\limits_{i=1}^mD_i^+(\dot{Y}_i(t))$ is the relative interior point of any polyhedron $D_i^+(\dot{Y}_i(t))$ for $i=1,...,m$ and $t\in I_{P_1,...,P_m}''$.

Let $I_{P_1,...,P_m}=(t_1-\varepsilon_5;t_1+\varepsilon_5)\cap I_{P_1,...,P_m}'\cap I_{P_1,...,P_m}''$.
The interval $I_{P_1,...,P_m}$ is defined.

Denote by $I_9$ the intersection of all intervals $I_{P_1,...,P_m}$ for all good families $P_1,...,P_m$.
Clearly, $t_1$ is border point of the interval $I_9$.
Take a point $t_4\in I_9$. %which does not equal to any elements of $[0;1]$ given by Lemma~\ref{finite_number_of_points}.

Consider a good family $P_1,...,P_m$ such that the distance from the intersection $D_{i_1}^+(P_{i_1}\cup\{X_{t_4}\})\cap\bigcap\limits_{i\neq i_1}D_{i}^+(P_i)$ to the hyperplane $\R^d\times0$ is minimal.

Use notations $T_2, \dot{Y}_{i}(t), U_1, i_4, \ddot{Y}_{i}(t)$ from above:

$\bullet$ by $T_2$ the unique point of the set $Y\setminus\{P_1\cup ...\cup P_m\}$;

$\bullet$ $\dot{Y}_{i_1}(t)=P_{i_1}\cup \{X_t\}$ for any $t\in [0;1]$;

$\bullet$ $\dot{Y}_{i}(t)=P_i$ for $i\neq i_1$, $t\in[0;1]$;

$\bullet$ $U_1$ one of open half-spaces $\R^d\times(-\infty;0)$, $\R^d\times(0;+\infty)$ which do not contains the union $\bigcup\limits_{t\in I_8}(\bigcap\limits_{i=1}^mD_i^+(\dot{Y}_{i}(t)))$;

$\bullet$ $ D^+(\{T_2\}\cup \aff(\dot{Y}_{i_4}(t_1)\cup L_{i_4}))\cap\bigcap\limits_{i\neq i_4} \aff(\dot{Y}_{i}(t_1)\cup L_i)$ intersect the set $U_1$;

$\bullet$ $\ddot{Y}_{i_4}(t)=\dot{Y}_{i_4}(t)\cup\{T_2\}$ for $t\in[0;1]$;

$\bullet$ $\ddot{Y}_{i}(t)=\dot{Y}_{i}(t)$ for $i\neq i_4$, $t\in[0;1]$.

For $t\in I_9$ the intersection $(\R^d\times0)\cap D_{i_1}^+(\ddot{Y}_{i_1}(t))\cap\bigcap\limits_{i\neq i_1}\aff(\ddot{Y}_i(t)\cup L_i)$ is unique point.
Denote by $Z_1(t)$ this point.
The point $Z_1(t)$ continuously depend on $t$.

%\begin{gather*}
%\dim\aff(\ddot{Y}_{i_1}(t)\cup L_{i_1})=1+\dim\aff(\dot{Y}_{i_1}(t)\cup L_{i_1})=\\
%=\dim\aff(P_{i_1}\cup\{X_0;X_1\}\cup L_{i_1})=1+\dim(\R^d\times0)\cap\aff(P_{i_1}\cup\{X_0;X_1\}\cup L_{i_1})
%\end{gather*}

Since
$\dim\aff (P_{i_1}\cup L_{i_1}\cup X_{t_1})+\sum\limits_{i\neq i_1}\dim\aff (P_i\cup L_i)=(m-1)(d+1)$, we have
$\dim\aff (P_{i_1}\cup L_{i_1}\cup\{X_0;X_1\})+\sum\limits_{i\neq i_1}\dim\aff (P_i\cup L_i)=1+(m-1)(d+1)$.
Hence
\begin{gather*}
 \dim(\R^d\times0)\cap\aff (P_{i_1}\cup L_{i_1}\cup\{X_0;X_1\})+\sum\limits_{i\neq i_1}\dim(\R^d\times0)\cap\aff (P_i\cup L_i)=\\
 =(m-1)d.
\end{gather*}
According to strong sufficiently general position
\begin{gather*}
 \dim(\R^d\times0)\cap\aff (P_{i_1}\cup L_{i_1}\cup\{X_0;X_1\})\cap\bigcap\limits_{i\neq i_1}\aff (P_i\cup L_i)=0.
\end{gather*}
Then the interestion $$(\R^d\times0)\cap\aff (P_{i_1}\cup L_{i_1}\cup\{X_0;X_1\})\cap\bigcap\limits_{i\neq i_1}\aff (P_i\cup L_i)$$ is unique point set which does not belong to any subspace $(\R^d\times0)\cap\aff (P_i'\cup L_i)$ for $i\neq i_1$ and proper subset $P_i'\subset P_i$.
The interestion $$(\R^d\times0)\cap\aff (P_{i_1}\cup L_{i_1}\cup\{X_0;X_1\})\cap\bigcap\limits_{i\neq i_1}\aff (P_i\cup L_i)$$
equal to $\{Z_1\}$ because
\begin{gather*}
 \{Z_1\}=D_{i_1}^+(P_{i_1}\cup X_{t_1})\cap\bigcap\limits_{i\neq i_1}D_i^+(P_i)\subset
 \\\subset(\R^d\times0)\cap\aff (P_{i_1}\cup L_{i_1}\cup\{X_0;X_1\})\cap\bigcap\limits_{i\neq i_1}\aff (P_i\cup L_i).
\end{gather*}
Therefore $Z_1$ is the relative interior point of the polyhedron $(\R^d\times0)\cap D_i^+(P_i)$ for $i\neq i_1$.

If $i_4=i_1$ then $P_i=\ddot{Y}_i(t_1)$ for any $i\neq i_1$ and $t\in I_9$.
Hence $Z_1(t_1)=Z_1$ is relative interior point of $(\R^d\times0)\cap D_i^+(\ddot{Y}_i(t_1))$ for $i\neq i_1$.
%Therefore the point $Z_1$ is relative interior point of the polyhedron $\bigcap\limits_{i\neq i_1}D_i^+(\ddot{Y}_i(t_1))$.
Since $$\aff D_i^+(\ddot{Y}_i(t_1))=\aff(\ddot{Y}_i(t_1)\cup L_i)$$ for $i\neq i_1$,
there is a neighbourhood of the point $t_1$ such that $Z_1(t)\in D_i^+(\ddot{Y}_i(t_1))$ for any $t$ from this neighbourhood and $i\neq i_1$.
Take some point $t_5$ from the interestion of this neighbourhood and $I_8$.
Then
\begin{gather*}
 Z_1(t_5)\in D_{i_1}^+(P_{i_1}\cup X_{t_5})\cap\bigcap\limits_{i\neq i_1}D_i^+(\ddot{Y}_i(t_1))=\\
 =D_{i_1}^+(P_{i_1}\cup X_{t_5})\cap\bigcap\limits_{i\neq i_1}D_i^+(P_i)=D_{i_1}^+(P_{i_1}\cup X_{t_5})\cap\bigcap\limits_{i\neq i_1}D_i^+(\ddot{Y}_i(t_5))\subset\bigcap\limits_{i=1}^m D_i^+(\ddot{Y}_i(t_5)).
\end{gather*}
Then the polyhedron $\bigcap\limits_{i=1}^m D_i^+(\ddot{Y}_i(t_5))$ intersect the hyperplane $\R^d\times0$.
A contradiction with the choice of $I_8$ because the polyhedron $\bigcap\limits_{i=1}^m D_i^+(\ddot{Y}_i(t_5))$ does not contains in the hyperplane $\R^d\times0$.
Hence $i_4\neq i_1$.
Therefore $\ddot{Y}_{i_1}(t)=\dot{Y}_{i_1}(t)=P_{i_1}\cup\{X_{t}\}$ for any $t\in [0;1]$.

According to the choice of $I_8$ the polyhedron $$\bigcap\limits_{i=1}^mD_i^+(\ddot{Y}_i(t))=D_{i_1}^+(P_{i_1}\cup\{X_{t}\})\cap\bigcap\limits_{i\neq i_1}D_i^+(\ddot{Y}_i(t))$$ does not intersect the hyperplane $\R^d\times0$ for any $t\in I_8$.
Then the polyhedron $D_{i_1}^+(P_{i_1}\cup\{X_{t_1},X_{t_4}\})\cap\bigcap\limits_{i\neq i_1}D_i^+(\ddot{Y}_i(t))$ does not intersect the half-space $U_1$.
The point $Z_1$ is unique point of the intersection of hyperplane $\R^d\times0$ and polyhedron $D_{i_1}^+(P_{i_1}\cup\{X_{t_1},X_{t_4}\})\cap\bigcap\limits_{i\neq i_1}D_i^+(\ddot{Y}_i(t))$.

According to the choice of the point $t_4\in I_9\subset I_{P_1,...,P_m}$ the interestion $\bigcap\limits_{i=1}^mD_i^+(\dot{Y}_i(t_4))$ is relative interior point of any polyhedron $D_i^+(\dot{Y}_i(t_4))$ for $i=1,...,m$.
Then this point is endpoint of both 1-dimensional polyhedra: $$\bigcap\limits_{i=1}^m D_i^+(\ddot{Y}_i(t_4))$$
and
$$D_{i_1}^+(\ddot{Y}_{i_1}(t_4))\cap\bigcap\limits_{i\neq i_1} \aff(\ddot{Y}_i(t_4)\cup L_{i}).$$
The point $\bigcap\limits_{i=1}^m D_i^+(\dot{Y}_i(t_4))$ is the most far point of the polyhedron $\bigcap\limits_{i=1}^m D_i^+(\ddot{Y}_i(t_4))$ from the hyperplane $\R^d\times0$ because $$(\R^d\times0)\cap D^+(\{T_2\}\cup D_{i_4}^+(\dot{Y}_{i_4}(t_4)))\cap\bigcap\limits_{i\neq i_4} \aff(\ddot{Y}_i(t_4)\cup L_{i})\neq\varnothing.$$
According to the choice of interval $I_8$ the polyhedron $\bigcap\limits_{i=1}^m D_i^+(\ddot{Y}_i(t_4))$ does not intersect the hyperplane $\R^d\times0$.
Then the polyhedron $\bigcap\limits_{i=1}^m D_i^+(\ddot{Y}_i(t_4))$ is a segment.
Denote by $Z_2$ the endpoint of the segment $\bigcap\limits_{i=1}^m D_i^+(\ddot{Y}_i(t_4))$ which does not equal to $\bigcap\limits_{i=1}^m D_i^+(\dot{Y}_i(t_4))$.
Since $(\R^d\times0)\cap\bigcap\limits_{i=1}^m D_i^+(\ddot{Y}_i(t_4))=\varnothing$, the points $Z_2$ and $\bigcap\limits_{i=1}^m D_i^+(\dot{Y}_i(t_4))$ does not separated by the hyperplane $\R^d\times0$.

The polyhedron $ D^+(\{T_2\}\cup \aff(\dot{Y}_{i_4}(t_1)\cup L_{i_4}))\cap\bigcap\limits_{i\neq i_4} \aff(\dot{Y}_{i}(t_1)\cup L_i)$:

$\bullet$ is 1-dimensional;

$\bullet$ contains the point $\{Z_1\}=\bigcap\limits_{i=1}^m\aff(\dot{Y}_{i}(t_1)\cup L_i)\subset\R^d\times0$ as border point;

$\bullet$ intersect the half-space $U_1$.

Then this polyhedron:

$\bullet$ contains in the closed half-space $U_1\cup(\R^d\times0)$;

$\bullet$ intersect the hyperplane $\R^d\times0$.

Hence the  polyhedron $$\bigcap\limits_{i=1}^m D_i^+(\ddot{Y}_i(t_1))\subset D^+(\{T_2\}\cup \aff(\dot{Y}_{i_4}(t_1)\cup L_{i_4}))\cap\bigcap\limits_{i\neq i_4} \aff(\dot{Y}_{i}(t_1)\cup L_i):$$

$\bullet$ intersect the hyperplane $\R^d\times0$ by the unique point $Z_1$

$\bullet$ contains in the closed half-space $U_1\cup(\R^d\times0)$.

If the polyhedron $\bigcap\limits_{i=1}^m D_i^+(\ddot{Y}_i(t_1))$ is 1-dimensional then it intersect the set $U_1$.
Hence the polyhedron $$D_{i_1}^+(P_{i_1}\cup\{X_{t_1},X_{t_4}\})\cap\bigcap\limits_{i\neq i_1} D_i^+(\ddot{Y}_i(t_1))$$
both intersect the set $U_1$ and contains the point $Z_2\notin U_1$.
Therefore there is a point $t_6$ between the points $t_1$ and $t_4$ such that the polyhedron $$D_{i_1}^+(P_{i_1}\cup\{X_{t_6}\})\cap\bigcap\limits_{i\neq i_1} D_i^+(\ddot{Y}_i(t_1))$$
 intersects both sets $U_1$ and $\R^d\times0$.
 A contradiction because $t_6\in I_8$.
Then the polyhedron $\bigcap\limits_{i=1}^m D_i^+(\ddot{Y}_i(t_1))$ is 0-dimensional.
Hence $\{Z_1\}=\bigcap\limits_{i=1}^m D_i^+(\ddot{Y}_i(t_1))$.
Therefore the point $Z_1$ is a vertex of the polyhedron $$D_{i_1}^+(P_{i_1}\cup\{X_0,X_1\})\cap\bigcap\limits_{i\neq i_1}D_i^+(\ddot{Y}_i(t_1)).$$
Since $\{Z_1\}=\bigcap\limits_{i=1}^m D_i^+(\ddot{Y}_i(t_1))$, the edges of the polyhedron $$D_{i_1}^+(P_{i_1}\cup\{X_0,X_1\})\cap\bigcap\limits_{i\neq i_1}D_i^+(\ddot{Y}_i(t_1))$$ which contains the vertex $Z_1$ are not separated by $D_{i_1}^+(P_{i_1}\cup\{X_{t_1}\})=D_{i_1}^+(\ddot{Y}_{i_1}).$
Since the polyhedron $$D_{i_1}^+(P_{i_1}\cup\{X_0,X_1\})\cap\bigcap\limits_{i\neq i_1}D_i^+(\ddot{Y}_i(t_1))$$ is 2-dimensional, there is two edges of this polyhedron which  contains the vertex $Z_1$.
One of these edges equals to $$D_{i_1}^+(P_{i_1}\cup\{X_0,X_1\})\cap\bigcap\limits_{i\neq i_1}D_i^+(\dot{Y}_i(t_1)).$$
Other edge equal to $$D_{i_1}^+(P_{i_1}\cup\{X_0,X_1\})\cap\bigcap\limits_{i\neq i_1}D_i^+(P_i'')$$
for some subsets $P_{i}''\subset \ddot{Y}_i(t_1)$ for $i\neq i_1$.
Then $$\{Z_1\}= D_{i_1}^+(P_{i_1}\cup X_{t_1})\cap\bigcap\limits_{i\neq i_1}D_i^+(P_i'').$$
According to sufficiently general position $$\dim\aff (P_{i_1}\cup L_{i_1}\cup \{X_{0},X_1\})+\sum\limits_{i\neq i_1}\dim\aff (P_i''\cup L_i)=1+(m-1)(d+1)$$ because the polyhedron $$D_{i_1}^+(P_{i_1}\cup\{X_0,X_1\})\cap\bigcap\limits_{i\neq i_1}D_i^+(P_i'')$$ is 1-dimensional.
Hence $$\dim\aff (P_{i_1}\cup L_{i_1}\cup X_{t_1})+\sum\limits_{i\neq i_1}\dim\aff (P_i''\cup L_i)=(m-1)(d+1).$$
Since the edges of the polyhedron $$D_{i_1}^+(P_{i_1}\cup\{X_0,X_1\})\cap\bigcap\limits_{i\neq i_1}D_i^+(\ddot{Y}_i(t_1))$$ which contains the point $Z_1$ does not separated by $\aff (P_{i_1}\cup L_{i_1}\cup X_{t_1})$ we have
 $\aff (P_{i_1}\cup L_{i_1}\cup X_{t})\cap\bigcap\limits_{i\neq i_1}D_i^+ (P_i'')\neq \varnothing$ for some $t\in I_9\subset I_8$.
Denote $P_{i_1}''=P_{i_1}$.
 The set $\{X_0,X_1\}\cup\bigcup\limits_{i=1}^mP_i''$ is sufficiently divisible because the polyhedron $$D_{i_1}^+(P_{i_1}\cup\{X_0,X_1\})\cap\bigcap\limits_{i\neq i_1}D_i^+(P_i'')$$
 does not contains in the hyperplane $\R^d\times0$ but intersect it.
 Then $$1+\sum\limits_{i=1}^m|(\R^d\times(0;+\infty))\cap P_i''|>|(\R^d\times(0;+\infty))\cap Y|.$$
 Hence $$\sum\limits_{i=1}^m|(\R^d\times(0;+\infty))\cap P_i''|=|(\R^d\times(0;+\infty))\cap Y|.$$
 Therefore $Y\setminus\{\bigcup\limits_{i=1}^mP_i''\}\subset \R^d\times0$.
 Since $Z_1$ is the relative internal point of the polyhedron $D_{i_1}^+(P_{i_1}\cup \{X_{t_1}\})$, the family $P_1'',...,P_m''$ is good.
 According to the choice of point $t_4\in I_9$ the interestion $D_{i_1}^+(P_{i_1}\cup\{X_{t_4}\})\cap\bigcap\limits_{i\neq i_1}D_i^+(P_i'')$ is not empty. %relative interior point of
 The polyhedron $D_{i_1}^+(P_{i_1}\cup\{X_0, X_1\})\cap\bigcap\limits_{i\neq i_1}D_i^+(P_i'')$ is proper face of the polyhedron $D_{i_1}^+(P_{i_1}\cup\{X_0, X_1\})\cap\bigcap\limits_{i\neq i_1}D_i^+(\ddot{Y}_i(t_1))$.
 Then the interestion $D_{i_1}^+(P_{i_1}\cup\{X_{t_4}\})\cap\bigcap\limits_{i\neq i_1}D_i^+(P_i'')$ is endpoint of the segment $$D_{i_1}^+(P_{i_1}\cup\{X_{t_4}\})\cap\bigcap\limits_{i\neq i_1}D_i^+(\ddot{Y}_i(t_1))=\bigcap\limits_{i=1}^mD_i^+(\ddot{Y}_i(t_4)).$$
Hence $\{Z_2\}=D_{i_1}^+(P_{i_1}\cup\{X_{t_4}\})\cap\bigcap\limits_{i\neq i_1}D_i^+(P_i'')$.
A contradiction with the choice of the family $P_1, ..., P_m$ because the point $Z_2$ closer to the hyperplane $\R^d\times0$ than the point $D_{i_1}^+(P_{i_1}\cup X_{t_4})\cap\bigcap\limits_{i\neq i_1}D_i^+(P_i)$

The proof of the existence of a family of good sets is required.
Take a partition $\tilde{Y}_1,...,\tilde{Y}_m$ of $Y$ such that:

$\bullet$ $Z_1\in \partial D_{i_1}^+(\tilde{Y}_{i_1}\cup\{X_0,X_1\})\cap\bigcap\limits_{i\neq i_1}D_i^+(\tilde{Y}_i)$;

$\bullet$ the interestion $D_{i_1}^+(\tilde{Y}_{i_1}\cup\{X_{t}\})\cap\bigcap\limits_{i\neq i_1}D_i^+(\tilde{Y}_i)$ is not empty for some $t''\in I_8$.

The existence of a partition $\tilde{Y}_1,...,\tilde{Y}_m$  will be proved later.
For $t$ between $t_1$ and $t''$ denote by $Z_3(t)$ most close point of the polyhedron $D_{i_1}^+(\tilde{Y}_{i_1}\cup\{X_t\})\cap\bigcap\limits_{i\neq i_1}D_i^+(\tilde{Y}_{i})$ to the hyperplane $\R^d\times0$.
Clearly:

$\bullet$ $\lim\limits_{t\rightarrow t_1} Z_3(t)=Z_1$;

$\bullet$ for any $t\in [0;1]$ the point $Z_3(t)$ belong to $\partial D_{i_1}^+(\tilde{Y}_{i_1}\cup\{X_0,X_1\})\cap\bigcap\limits_{i\neq i_1}D_i^+(\tilde{Y}_{i})$.

Then there is an egde $\mathfrak{E}_1$ of $ D_{i_1}^+(\tilde{Y}_{i_1}\cup\{X_0,X_1\})\cap\bigcap\limits_{i\neq i_1}D_i^+(\tilde{Y}_{i})$ which contains the points $Z_1$ and $Z_3(t)$ for some $t\in I_8$.
Then $$\mathfrak{E}_1\subset \partial D_{i_1}^+(\tilde{Y}_{i_1}\cup\{X_0,X_1\})\cap\bigcap\limits_{i\neq i_1}D_i^+(\tilde{Y}_{i})$$
From Lemma~\ref{polyhedron_border} follows that
$$\mathfrak{E}_1= D_{i_1}^+(P_{i_1}\cup\{X_0,X_1\})\cap\bigcap\limits_{i\neq i_1}D_i^+(P_{i}(t_1))$$
for some $P_{i_1}\subset \tilde{Y}_{i_1}$, $P_i\subset \tilde{Y}_{i}$ for $i\neq i_1$.
According to sufficiently general position $$\dim\aff (P_{i_1}\cup L_{i_1}\cup \{X_{0},X_1\})+\sum\limits_{i\neq i_1}\dim\aff (P_i\cup L_i)=1+(m-1)(d+1)$$ because $\mathfrak{E}_1\not\subset\R^d\times0$.
Hence $$\dim\aff (P_{i_1}\cup L_{i_1}\cup X_{t_1})+\sum\limits_{i\neq i_1}\dim\aff (P_i\cup L_i)=(m-1)(d+1).$$
Since \begin{gather*}
\{Z_1\}=(\R^d\times0)\cap\mathfrak{E}_1=\\
=(\R^d\times0)\cap D_{i_1}^+(P_{i_1}\cup\{X_0,X_1\})\cap\bigcap\limits_{i\neq i_1}D_i^+(P_{i}(t_1)),
\end{gather*}
from strong sufficiently general position follows that the point $Z_1$ belong to the relative interior of the polyhedron $(\R^d\times0)\cap D_{i_1}^+(P_{i_1}\cup\{X_0,X_1\})$.
Therefore the point $Z_1$ belong to the relative interior of the polyhedron $ D_{i_1}^+(P_{i_1}\cup\{X_0,X_1\})$ because this polyhedron separated by the hyperplane $\R^d\times0$.
Then the point $Z_1$ is the relative internal point of the polyhedron $D_{i_1}^+(P_{i_1}\cup X_{t_1})$.
Clearly $\{Z_1\}= D_{i_1}^+(P_{i_1}\cup X_{t_1})\cap\bigcap\limits_{i\neq i_1}D_i^+(P_i)$;
For $t\in I_8$ such that $Z_3(t)\in \mathfrak{E}_1$ we have $\aff (P_{i_1}\cup L_{i_1}\cup X_{t})\cap\bigcap\limits_{i\neq i_1}D_i^+ (P_i)\neq \varnothing$.
Since $\mathfrak{E}_1$ does not belong to hyperplane $\R^d\times0$ but intersect it, the set $\{X_0,X_1\}\cup P_1\cup...\cup P_m$ is sufficiently divisible.
Then $1+\sum\limits_{i=1}^m|(\R^d\times(0;+\infty))\cap P_i|>|(\R^d\times(0;+\infty))\cap Y|$.
Hence $$\sum\limits_{i=1}^m|(\R^d\times(0;+\infty))\cap P_i|=|(\R^d\times(0;+\infty))\cap Y|.$$
Therefore $Y\setminus\{P_1\cup ...\cup P_m\}\subset\R^d\times0$.

The family $P_1, ..., P_m$ is good.
The existence of a good family is proved.

The proof of the existence of a partition $\tilde{Y}_1,...,\tilde{Y}_m$ of $Y$ is required.
Let us consider a subcases:

{\it Case 2.1: The polyhedron $\bigcap\limits_{i=1}^mD_i^+(Y_i(t_1))$ does not contains in any proper face of the polyhedra $D_i^+(Y_{i}(t_1))$, $i=1,...,m$.}

Since $Z_1\in\partial \bigcap\limits_{i=1}^mD_i^+(Y_i(t_1))$, we have $$Z_1\in\partial D_{i_1}^+((Y_{i_1}(t_1)\setminus\{X_{t_1}\})\cup\{X_0,X_1\})\cap\bigcap\limits_{i\neq i_1}D_i^+(Y_i(t_1)).$$
According to Proposition~\ref{polyhedron_intersect_interior} there is a neighbourhood of the point $t_1$ such that for any $t$ from this neighbourhood the interestion $(\R^d\times(\R\setminus\{0\}))\cap\bigcap\limits_{i=1}^m D_i^+(\bar{Y}_i(t))$ is not empty.
Take a point $t''$ from intersection of this neighbourhood and the interval $I_8$.
Take required family as:
 $\tilde{Y}_{i}=Y_{i}(t_1)\cap Y$ for $i=1,...,m$.

{\it Case 2.2: The polyhedron $\bigcap\limits_{i=1}^mD_i^+(Y_i(t_1))$ contains in some proper face of the polyhedra $D_i^+(Y_{i}(t_1))$, $i=1,...,m$.}

Use notations from Propositions~\ref{polyhedron_in_border} and~\ref{polyhedron_in_border_2}.
According to Proposition~\ref{polyhedron_in_border_2} there is a neighbourhood of $t_1$ such that for any $t$ from this neighbourhood there is a number $i_2$ such that the polyhedron:
$$D_{i_2}^+(Q_{i_2}(t)\cup\{T_1\})\cap\bigcap\limits_{i\neq i_2}D_i^+(Q_i(t))$$
intersect the set $\R^d\times(\R\setminus\{0\})$.
Take a point $t''$ from intersection of this neighbourhood and the interval $I_8$.
Take required family as:
 $\tilde{Y}_{i}=Y_{i}(t_1)\cap Y$ for $i=1,...,m$.

 The partition $\tilde{Y}_1,...,\tilde{Y}_m$ of $Y$ is constructed in both cases 2.1 and 2.2.
 Lemma is proved.

\begin{lemma}\label{moving_the_point_2}
 In the conditions of Lemma~\ref{moving_the_point}
 for any border point of the set $T$ there is a punctured neighbourhood of this point which contains in the set $T$.
 Neither border point of $T$ is endpoint of the segment $[0;1]$.
\end{lemma}

{\it Proof.}
Denote by $t_0$ some border point of $T$.
There is a sequence of points $t_1, t_2, ...$ of the set $T$ such that $\lim\limits_{i\rightarrow\infty}t_i=t_0$.
For $t\in T$ take a partition $Y_1(t), ..., Y_m(t)$ of the set $ Y\cup\{X_t\}$ according to Definition~\ref{sufficiently_divisible}.
There is a subsequence $t_1', t_2', ...$ of the sequence $t_1, t_2, ...$ such that:

$\bullet$ the interestion $Y_i(t_j')\cap Y$ does not depend on $j$ for any $i=1,...,m$;

$\bullet$ there is a number $i_1$ such that $X_{t_j'}\in Y_{i_1}(t_j')$ for any $j$.

%Without loss of generality $(\R^d\times0)\cap\bigcap\limits_{i=1}^mD_i^+(Y_{i}(t_j))\neq\varnothing$ for any $j$.
%These interestions are unique point sets.
Denote by $Z_j$ unique point of interestion $(\R^d\times0)\cap\bigcap\limits_{i=1}^m\aff(Y_{i}(t_j)\cup L_i)$ for any $j$.
Since $t_j\in I$ for $j=1, 2, ...$, the point $Z_j$ belong to either
$$(\R^d\times0)\cap\bigcap\limits_{i=1}^mD_i^+(Y_{i}(t_j))$$
or
$$(\R^d\times0)\cap\bigcap\limits_{i=1}^mD_i^-(Y_{i}(t_j)).$$

Suppose that the sequence $Z_1, Z_2, ...$ has a limit.
Denote it by $Z_0$.
Hence $Z_0$ belong to either
$$(\R^d\times0)\cap D_i^+(Y_{i_1}\cup\{X_{t_0}\})\cap\bigcap\limits_{i\neq i_1}D_i^+(Y_{i})$$
or
$$(\R^d\times0)\cap D_i^-(Y_{i_1}\cup\{X_{t_0}\})\cap\bigcap\limits_{i\neq i_1}D_i^-(Y_{i}).$$
These polyhedra does not contains in the hyperplane.
A contradiction because $t_0\notin T$.
Therefore the sequence $Z_1, Z_2, ...$ has no limit.

Then the line
$$\aff(Y_{i_1}\cup L_{i_1}\cup\{X_{t_0}\})\cap\bigcap\limits_{i\neq i_1}\aff(Y_i\cup L_i)$$
does not intersect the hyperplane $\R^d\times0$.

The set $Y\cup\{X_{t_0}\}$ does not in strong general position with the subspaces with affine subspaces $L_1, ..., L_m$ because
\begin{multline*}
 \dim(\R^d\times0)\cap \aff(Y_{i_1}\cup L_{i_1}\cup\{X_{t_0}\})+\sum\limits_{i\neq i_1}\dim(\R^d\times0)\cap \aff(Y_i\cup L_i)\neq\\
 \neq-1=\dim\varnothing=\dim(\R^d\times0)\cap\aff(Y_{i_1}\cup L_{i_1}\cup\{X_{t_0}\})\cap\bigcap\limits_{i\neq i_1}\aff(Y_i\cup L_i).
\end{multline*}
Hence $X_{t_0}\notin\{X_0,X_1\}$.
Therefore $t_0\notin\{0,1\}$.

The line
$$\aff(Y_{i_1}\cup L_{i_1}\cup\{X_{t_0}\})\cap\bigcap\limits_{i\neq i_1}\aff(Y_i\cup L_i)$$
is not parallel to the line
$$\aff(Y_{i_1}\cup L_{i_1}\cup\{X_{t_1}\})\cap\bigcap\limits_{i\neq i_1}\aff(Y_i\cup L_i)$$
because it intersect the hyperplane $\R^d\times0$.
The subspace
$$\aff(Y_{i_1}\cup L_{i_1}\cup\{X_{0},X_1\})\cap\bigcap\limits_{i\neq i_1}\aff(Y_i\cup L_i)$$
is 2 -dimensional and contains these lines.
Then the intersection of the lines
$$\aff(Y_{i_1}\cup L_{i_1}\cup\{X_{t_0}\})\cap\bigcap\limits_{i\neq i_1}\aff(Y_i\cup L_i)$$
and
$$\aff(Y_{i_1}\cup L_{i_1}\cup\{X_{t_1}\})\cap\bigcap\limits_{i\neq i_1}\aff(Y_i\cup L_i)$$
is not empty.
Hence the intersection
$\bigcap\limits_{i=1}^m\aff(Y_i\cup L_i)$
is unique point.
For any $t\in[0;1]\setminus\{t_0\}$ the affine hull of the lines
$$\aff(Y_{i_1}\cup L_{i_1}\cup\{X_{t_0}\})\cap\bigcap\limits_{i\neq i_1}\aff(Y_i\cup L_i)$$
and
$$\aff(Y_{i_1}\cup L_{i_1}\cup\{X_{t}\})\cap\bigcap\limits_{i\neq i_1}\aff(Y_i\cup L_i)$$
equal to the subspace
$$\aff(Y_{i_1}\cup L_{i_1}\cup\{X_{0},X_1\})\cap\bigcap\limits_{i\neq i_1}\aff(Y_i\cup L_i)$$ which is not parallel to the hyperplane $\R^d\times0$.
Therefore the line
$$\aff(Y_{i_1}\cup L_{i_1}\cup\{X_{t}\})\cap\bigcap\limits_{i\neq i_1}\aff(Y_i\cup L_i)$$ is not parallel to to the hyperplane $\R^d\times0$.

Then there is a punctured neighbourhood of $t_0$ such that the union of the rays:
$$ D_i^+(Y_{i_1}\cup\{X_{t}\})\cap\bigcap\limits_{i\neq i_1}D_i^+(Y_{i})$$
or
$$ D_i^-(Y_{i_1}\cup\{X_{t}\})\cap\bigcap\limits_{i\neq i_1}D_i^-(Y_{i})$$
interest the hyperplane $\R^d\times0$.
Therefore this punctured neighbourhood contains in the set $T$.
Lemma is proved.

\section{Proof by moving the point}

In [BS16] Theorem~\ref{th:pm} proved in stronger form: Theorem~\ref{th:tverbM}. Similarly, we prove  Theorem~\ref{main_result_2} in a stronger form: Theorem~\ref{partition_of_separated_set}. %Theorem~\ref{th:tverbM} follows from
 % The case $L_1=...=L_m=\varnothing$ in Theorem~\ref{partition_of_separated_set} implies Theorem~\ref{th:tverbM} similarly to

\begin{theorem}\label{partition_of_separated_set}
 Let

 $\bullet$ the subset $A\subset\R^{d}$  is  strongly independent with affine subspaces $L_1, ..., L_m\subset \R^d$;

 $\bullet$ $|A|+\sum\limits_{i=1}^{m}(1+\dim L_i)= (m-1)(d+1)+1$;

%$\bullet$ $k$ is integer number such that $0\leq k\leq |A|-|\{i|L_i=\varnothing\}|$;

 $\bullet$ $H$ is hyperplane in $R^d$ which divide it into subspaces $H^+$, $H^-$;

 $\bullet$   $A\cap H=\varnothing$.

Then there are a partitions $A_1, A_2, ..., A_m$ of $A\cap H^+$ and $B_1, B_2, ..., B_m$ of $A\cap H^-$ such that at least one of the interestions $\bigcap\limits_{i=1}^mC(L_i,A_i,B_i)$ and $\bigcap\limits_{i=1}^mC(L_i,B_i,A_i)$ is non-empty.
 \end{theorem}

 {\it Proof.}

Note that if $\bigcap\limits_{i=1}^mL_i\neq\varnothing$ then statement of Theorem~\ref{main_result} holds.
Next, we suppose that  $\bigcap\limits_{i=1}^mL_i=\varnothing$.
Without loss of generality $|A\cap H^-|\leq|A\cap H^+|$.

Denote $k=|A\cap H^-|$.
Then $k\leq\min(|A|-|\{i|L_i=\varnothing\}|;(m-1)d-\sum\limits_{i=1}^m\dim L_i)$.

Consider $\R^d$ as $\R^d\times0\subset\R^{d+1}$.
Apply Lemma~\ref{taking_higher_points} for the subset $\mathfrak{U}=A\cap H^-\subset A$ and $|\mathfrak{U}|=k$, we get:

 $\bullet$ a finite subset $A_{higher}\subset \R^d\times(1;+\infty)$;

 $\bullet$ a partition $\mathfrak{U}_1,...,\mathfrak{U}_m$ of $\mathfrak{U}$;

 $\bullet$ a pairwise disjoint subsets $\mathfrak{U}_1',...,\mathfrak{U}_m'$ of $A_{higher}\cup(A\cap H^+)$;

 $\bullet$ $|(\R^d\times(0;1))\cap\bigcap\limits_{i=1}^mC(L_i,\mathfrak{U}_i', \mathfrak{U}_i)|=1$.

There is an affine transformation which:

 $\bullet$ save the last coordinate in $\R^{d+1}$;

 $\bullet$ fix all points in the hyperplane $\R^d\times0$;

 $\bullet$ image of the set $A_{higher}$ belongs to the half-space $H^+\times\R\subset \R^{d+1}$.

 Apply this affine transformation to the set $A_{higher}$
 denote image by $A_{higher}'$.
 According to Lemma~\ref{sufficiently_general_position_lemma1} there are a points $W_1^{**}, ..., W_k^{**}\in H^-\times(1;+\infty)$ such that the set $A\cup A_{higher}'\cup\{W_1^{**}, ..., W_k^{**}\}$ in strong general position with the subspaces $L_1, ..., L_m$.

 Prove that there is a sufficiently divisible subset of $A\cup A_{higher}'\cup\{W_1^{**}, ..., W_k^{**}\}$ which contains exactly $k$ points from $H^-\times\R$.
 Denote by $f$ the last coordinate functional on $\R^{d+1}$.

  For $q=0,..., k$ take a pairwise disjoint subsets $Y_1^{(q)+}, ..., Y_m^{(q)+}, Y_1^{(q)-}, ..., Y_m^{(q)-}$ of $A\cup A_{higher}'\cup\{W_1^{**}, ..., W_k^{**}\}$, if it exists, such that:

$\bullet$ $Y_j^{(q)-}\subset (A\cap H^-)\cup\{W_1^{**}, ..., W_k^{**}\}$ for $j=1, ..., m$;

$\bullet$ $\sum\limits_{j=1}^{m}|Y_j^{(q)-}|=q$;

$\bullet$ $Y_j^{(q)+}\subset (A\cap H^+)\cup A_{higher}'$ for $j=1, ..., m$;

$\bullet$ $\bigcap\limits_{j=1}^{m}C(L_j, Y_j^{(q)+}, Y_j^{(q)-})$ is an unique point set; denote this point by $Z^{(q)}$;

$\bullet$ $f(Z^{(q)})>0$;

$\bullet$ $f(Z^{(q)})$ is minimal with respect to previous conditions.

Clearly, the point $Z^{(k)}$ exists and $Z^{(k)}\in\R^d\times(0;1)$.

Since $f(Z^{(k)})>0$, we have $C(L_i, Y_i^{(k)+}, Y_i^{(k)-})\not\subset\R^d\times0$ for $i=1,...,m$.
From sufficiently general position follows that $$\sum\limits_{j=1}^{m}(|Y_j^{(k)+}|+|Y_j^{(k)-}|)=(m-1)(d+1)-\sum\limits_{j=1}^{m}\dim L_j.$$
Hence $$\sum\limits_{j=1}^{m}(|Y_j^{(k)+}\cap A|+|Y_j^{(k)-}\cap A|)\leq (m-1)(d+1)-\sum\limits_{j=1}^{m}(1+\dim L_j).$$
Therefore there is a point $W\in A\setminus\{Y_1^{(k)+}\cup ...\cup Y_m^{(k)+}\cup Y_1^{(k)-}\cup ...\cup Y_m^{(k)-}\}$.

Define a subsets $Y_1^{+*}, ...,Y_m^{+*}, Y_1^{-*}, ..., Y_m^{-*}\subset A\cup A_{higher}'\cup\{W_1^{**}, ..., W_k^{**}\}$ by follows:

In the case, where there is no point $Z^{(k-1)}$ in $\R^d\times(0;f(Z^{(k)}))$:

According Lemma~\ref{linear_optimization}
there is a number $j_c$ such that the point $Z^{(k)}$ is not a minimal point of $f$ on the intersection $$C(L_{j_c}, \{W\}\cup Y_{j_c}^{(k)+}, Y_{j_c}^{(k)-})\cap\bigcap\limits_{j\neq j_c}C(L_j, Y_j^{(k)+}, Y_j^{(k)-}).$$

In this case take $Y_{j_c}^{+*}=\{W\}\cup Y_{j_c}^{(k)+}$ , $Y_{j_c}^{-*}=Y_{j_c}^{(k)-}$ and $$Y_{j}^{+*}=Y_{j}^{(k)+}, Y_{j}^{-*}=Y_{j}^{(k)-}$$
for $j\neq j_c$.
$\min\limits_{\bigcap\limits_{i=1}^mC(L_i,Y_i^{+*},,Y_i^{-*})}f<f(Z^{(k)})$.

In the case, where there is a point $Z^{(k-1)}$ in $\R^d\times(0;f(Z^{(k)}))$:

Since $\sum\limits_{i=1}^m|Y_i^{(k-1)-}|=k-1$, there is a point $W_2\in \{W_1^{**}, ..., W_k^{**}\}\setminus(\bigcup\limits_{i=1}^mY_i^{(k-1)-})$.
According Lemma~\ref{linear_optimization}
there is a number $j_c$ such that the point $Z^{(k-1)}$ is not a minimal point of $f$ on the intersection $$C(L_{j_c}, Y_{j_c}^{(k-1)+}, \{W_2\}\cup Y_{j_c}^{(k-1)-})\cap\bigcap\limits_{j\neq j_c}C(L_j, Y_j^{(k-1)+}, Y_j^{(k-1)-}).$$
In this case take $Y_{j_c}^{+*}=Y_{j_c}^{(k-1)+}$ , $Y_{j_c}^{-*}=\{W\}\cup Y_{j_c}^{(k-1)-}$ and $$Y_{j}^{+*}=Y_{j}^{(k-1)+}, Y_{j}^{-*}=Y_{j}^{(k-1)-}$$
for $j\neq j_c$.
Then $\min\limits_{\bigcap\limits_{i=1}^mC(L_i,Y_i^{+*},Y_i^{-*})}f<f(Z^{(k-1)})<f(Z^{(k)})$.

In both these cases we have $\min\limits_{\bigcap\limits_{i=1}^mC(L_i,Y_i^{+*},Y_i^{-*})}f<\min(f(Z^{(k-1)}),f(Z^{(k)}))$ and
$|(H^-\times\R)\cap\bigcup\limits_{i=1}^mY_i^{-*}|=k$.

If the polyhedron $\bigcap\limits_{i=1}^mC(L_i,Y_i^{+*},Y_i^{-*})$
does not interest the hyperplane $\R^d\times0$ then $f>0$ on this polyhedron.
Take a minimal point $Z^{*}$ of the functional $f$ on this polyhedron.
Then $f(Z^{*})>0$ and $Z^{*}$ is border point of this polyhedron.
From Lemma~\ref{polyhedron_border} follows that $\{Z^{*}\}=\bigcap\limits_{j=1}^mC(L_j, Y_j^{+**}, Y_j^{-**})$ for some subsets $$Y_j^{+**}\subset Y_j^{+*}, Y_j^{-**}\subset Y_j^{-*}$$ for $j=1,...,m$.
From sufficiently general position follows that $$|\bigcup\limits_{j=1}^m(Y_j^{+**}\cup Y_j^{-**})|=-1+|\bigcup\limits_{j=1}^m(Y_j^{+*}\cup Y_j^{-*})|.$$
Hence $\sum\limits_{i=1}^m|Y_j^{-**}|\in\{k-1,k\}$.
A contradiction with the choice of points $Z^{(k)}, Z^{(k-1)}$  because $f(Z^*)=\min\limits_{\bigcap\limits_{i=1}^mC(L_i,Y_i^{+*},Y_i^{-*})}f<\min(f(Z^{(k-1)}),f(Z^{(k)}))$.

Therefore $$\bigcap\limits_{i=1}^mC(L_i,Y_i^{+*},Y_i^{-*})\cap(\R^d\times0)\neq\varnothing.$$
Then the set $\bigcup\limits_{j=1}^m(Y_j^{+*}\cup Y_j^{-*})$ is sufficiently divisible.
The existence a sufficiently divisible subset of $A\cup A_{higher}'\cup\{W_1^{**}, ..., W_k^{**}\}$ which contains exactly $k$ points from $H^-\times\R$ is proved.

Take a  subset of $A''\subset A\cup A_{higher}'\cup\{W_1^{**}, ..., W_k^{**}\}$ such that:

$\bullet$ $A''$ is sufficiently divisible;

$\bullet$ $|A''\cap H^-|=k$;

$\bullet$ $|A''\cap (\R^d\times(0;+\infty))|$ is minimal with respect of previous conditions.

Suppose that $|A''\cap (\R^d\times(0;+\infty))|\geq m+1$.
Then there are a points:

$\bullet$ $X_1\in A\setminus A''$;

$\bullet$ $X_0\in A''\cap (\R^d\times(0;+\infty))$ which does not separated from the point $X_1$ by the hyperplane $H\times \R$.
Denote by $T$ a subset of numbers $t\in[0;1]$ such that the set $(A''\setminus\{X_0\})\cup \{tX_1+(1-t)X_0\}$ is sufficiently divisible.
From the choice of $A''$ follows that $0\in T$.
Hence $T\neq\varnothing$.
From Lemma~\ref{moving_the_point} follows that the subset $T\subset[0;1]$ is open.
From both Lemmas~\ref{moving_the_point} and~\ref{moving_the_point_2} follows that the number $\sup T$ has a punctured neighbourhood which contains in $T$.
Then $\sup T=1$.
From Lemma~\ref{moving_the_point_2} follows that the number 1 is not border point of $T$.
Hence $1\in T$.
Therefore the set $(A''\setminus\{X_0\})\cup \{X_1\}$ is sufficiently divisible.
A contradiction with choice of $A''$ because $$|((A''\setminus\{X_0\})\cup\{X_1\})\cap (\R^d\times(0;+\infty))|<|A''\cap (\R^d\times(0;+\infty))|.$$

Then $|A''\cap (\R^d\times(0;+\infty))|\leq m$.
Take a partition $A_1'', ..., A_m''$ of $A''$ such that:

$\bullet$ $|A_i''|\leq d+1-\dim L_i$ for any $i=1,...,m$;

$\bullet$ either of polyhedra $$\bigcap\limits_{i=1}^mC(L_i\cap H, A_i''\cap (H^+\times\R), A_i''\cap (H^-\times\R))$$ or $$\bigcap\limits_{i=1}^mC(L_i\cap H, A_i''\cap (H^-\times\R), A_i''\cap (H^+\times\R))$$ does not contains in the hyperplane $\R^d\times0$ but intersect it. %Denote this polyhedron by $P_1$.

Hence $A_i''\not\subset\R^d\times0$ for $i=1,...,m$.
Therefore $|A_i''\cap(\R^d\times(0;+\infty))|=1$ for $i=1,...,m$.
Then both:
$$(\R^d\times0)\cap C(L_i\cap H, A_i''\cap (H^+\times\R), A_i''\cap (H^-\times\R))=C(L_i\cap H, A_i''\cap H^+, A_i''\cap H^-)$$ and $$(\R^d\times0)\cap C(L_i\cap H, A_i''\cap (H^-\times\R), A_i''\cap (H^+\times\R))= C(L_i\cap H, A_i''\cap H^-, A_i''\cap H^+)$$
for $i=1,...,m$.
Take $A_i=A_i''\cap H^+$, $B_i=B_i''\cap H^-$.
Then at least one of the interestions $\bigcap\limits_{i=1}^mC(L_i,A_i,B_i)$ and $\bigcap\limits_{i=1}^mC(L_i,B_i,A_i)$ is non-empty.

 {\it Proof of Theorem~\ref{main_result_2} in conditions 1,3.}
 Take a  hyperplane $H\subset\R^d$ which separate $k$ points of the set $A$ from other points.
 Apply Theorem~\ref{partition_of_separated_set}.
 Take a partition $A_1, A_2, ..., A_m, B_1, B_2, ..., B_m$ of the set $A$.
 %Without loss of generality $\sum\limits_{i=1}^m|B_i|=k$.

 If $k+|\{i|L_i\neq\varnothing\}|\leq m-1$, then there is a number $i'$ such that $B_{i'}=L_{i'}=\varnothing$.
 Hence the polyhedron $C(L_{i'},B_{i'},A_{i'})$ does not exists.
 Therefore $\bigcap\limits_{i=1}^mC(L_i,A_i,B_i)\neq\varnothing$.
 Theorem~\ref{main_result_2} proved in the condition 1.

 If $k=\frac{1}{2}|A|$, then $\sum\limits_{i=1}^m|A_i|=|A|-k=k=\sum\limits_{i=1}^m|B_i|$. Take non-empty interestion from $\bigcap\limits_{i=1}^mC(L_i,A_i,B_i)$ and $\bigcap\limits_{i=1}^mC(L_i,B_i,A_i)$.
 Theorem~\ref{main_result_2} proved in the condition 3.

 {\it Proof of Theorem~\ref{main_result_2} in condition 2.}

 We proove Theorem~\ref{main_result_2} in the condition 2 by induction on $k$.

 Base of induction: $k\in\{\lfloor\frac{|A|}{2}\rfloor, 1+\lfloor\frac{|A|}{2}\rfloor\}$.

 Take a hyperplane $H\subset\R^d$ which separate some $\lfloor\frac{|A|}{2}\rfloor$ point of the set $A$ from other $|A|-\lfloor\frac{|A|}{2}\rfloor$ points.
Apply Theorem~\ref{partition_of_separated_set}.
 Take a partition $A_1, A_2, ..., A_m, B_1, B_2, ..., B_m$ of the set $A$.
 Without loss of generality $\sum\limits_{i=1}^m|B_i|=\lfloor\frac{|A|}{2}\rfloor$.
 Then $\sum\limits_{i=1}^m|A_i|=|A|-\lfloor\frac{|A|}{2}\rfloor\in\{\lfloor\frac{|A|}{2}\rfloor, 1+\lfloor\frac{|A|}{2}\rfloor\}$.
Take non-empty interestion from $\bigcap\limits_{i=1}^mC(L_i,A_i,B_i)$ and $\bigcap\limits_{i=1}^mC(L_i,B_i,A_i)$.
Theorem~\ref{main_result_2} in the condition 2 proved for some $k\in\{\lfloor\frac{|A|}{2}\rfloor, 1+\lfloor\frac{|A|}{2}\rfloor\}$.

Step of induction: from $k$ to $k-1$.
According to hypothesis there is partition $A_1, A_2, B_1, B_2$ of $A$ such that:

$\bullet$ $C(L_1,A_1,B_1)\cap C(L_2,A_2,B_2)\neq \varnothing$;

$\bullet$ $|B_1|+|B_2|=k$.

Without loss of generality $B_2\neq\varnothing$.
Take a point $X\in B_2$.
Since $$C(L_1,A_1,B_1)\cap C(\varnothing,C(L_2,A_2,B_2\setminus\{X\}),\{X\})\neq \varnothing,$$
the interestion $C(L_1,A_1\cup\{X\},B_1)\cap C(L_2,A_2,B_2\setminus\{X\})$ is not empty.

Step is proved.
Hence Theorem~\ref{main_result_2} in the condition 2 is proved.


\begin{thebibliography}{RSS95}



\bibitem[BS16]{BS16}  \emph{I. Barany, P. Soberon} (2018).  Tverberg Plus Minus. Discrete Comput Geom 60, 588-598. https://doi.org/10.1007/s00454-017-9960-1

\bibitem[Ro01]{Ro01} \emph{J-P. Roudneff} (2001) Partitions of points into simplices with k-dimensional intersection. I. The conic Tverberg's theorem. European J. Combin., 733-743.

\bibitem[Re68]{Re68} \emph{J. R. Reay}  An extension of Radon's theorem, Illinois J. Math., 12 (1968), 184-189.

\bibitem[Tv66]{Tv66} \emph{H. Tverberg} (1966), A Generalization of Radon's Theorem. Journal of the London Mathematical Society, s1-41: 123-128. doi:10.1112/jlms/s1-41.1.123

\end{thebibliography}
\end{document}